\documentclass[12pt,notitlepage]{article}

\usepackage[english]{babel}
\usepackage[utf8]{inputenc}
\usepackage{color}
\usepackage{colortbl}
\usepackage[dvipsnames]{xcolor}
\usepackage[a4paper,lmargin={2cm},rmargin={2cm},
tmargin={2.5cm},bmargin = {2.5cm}]{geometry}
\usepackage{amssymb}
\usepackage{amsthm}
\usepackage{amsmath}
\usepackage{mathtools}
\usepackage{graphicx}
\usepackage{hyperref}
\usepackage{pdfpages}
\usepackage{appendix}
\usepackage{listings}
\usepackage{tikz}
\usepackage{dsfont}
\usepackage{subcaption}
\usepackage{multicol}
\usepackage{caption}
\usepackage{float}
\usepackage{etoolbox}
\usepackage{authblk}
\usepackage{setspace}
\usepackage{titlesec}
\usepackage{changepage}
\usepackage{authblk}

\usepackage{algorithm}
\usepackage{algorithmic}

\usepackage[backend=biber, 
            style=numeric,   
            doi=true,        
            url=true,        
            hyperref=true    
]{biblatex}
\usepackage{csquotes}
\title{\LARGE Domain decomposition dynamical low-rank for multi-dimensional radiative transfer equations}
\author[1]{Stefan Brunner}
\author[1]{Lukas Einkemmer}
\author[2]{Terry Haut}
\affil[1]{Department of Mathematics, University of Innsbruck, Austria}
\affil[2]{Lawrence Livermore National Laboratory, USA}
\date{July 23, 2026}

\titleformat{\section}{\large\bfseries}{\thesection}{1em}{}

\newcommand\inner[3]{\left\langle #1, #2 \right\rangle_{#3}}
\newcommand\innersingle[2]{\left\langle #1 \right\rangle_{#2}}

\newcommand\suml[3]{\sum\limits_{#1}^{#2} #3}
\newcommand*\diff{\mathop{}\!\mathrm{d}}
\newcommand{\norm}[1]{\left\lVert#1\right\rVert}

\begin{document}

\maketitle

\begin{center}
    \textbf{Abstract}
\end{center}
  
\begin{adjustwidth}{1cm}{1cm} 
    \small 
    In this paper, we propose a domain decomposition dynamical low-rank method to solve high-dimensional radiative transfer problems and similar kinetic equations. The algorithm uses a separate low-rank approximation on each spatial subdomain, which means that, for a given accuracy, we can often use {\color{black}fewer degrees of freedom} compared to classic dynamical low-rank methods. In particular, we can solve problems with point sources efficiently, that for classic algorithms require almost full rank. Our algorithm only transfers boundary data between subdomains and is thus very attractive for distributed memory parallelization, where classic dynamical low-rank algorithms suffer from global data dependency. We demonstrate the efficiency of our algorithm by a number of challenging test examples that have both very optical thin and thick regions.
\end{adjustwidth}

\section{Introduction}\label{sec::introduction}

Accurately simulating radiation transport phenomena is vital in many fields. For example, in inertial confinement fusion, atmospheric science and astrophysics. Radiation transport can be described mathematically by the radiative transfer equation (RTE). It models the movement of particles through physical materials, taking absorption, emission and scattering processes into account. Since the distribution of the radiation intensity depends on phase space, i.e.~on both spatial and angular variables, we have a high-dimensional problem whose numerical solution suffers from the curse of dimensionality. Thus, finding an efficient solver is difficult.

Commonly, numerical methods based on a phase space discretization are used. This includes the discrete ordinates method (called $S_N$){\color{black}\cite{Haut2019,Mathews1999}}, which evaluates the distribution function pointwise, and methods based on variants of spherical harmonics in velocity space (called $P_N$){\color{black}\cite{McClarren2007,BRUNNER2005}}, i.e.~methods which expand the velocity distribution in a global basis. Both methods do suffer from the curse of dimensionality. That is, the number of degrees of freedom scales unfavorable for higher dimensional problems. Another approach, which has also been widely used, are implicit Monte Carlo methods{\color{black}\cite{Dupree2002}}. This approach is based on repeated random sampling of particle trajectories in order to obtain an estimate of the solution. While convergence is technically not dimension dependent, it is usually very slow (square root in the number of particles). Thus, usually this approach also comes with high computational cost, especially for problems with regions that are optically thick. For more details we refer the reader to, e.g., \cite{Brunner_testproblem}.

Recently, dynamical low-rank approximations (DLRA) have seen increasing interest for kinetic problems in general and radiation transport more specifically (see, e.g., the review article \cite{Einkemmer2024c}). The main idea is to approximate the distribution function using lower dimensional low-rank factors that only depend on either space or angular variables, thereby significantly reducing memory and computational cost. For these low-rank factors, (partial differential) equations of motions are then derived that can be solved numerically. Some early work on the radiative transfer equation can be found in \cite{Peng2020,Ding2021,Einkemmer2021} and applications such as radiation therapy \cite{Kusch2021a} have been pursued. In addition, much work in recent years has resulted in conservative variants of these schemes (see, e.g., \cite{Einkemmer2021a,Coughlin2023,Guo2022b,Peng2021}) and (implicit or semi-implicit) schemes that preserve the underlying diffusive or fluid limit of the kinetic equation (see, e.g., \cite{Einkemmer2021,Einkemmer2022a,Patwardhan2024}).

Dynamical low-rank methods can be extremely effective for such problems and often result in drastically reduced memory and computational cost. This is particularly true on modern computing architectures, such as graphic processing units (GPUs; see \cite{Cassini2021}). However, the dynamical low-rank approach fundamentally uses an approximation with global basis functions. This results in global data dependency, which is undesirable when parallelizing such methods to distributed memory supercomputers. Moreover, there are some problems where there is no global low-rank structure (e.g.~a point source). In this case classic low-rank methods do require a large rank to obtain accurate results.

In this paper, we propose a domain decomposition dynamical low-rank approach that uses a separate low-rank approximation on each subdomain. In particular, 
\begin{itemize}
  \item the method only exchanges boundary data between subdomains. That is, there is only local data dependency between subdomains;
  \item since we have a low-rank approximation on each subdomain, we can solve problems that are of high rank locally and problems, which have parts of the domain where the rank is low (e.g.~due to strong absorption), much more efficiently than classic dynamical low-rank methods.
\end{itemize}
As we will see, a main difficulty is to obtain a domain decomposition method that does not suffer from excessive increase of rank in intermediate steps (which would increase memory cost and decrease computational efficiency significantly). Our approach, which is based on the projector splitting integrator, requires some augmentation but not more so than classic Basis Update and Galerkin (BUG) integrators (see, e.g., \cite{BUG_Gianluca_augmented}).

As related work, we mention \cite{McClarren_testproblem}, where the authors decompose the angular domain into multiple parts and apply DLRA to each of those parts separately. This is an easier problem, because if the subdomains are chosen as quadrants/octants of angular space, no coupling between the subdomains is necessary. However, the approach is also more limiting since usually much more grid points are used in the spatial domain and there is thus more potential for parallelization. Moreover, different parts in the spatial domain can behave very differently. For example, in inertial confinement fusion it is very common to have parts of the domain which are very optically thick and are thus close to the diffusion limit. Such regions can be represented easily with a very low rank. Our spatial domain decomposition approach can easily exploit such a structure of the problem. The approach of \cite{McClarren_testproblem}, on the other hand, would still require a low-rank approximation that is global in the spatial variables. {\color{black} A potential limitation of the current implementation for our spatial domain decomposition approach are load-balancing issues caused by varying numerical ranks across subdomains.}

We also mention \cite{adapt_DLR_nonlin_boltz}, which proposes a low-rank approximation that treats inflow/outflow boundary conditions. While the authors do not directly consider domain decomposition in their work, the problem of imposing inflow/outflow boundary conditions is very much related to our approach, which only communicates boundary information between two subdomains.

The remainder of the paper is structured as follows. In section \ref{sec::dlra}, we introduce our model of the radiative transfer equation and give a short recap of how classic dynamical low-rank methods are used in that context. In section \ref{sec::algorithm}, we present the proposed domain decomposition algorithm. Finally, in section \ref{sec::num_results}, we show numerical results for some known test problems, comparing our new algorithm with the classic DLRA approach.

\section{Dynamical low-rank approximation of the radiative transfer equation}\label{sec::dlra}

In this paper we consider the following $2\times1$ dimensional radiative transfer equation (see, e.g.,~\cite{Brunner_testproblem})
\begin{equation}\label{eq::rt}
\begin{aligned}
  \partial_t f(t,x,y,\phi) &+ c_{\text{adv}}v \cdot \nabla_{x,y}f(t,x,y,\phi) \\
  &= c_{\text{s}}(x,y)\frac{1}{2\pi}\rho(t,x,y) - c_{\text{t}}(x,y)f(t,x,y,\phi) + Q(x,y),
\end{aligned}
\end{equation}
where $(x,y) \in \Omega = [0,1]\times[0,1]$, the angle $\phi \in [0,2\pi)$ and $f: \mathbb{R}^+ \times \Omega \times [0,2\pi) \rightarrow \mathbb{R}^+$ is the distribution function. It describes the radiation energy
density flowing in a particular direction specified by the angle $\phi$. The velocity $v$ is given by $v = (\cos \phi, \sin \phi)$. 
The coefficient $c_{\text{adv}}$ is the advection speed, which is equal to the value of the speed of light in the chosen unit system. The other coefficients, $c_{\text{s}}$ and $c_{\text{t}}$, 
are the scattering and the total coefficients. The total coefficient is defined as the sum of the absorption and scattering coefficient, i.e. $c_{\text{t}} = c_{\text{a}} + c_{\text{s}}$. Lastly, $Q$ is an external source term.

The most important physical observable is the density $\rho$, which can be computed from the distribution function as follows
\[ \rho(t,x,y) = \int f(t,x,y,\phi) \diff \phi = \innersingle{f}{\phi}. \]

Due to the high dimensionality of the problem, we want to approximate the distribution function $f(t,x,y,\phi)$ by a low-rank approximation
\begin{equation}
  f(t,x,y,\phi) = \suml{i,j=1}{r}{U_i(t,x,y)S_{ij}(t)V_j(t,\phi)}.
\end{equation}
Note that this is similar to a singular value decomposition (SVD), but formulated here in the continuous setting. The functions $U_i$, $S_{ij}$ and $V_j$ are the low-rank factors and satisfy the orthonormality conditions $\inner{U_i}{U_k}{x,y}=\delta_{ik}$ and $\inner{V_j}{V_l}{\phi} = \delta_{jl}$, where $\inner{\cdot}{\cdot}{}$ denotes the standard $L^2$ inner product. From now on we will denote sums starting at index $1$ and ending at index $r$ simply by $\suml{i,j}{}{}$, meaning we only write the summation indices but not the bounds.

Following \cite{Koch2007,PSI_Lukas}, a system of evolution equations for the low-rank factors can be derived. The resulting system, however, is ill-conditioned in the presence of small singular values. To avoid this problem so-called robust integrators have been developed. The most common ones are the projector splitting integrator (PSI) \cite{Lubich2014,PSI_Lukas} and the basis update and Galerkin integrator (BUG) \cite{BUG_Gianluca,BUG_Gianluca_augmented}, but also a parallel integrator \cite{Parallel_Gianluca_testproblem} has been proposed.

In this work we focus, for reasons that will be explained in more detail later, on the projector splitting integrator. The PSI splits the projector onto the tangent space of the low-rank manifold  into three parts corresponding to the three low-rank factors. We will now summarize the method (for more details, see \cite{BUG_Gianluca, PSI_Lukas, Lubich2014}). First, we write equation \eqref{eq::rt} as follows
\begin{equation}
    \partial_t f  = F(f), \qquad
    F(f)= -c_{\text{adv}} v \cdot \nabla_{x,y}f + c_{\text{s}} \frac{1}{2\pi}\rho- c_{\text{t}}f + Q. \label{eq:rhs}
\end{equation}
Starting with a low-rank representation at time $t_n$
\[ f^n(x,y,\phi) = \suml{i,j}{}{U_i^n(x,y)S_{ij}^n V_j^n(\phi)}, \]
the (first order) PSI proceeds as follows:
\begin{enumerate}
    \item \textbf{K-step:} Integrate from $t_n$ to $t_{n+1}$ the differential equation
    \begin{equation}\label{eq::Kstep_general}
    \begin{aligned}
        &\partial_t K_j (t,x,y) = \inner{V_j^n(\phi)}{F\left(\suml{k}{}{K_k(t,x,y)V_k^n(\phi)} \right)}{\phi}, \\
        & K_j(t_n,x,y) = \suml{i}{}{U_i^n(x,y)S_{ij}^n},
    \end{aligned}
    \end{equation}
        to obtain $K_j^{n+1}(x,y) = K_j(t_{n+1},x,y)$. For the space spanned by the functions $K_j^{n+1}(x,y)$ for $j\in \{ 1, \dots , r\}$, compute an orthonormal basis (e.g. by a QR decomposition) $U_i^{n+1}(x,y)$ for $i\in \{ 1, \dots , r\}$, such that $K_j^{n+1}(x,y) =\suml{i}{}{U_i^{n+1}(x,y) \hat{S}_{ij}}$.

    \item \textbf{S-step:} Integrate from $t_n$ to $t_{n+1}$ the differential equation
    \begin{equation}\label{eq::Sstep_general}
    \begin{aligned}
        &\partial_t S_{ij}(t) = -\inner{U_i^{n+1}(x,y)V_j^n(\phi)}{F\left(\suml{k,l}{}{U_l^{n+1}(x,y)S_{lk}(t)V_k^n(\phi)}\right)}{x,y,\phi}, \\
        &S_{ij}(t_n) = \hat{S}_{ij}.
    \end{aligned}
    \end{equation}
    
        Set $\Tilde{S}_{ij}=S_{ij}(t_{n+1})$.

    \item \textbf{L-step:} Integrate from $t_n$ to $t_{n+1}$ the differential equation
    \begin{equation}\label{eq::Lstep_general}
    \begin{aligned}
        &\partial_t L_i(t,\phi) = \inner{U_i^{n+1}(x,y)}{F\left(\suml{l}{}{U_l^{n+1}(x,y)L_l(t,\phi)} \right)}{x,y}, \\
        &L_i(t_n,\phi) = \suml{j}{}{V_j^n(\phi) \Tilde{S}_{ij}},
    \end{aligned}
    \end{equation}
        to obtain $L^{n+1}_i(\phi) = L(t_{n+1},\phi)$. For the space spanned by the functions $L_i^{n+1}(\phi)$ for $i\in \{ 1, \dots , r\}$, compute an orthonormal basis $V_j^{n+1}(\phi)$ for $j\in \{ 1, \dots , r\}$, such that $L_i^{n+1}(\phi)=\suml{j}{}{V_j^{n+1}(\phi) S_{ij}^{n+1}}$.
\end{enumerate}

The approximation at time $t_{n+1}$ is then given by 
\begin{equation*}
f^{n+1}(x,y,\phi) = \suml{i,j}{}{U_i^{n+1}(x,y)S_{ij}^{n+1}V_j^{n+1}(\phi)}.
\end{equation*}

Next, apply the projector splitting integrator to the radiative transfer equation \eqref{eq::rt}. For the sake of simplicity, we will omit variable dependency. Plugging equation \eqref{eq:rhs} into equation \eqref{eq::Kstep_general} we get the {\color{black}K-step}
\begin{equation}\label{eq::Kstep}
  \partial_t K_j =-c_{\text{adv}}\suml{k}{}{\inner{V_j^n}{vV_k^n}{\phi}\nabla_{x,y}K_k} + c_{\text{s}}\frac{1}{2\pi}\suml{k}{}{K_k\innersingle{V_k^n}{\phi}\innersingle{V_j^n}{\phi}} - c_{\text{t}}K_j + Q \innersingle{V_j^n}{\phi}.
\end{equation}
Similarly, we get the evolution equation for the S-step
\begin{equation}\label{eq::Sstep}
\begin{aligned}
  \partial_tS_{ij} &= c_{\text{adv}}\suml{k,l}{}{\inner{U_i^{n+1}}{\nabla_{x,y}U_l^{n+1}}{x,y}\inner{V_j^n}{vV_k^n}{\phi}S_{lk}} \\
  &- \frac{1}{2\pi}\suml{k,l}{}{\inner{U_i^{n+1}}{c_{\text{s}}U_l^{n+1}}{x,y}\innersingle{V_k^n}{\phi}\innersingle{V_j^n}{\phi}S_{lk}} \\
  &+ \suml{l}{}{\inner{U_i^{n+1}}{c_{\text{t}}U_l^{n+1}}{x,y}S_{lj}}-\inner{U_i^{n+1}}{Q}{x,y}\innersingle{V_j^n}{\phi}
\end{aligned}
\end{equation}
and for the L-step
\begin{equation}\label{eq::Lstep}
\begin{aligned}
  \partial_tL_i&= -c_{\text{adv}}\suml{l}{}{\inner{U_i^{n+1}}{\nabla_{x,y}U_l^{n+1}}{x,y}vL_l} + \frac{1}{2\pi}\suml{l}{}{\inner{U_i^{n+1}}{c_{\text{s}}U_l^{n+1}}{x,y}\innersingle{L_l}{\phi}} \\
  &- \suml{l}{}{\inner{U_i^{n+1}}{c_{\text{t}}U_l^{n+1}}{x,y}L_l} + \inner{U_i^{n+1}}{Q}{x,y}.
\end{aligned}
\end{equation}
Note that as $U$ and $K$ depend only on spatial variables, $V$ and $L$ only depend on the angle, and $S$ does not depend on either, storing the low-rank factors is much cheaper than storing the entire distribution function. Moreover, equations \eqref{eq::Kstep}-\eqref{eq::Lstep} can be evaluated efficiently and the computational complexity of doing so scales at most with either the number of spatial grid points or the number of angles (but not the product of both). For a more detailed computational and memory complexity analysis we refer the reader, e.g., to \cite{PSI_Lukas}.

\section{Domain decomposition low-rank algorithm}\label{sec::algorithm}

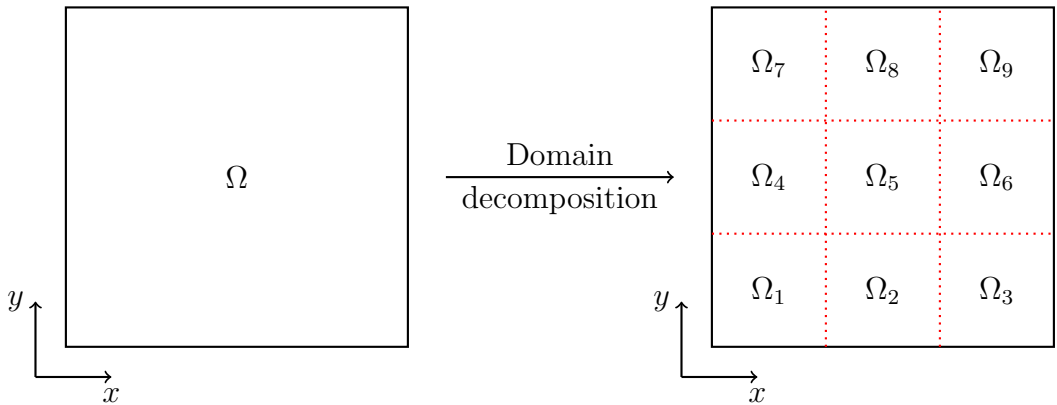
\begin{figure}[H]
  \centering
  \begin{tikzpicture}
    \draw[thick] (0,0) rectangle (4.5,4.5);
    
    \node at (2.25,2.25) {$\Omega$};

    \draw[->, thick] (-0.4,-0.4) -- (0.6,-0.4) node[below] {$x$};
    \draw[->, thick] (-0.4,-0.4) -- (-0.4,0.6) node[left] {$y$};

    \draw[->, thick] (5,2.25) -- (8,2.25);
    \node at (6.5,2.55) {Domain};
    \node at (6.5,1.95) {decomposition};

    \begin{scope}[shift={(8.5,0)}] 
      \draw[thick] (0,0) rectangle (4.5,4.5);
      
      \draw[red, thick, dotted] (1.5,0) -- (1.5,4.5);
      \draw[red, thick, dotted] (3,0) -- (3,4.5);
      
      \draw[red, thick, dotted] (0,1.5) -- (4.5,1.5);
      \draw[red, thick, dotted] (0,3) -- (4.5,3);
      
      \node at (0.75,0.75) {$\Omega_{1}$};
      \node at (2.25,0.75) {$\Omega_{2}$};
      \node at (3.75,0.75) {$\Omega_{3}$};
      
      \node at (0.75,2.25) {$\Omega_{4}$};
      \node at (2.25,2.25) {$\Omega_{5}$};
      \node at (3.75,2.25) {$\Omega_{6}$};
      
      \node at (0.75,3.75) {$\Omega_{7}$};
      \node at (2.25,3.75) {$\Omega_{8}$};
      \node at (3.75,3.75) {$\Omega_{9}$};
      
      \draw[->, thick] (-0.4,-0.4) -- (0.6,-0.4) node[below] {$x$};
      \draw[->, thick] (-0.4,-0.4) -- (-0.4,0.6) node[left] {$y$};
    \end{scope}
  \end{tikzpicture}
  \caption{Illustration of block domain decomposition.}
  \label{fig::block_domain_decomposition}
\end{figure}

In this section, we describe the proposed algorithm, which combines domain decomposition with the dynamical low-rank approach to solve equation \eqref{eq::rt}. A block based domain decomposition, as illustrated in Figure \ref{fig::block_domain_decomposition}, is considered. Our approach is based on imposing inflow/outflow conditions on each subdomain. That is, the outflow of a neighboring domain becomes the inflow for the current domain and vice versa. {\color{black}Such boundary conditions have already been considered in \cite{adapt_DLR_nonlin_boltz}}. An issue here is that potentially new information enters via the boundary as the neighboring subdomain uses a low-rank approximation with, in general, different basis functions. That is, the low-rank approximation does not necessarily have the correct basis to accurately represent the information from the boundary. This was already recognized in \cite{adapt_DLR_nonlin_boltz} and, in fact, a naive algorithm gives incorrect results. The obvious solution is to augment the basis with information from the neighboring domains. 

However, doing this in a naive way would be very costly as in two dimensions we have potentially eight neighboring subdomains which all, in general, have a distinct low-rank representation. This would increase memory usage drastically and is even worse in terms of computational effort as the computational cost scales with $r^2$. We solve this problem in the following way:
\begin{itemize}
  \item First, we perform a splitting of the advection such that only one-dimensional advections need to be considered. This has the benefit that we only require augmentation by at most two neighboring subdomains in each step of the algorithm. Splitting into different directions is a very common approach for kinetic equations (see, e.g., \cite{Cheng1976,Sonnendruecker1999,Einkemmer2016a}), although in these works it is done for a somewhat different reason (to avoid the complexity of multi-dimensional interpolation).
  \item Since there is usually a significant overlap of the information contained in the basis functions of neighboring domains, we add only information that is not redundant. That is, we only add information that is not yet contained in the low-rank representation of the current subdomain.
\end{itemize}

While, in principle, doing such an augmentation is very natural for the BUG integrators (see, e.g., \cite{BUG_Gianluca_augmented}), doing so would result in a relatively large number of basis functions, {\color{black}i.e., a rank of $3r$ or even more}. The reason for this is that modern BUG integrators have an intrinsic augmentation that needs to be done in addition to the augmentation required by the inflow. The larger the system size, the more problematic this becomes due to the at least quadratic scaling in rank that holds for all of these methods. Thus, we will focus here on the projector splitting integrator. This allows us to obtain an algorithm (as we will see in the numerical results) that does not need to increase the rank for the intermediate computation to above roughly $2r$ (i.e., ~is about the same as for a BUG integrator {\color{black}without additional rank augmentation due to the inflow conditions}).

In the following we will explain our algorithm in detail for the two-dimensional case. However, due to the splitting the extension to three dimensions is immediate. We start with the radiative transfer equation \eqref{eq::rt} and treat the $x$ advection, the $y$ advection, and the combined absorption, scattering and the source term separately. Using $v_x = \cos{(\phi)}$ and $v_y = \sin{(\phi)}$ the splitting algorithm to update $f^n(x,y,\phi)$ to $f^{n+1}(x,y,\phi)$ reads as follows.
\begin{enumerate}
  \item We first consider the $x$ advection:
      \begin{align*}
    &\partial_t f(t,x,y,\phi) + c_{\text{adv}}v_x(\phi)\partial_xf(t,x,y,\phi) = 0, \\
    &f(t_n,x,y,\phi) = f^n(x,y,\phi).
  \end{align*}
  By solving this equation we obtain $f^{*}(x,y,\phi) = f(t_{n+1},x,y,\phi)$.
  \item Next, we treat the $y$ advection:
  \begin{align*}
    &\partial_t f(t,x,y,\phi) + c_{\text{adv}}v_y(\phi)\partial_yf(t,x,y,\phi) = 0, \\
    &f(t_n,x,y,\phi) = f^{*}(x,y,\phi).
  \end{align*}
  By solving the equation we obtain $f^{**}(x,y,\phi) = f(t_{n+1},x,y,\phi)$.
  \item Finally, we consider absorption, scattering and the source term:
  \begin{align*}
    &\partial_t f(t,x,y,\phi) = c_{\text{s}}(x,y)\frac{1}{2\pi}\rho(t,x,y) - c_{\text{t}}(x,y)f(t,x,y,\phi) + Q(x,y), \\
    &f(t_n,x,y,\phi) = f^{**}(x,y,\phi).
  \end{align*}
  By solving the equation we obtain $f^{n+1}(x,y,\phi) = f(t_{n+1},x,y,\phi)$.
\end{enumerate}
This is first order accurate in time. However, it can easily be raised to second or higher order by the usual splitting techniques (see, e.g., \cite{Lubich2014,Hairer2006}). In the following we apply the projector splitting integrator to each of these three steps separately and discuss how to impose the inflow boundary conditions. 

\subsection{Step 1: DLR advection in x}

We start with the advection in the $x$-direction, i.e. step 1. Similar to section \ref{sec::dlra} we get
\begin{align}
  \partial_t K_j &=-c_{\text{adv}}\suml{k}{}{\inner{V_j}{v_x V_k}{\phi} \partial_x K_k}, \label{advx-evolK-orig} \\
  \partial_t S_{ij} &= c_{\text{adv}} \suml{k,l}{}{\inner{U_i}{ \partial_ x U_l}{x,y}\inner{V_j}{v_x V_k}{\phi}S_{lk}},\label{advx-evolS-orig} \\
  \partial_t L_i &= -c_{\text{adv}}\suml{l}{}{\inner{U_i}{ \partial_x U_l}{x,y}v_x L_l}. \label{advx-evolL-orig}
\end{align}

\begin{figure}[H]
  \centering
  \begin{tikzpicture}
    \draw[thick] (0,0) rectangle (9,3);

    \draw[red, thick] (3,0) node[anchor=north] {$\partial\Omega_L$} -- (3,3) node[anchor=south] {$x=x_L$};
    \draw[red, thick] (6,0) node[anchor=north] {$\partial\Omega_R$} -- (6,3) node[anchor=south] {$x=x_R$};

    \node at (1.5,1.5) {$\Omega_{L}$};
    \node at (4.5,1.5) {$\Omega_{C}$};
    \node at (7.5,1.5) {$\Omega_{R}$};

    \foreach \x in {-0.5, -1, -1.5} {
      \filldraw (\x, 1.5) circle (1pt);
    }
    \foreach \x in {9.5, 10, 10.5} {
      \filldraw (\x, 1.5) circle (1pt);
    }
    \foreach \y in {3.5, 4, 4.5} {
      \filldraw (1.5, \y) circle (1pt);
      \filldraw (4.5, \y) circle (1pt);
      \filldraw (7.5, \y) circle (1pt);
    }
    \foreach \y in {-0.5, -1, -1.5} {
      \filldraw (1.5, \y) circle (1pt);
      \filldraw (4.5, \y) circle (1pt);
      \filldraw (7.5, \y) circle (1pt);
    }

    \draw[->, thick, blue] (3,1.5) -- (3.75,1.5);
    \draw[->, thick, blue] (3,2.25) -- (3.75,2.25);
    \draw[->, thick, blue] (3,0.75) -- (3.75,0.75);
    \draw[->, thick, blue] (6,1.5) -- (6.75,1.5);
    \draw[->, thick, blue] (6,2.25) -- (6.75,2.25);
    \draw[->, thick, blue] (6,0.75) -- (6.75,0.75);

    \draw[->, thick, Green] (3,1.5) -- (2.25,1.5);
    \draw[->, thick, Green] (3,2.25) -- (2.25,2.25);
    \draw[->, thick, Green] (3,0.75) -- (2.25,0.75);
    \draw[->, thick, Green] (6,1.5) -- (5.25,1.5);
    \draw[->, thick, Green] (6,2.25) -- (5.25,2.25);
    \draw[->, thick, Green] (6,0.75) -- (5.25,0.75);

  \end{tikzpicture}
  \caption{Inflow and outflow for subdomain $\Omega_C$ during the x advection step.}
  \label{fig::inflow_outflow}
\end{figure}
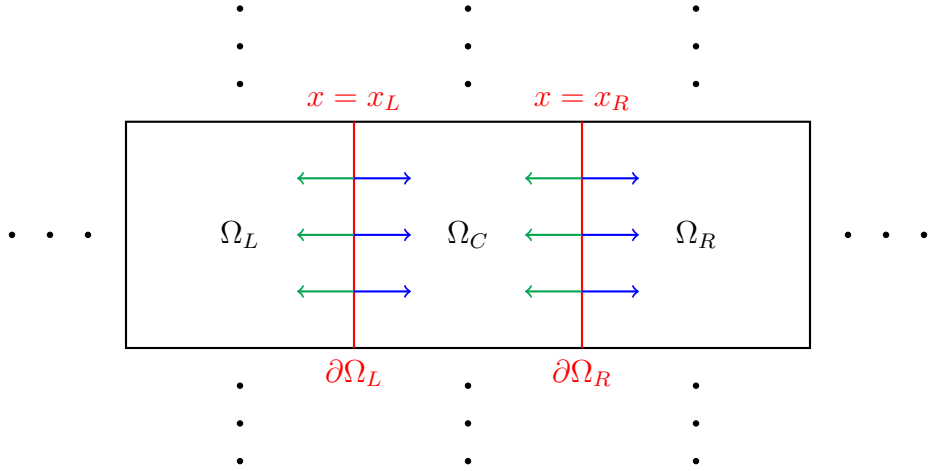

To solve equation \eqref{advx-evolK-orig} we need boundary information at the inflow boundaries. This data is obtained from the neighboring domains (see illustration in Figure \ref{fig::inflow_outflow}).  Let $\Omega_C$ be our current subdomain, $\Omega_L$ the subdomain to the left and $\Omega_R$ the subdomain to the right. If $\phi \leq \frac{\pi}{2}$ or $\phi \geq \frac{3\pi}{2}$, then our velocity in the $x$ direction is $\cos{(\phi)} \geq 0$, and thus we have flow to the right, indicated by the blue arrows. In this case, we have inflow at our left boundary $\partial\Omega_L$ and outflow on our right boundary $\partial\Omega_R$. If, on the other hand, $\frac{\pi}{2} < \phi < \frac{3\pi}{2}$, then the velocity in the $x$ direction is $\cos{(\phi)} < 0$, and we thus have flow to the left, indicated by the green arrows.

We denote the value of $f$ at the boundary by
\begin{equation}
  f^b: \mathbb{R}^+ \times \{x_L,x_R\} \times \Omega_C^Y \times [0,2\pi) \rightarrow \mathbb{R}^+.
\end{equation}
We then have on the left boundary
\begin{equation}\label{eq::fboundary_x_L}
  f^b(t,x_L,y,\phi) = \begin{cases}
    f_{\text{in}}(t,x_L,y,\phi), & \text{if }(\phi \leq \frac{\pi}{2}) \lor (\phi \geq \frac{3\pi}{2}), \\
    f(t,x_L,y,\phi), & \text{if } \frac{\pi}{2} < \phi < \frac{3\pi}{2},
  \end{cases}
\end{equation}
and on the right boundary
\begin{equation}\label{eq::fboundary_x_R}
  f^b(t,x_R,y,\phi) = \begin{cases}
    f_{\text{in}}(t,x_R,y,\phi), & \text{if } \frac{\pi}{2} < \phi < \frac{3\pi}{2}, \\
    f(t,x_R,y,\phi), & \text{if } (\phi \leq \frac{\pi}{2}) \lor (\phi \geq \frac{3\pi}{2}),
  \end{cases}
\end{equation}
where $f_{\text{in}}$ can be obtained by evaluating the low-rank representation of the neighboring subdomains at $x_L$ and $x_R$ (this is the inflow data). At the outflow boundary we simply have the solution provided by the low-rank representation on the current subdomain, i.e.~$f$ evaluated at $x_L$ and $x_R$, respectively.

Our goal, however, is to obtain boundary data for $K_j$, not for $f$. To accomplish this, we project $f^b(t,x,y,\phi)$ to the space spanned by the $V_j$ (as has been done in \cite{adapt_DLR_nonlin_boltz}). That is, for $x \in \{x_L, x_R\}$ we have
\begin{equation*}
  K_j^b(t,x,y) = \inner{f^b(t,x,y,\phi)}{V_j(t,\phi)}{\phi}.
\end{equation*}
Using equations \eqref{eq::fboundary_x_L} and \eqref{eq::fboundary_x_R} we get
\begin{align*}
  K_j^b(t,x_{L},y) &= \inner{f_{\text{in}}(t,x_{L},y,\phi) \mathds{1}_{(\phi \leq \frac{\pi}{2}) \lor (\phi \geq \frac{3\pi}{2})}}{V_j(t,\phi)}{\phi} \\
  &+ \inner{f(t,x_{L},y,\phi) \mathds{1}_{\frac{\pi}{2} < \phi < \frac{3\pi}{2}}}{V_j(t,\phi)}{\phi}
\end{align*}
and
\begin{align*}
  K_j^b(t,x_{R},y) &= \inner{f_{\text{in}}(t,x_{R},y,\phi) \mathds{1}_{\frac{\pi}{2} < \phi < \frac{3\pi}{2}}}{V_j(t,\phi)}{\phi} \\
  &+ \inner{f(t,x_{R},y,\phi) \mathds{1}_{(\phi \leq \frac{\pi}{2}) \lor (\phi \geq \frac{3\pi}{2})}}{V_j(t,\phi)}{\phi}.
\end{align*}
We can now use the low-rank representation $f=\suml{l}{}{K_l V_l}$ for the outflow boundary and the low-rank representations from the two neighboring domains, denoted by $f^L = \suml{l}{}{K_l^L V_l^L}$ and $f^R = \suml{l}{}{K_l^R V_l^R}$, respectively, to get
\begin{equation}\label{eq::generate_Kboundary}
    \begin{aligned}
    K_j^b(t,x_{L},y) &= \suml{l}{}{K_l^L(t,x_L,y) \inner{\mathds{1}_{(\phi \leq \frac{\pi}{2}) \lor (\phi \geq \frac{3\pi}{2})} V_l^L(t,\phi)}{ V_j(t,\phi)}{\phi}} \\
    &+ \suml{l}{}{{K}_l(t,x_{L},y)\inner{\mathds{1}_{\frac{\pi}{2} < \phi < \frac{3\pi}{2}} V_l(t,\phi)}{V_j(t,\phi)}{\phi}} \\
    K_j^b(t,x_{R},y) &= \suml{l}{}{K_l^R(t,x_R,y)\inner{ \mathds{1}_{\frac{\pi}{2} < \phi < \frac{3\pi}{2}} V_l^R(t,\phi) }{V_j(t,\phi)}{\phi}} \\
    &+ \suml{l}{}{{K}_l(t,x_{R},y)\inner{\mathds{1}_{(\phi \leq \frac{\pi}{2}) \lor (\phi \geq \frac{3\pi}{2})} V_l(t,\phi)}{V_j(t,\phi)}{\phi}},
    \end{aligned}
\end{equation}

Combining the boundary data for $K$ with equation \eqref{advx-evolK-orig} the K step of the projector splitting integrator for the $x$ advection is now given by
\begin{equation}\label{eq::advx_evolK}
  \partial_t K_j = -c_{\text{adv}}\suml{k}{}{\inner{V_j}{v_x V_k}{\phi}\partial_x K_k}, \qquad
  K_j \vert_{x \in \{x_L,x_R\}} = K_j^b.
\end{equation}

We slightly rewrite equation \eqref{advx-evolS-orig} for the S step such that we can directly use the boundary data which are expressed in $K$ and not in $U$. This gives
\begin{equation}\label{eq::advx_evolS}
  \partial_tS_{ij} = c_{\text{adv}} \suml{k}{}{\inner{U_i}{\partial_x K_k}{x,y}\inner{V_j}{v_xV_k}{\phi}}, \qquad
  K_j \vert_{x \in \{x_L,x_R\}} = K_j^b.
\end{equation}
We do the same for the L step, which gives
\begin{equation}\label{eq::advx_evolL}
  \partial_tL_i = -c_{\text{adv}}\suml{j}{}{\inner{U_i}{\partial_x K_j}{x,y}v_xV_j}, \qquad
  K_j \vert_{x \in \{x_L,x_R\}} = K_j^b.
\end{equation}
The corresponding projector splitting integrator for the $x$ advection is summarized in algorithm~\ref{alg:projec-advx}.

\begin{algorithm}[H]
  \caption{Projector splitting algorithm for $x$ advection. \label{alg:projec-advx}}
  \begin{algorithmic}
    \STATE \textbf{Input:} $U_i^0$, $S_{ij}^0$, $V_j^0$ (such that $f^0(x,y,\phi) = \suml{i,j}{}{U_i^0(x,y)S_{ij}^0 V_j^0(\phi)}$), $K_j^{L,0}$, $V_j^{L,0}$, $K_j^{R,0}$, $V_j^{R,0}$
    \STATE \textbf{Output:} $U_i^1$, $S_{ij}^3$, $V_j^1$ (such that $f^1(x,y,\phi) = \suml{i,j}{}{U_i^1(x,y)S_{ij}^3 V_j^1(\phi)}$)
  \end{algorithmic}
  \vspace{0.5em}
  \begin{algorithmic}[1]
    \STATE Solve equation \eqref{eq::advx_evolK} with initial value $\suml{i,j}{}{U_i^0 S_{ij}^0}$ and boundary value $K_j^b$ up to time $\Delta t$ to obtain $K_j^1$.
    \STATE Perform a QR decomposition of $K_j^1$ to get $U_i^1$ and $S_{ij}^1$.
    \STATE Solve equation \eqref{eq::advx_evolS} with initial value $S_{ij}^1$ and boundary value $K_j^b$ up to time $\Delta t$ to obtain $S_{ij}^2$.
    \STATE Solve equation \eqref{eq::advx_evolL} with initial value $\suml{j}{}{S_{ij}^2 V_j^0}$ and boundary value $K_j^b$ up to time $\Delta t$ to obtain $L_i^1$.
    \STATE Perform a QR decomposition of $L_i^1$ to get $V_j^1$ and $S_{ij}^3$.
  \end{algorithmic}
  \vspace{0.5em}
  The boundary values $K_j^b$ were calculated as in \eqref{eq::generate_Kboundary} using the low-rank representations from the neighboring subdomains $K_j^{L,0}$, $V_j^{L,0}$, $K_j^{R,0}$, and $V_j^{R,0}$.
\end{algorithm}

\subsection{Step 2: DLR advection in y}

The advection in the $y$ direction proceeds completely analogous to the advection in the $x$ direction described in the previous section, just that now the inflowing values come from the subdomains on the bottom and the top of our current subdomain. For completeness, we reproduce here the obtained evolution equations for the K step
\begin{equation}\label{eq::advy_evolK}
  \partial_t K_j = -c_{\text{adv}}\suml{k}{}{\inner{V_j}{v_yV_k}{\phi}\partial_y K_k},
  \qquad K_j \vert_{y \in \{y_B,y_T\}} = K^b_j,
\end{equation}
S step
\begin{equation}\label{eq::advy_evolS}
  \partial_tS_{ij} = c_{\text{adv}} \suml{k}{}{\inner{U_i}{\partial_y K_k}{x,y}\inner{V_j}{v_yV_k}{\phi}},
  \qquad K_j \vert_{y \in \{y_B,y_T\}} = K^b_j,
\end{equation}
and L step
\begin{equation}\label{eq::advy_evolL}
  \partial_tL_i = -c_{\text{adv}}\suml{j}{}{\inner{U_i}{\partial_y K_j}{x,y}v_yV_j},
  \qquad K_j \vert_{y \in \{y_B,y_T\}} = K^b_j.
\end{equation}
The boundary data $K^b_j$ is computed according to 
\begin{equation}\label{eq::generate_Kboundary_yadvection}
    \begin{aligned}
    K_j^b(t,x,y_B) &= \suml{l}{}{K_l^B(t,x,y_B) \inner{\mathds{1}_{0 \leq \phi \leq \pi} V_j^B(t,\phi)}{ V_j(t,\phi)}{\phi}} \\
    &+ \suml{l}{}{{K}_l(t,x,y_B)\inner{\mathds{1}_{\pi < \phi < 2\pi} V_l(t,\phi)}{V_j(t,\phi)}{\phi}} \\
    K_j^b(t,x,y_T) &= \suml{l}{}{K_l^T(t,x,y_T)\inner{ \mathds{1}_{\pi < \phi < 2\pi} V_l^T(t,\phi) }{V_j(t,\phi)}{\phi}} \\
    &+ \suml{l}{}{{K}_l(t,x,y_T)\inner{\mathds{1}_{0 \leq \phi \leq \pi} V_l(t,\phi)}{V_j(t,\phi)}{\phi}},
    \end{aligned}
\end{equation}
The corresponding projector splitting integrator for the $y$ advection is summarized in algorithm~\ref{alg:projec-advy}.

\begin{algorithm}[H]
  \caption{Projector splitting algorithm for $y$ advection. \label{alg:projec-advy}}
  \begin{algorithmic}
    \STATE \textbf{Input:} $U_i^0$, $S_{ij}^0$, $V_j^0$ (such that $f^0(x,y,\phi) = \suml{i,j}{}{U_i^0(x,y)S_{ij}^0 V_j^0(\phi)}$), $K_j^{B,0}$, $V_j^{B,0}$, $K_j^{T,0}$, $V_j^{T,0}$ 
    \STATE \textbf{Output:} $U_i^1$, $S_{ij}^3$, $V_j^1$ (such that $f^1(x,y,\phi) = \suml{i,j}{}{U_i^1(x,y)S_{ij}^3 V_j^1(\phi)}$)
  \end{algorithmic}
  \vspace{0.5em}
  \begin{algorithmic}[1]
    \STATE Solve equation \eqref{eq::advy_evolK} with initial value $\suml{i,j}{}{U_i^0 S_{ij}^0}$ and boundary value $K_j^b$ up to time $\Delta t$ to obtain $K_j^1$.
    \STATE Perform a QR decomposition of $K_j^1$ to get $U_i^1$ and $S_{ij}^1$.
    \STATE Solve equation \eqref{eq::advy_evolS} with initial value $S_{ij}^1$ and boundary value $K_j^b$ up to time $\Delta t$ to obtain $S_{ij}^2$.
    \STATE Solve equation \eqref{eq::advy_evolL} with initial value $\suml{j}{}{S_{ij}^2 V_j^0}$ and boundary value $K_j^b$ up to time $\Delta t$ to obtain $L_i^1$.
    \STATE Perform a QR decomposition of $L_i^1$ to get $V_j^1$ and $S_{ij}^3$.
  \end{algorithmic}
  \vspace{0.5em}
  The boundary values $K_j^b$ were calculated as in \eqref{eq::generate_Kboundary_yadvection} using the low-rank representations from the neighboring subdomains $K_j^{B,0}$, $V_j^{B,0}$, $K_j^{T,0}$, and $V_j^{T,0}$.
\end{algorithm}

\subsection{Step 3: DLR for absorption, scattering and the source term}

Here we can use a standard projector splitting operator as no boundary information is needed for the spatially pointwise absorption, scattering, and source terms. The equation for the K step is
\begin{equation}\label{eq::coll_evolK}
  \partial_t K_j = \frac{c_{\text{s}}}{2\pi}\suml{k}{}{K_k\innersingle{V_k}{\phi}\innersingle{V_j}{\phi}} - c_tK_j + Q \innersingle{V_j}{\phi}.
\end{equation}
For the S step we get
\begin{equation}\label{eq::coll_evolS}
  \partial_tS_{ij} = -\frac{1}{2\pi}\suml{k,l}{}{\inner{U_i}{c_{\text{s}}U_l}{x,y}\innersingle{V_k}{\phi}\innersingle{V_j}{\phi}S_{lk}} + \suml{l}{}{\inner{U_i}{c_{\text{t}}U_l}{x,y}S_{lj}}-\inner{U_i}{Q}{x,y}\innersingle{V_j}{\phi}
\end{equation}
and finally for the L step we have
\begin{equation}\label{eq::coll_evolL}
  \partial_tL_i = \frac{1}{2\pi}\suml{l}{}{\inner{U_i}{c_{\text{s}}U_l}{x,y}\innersingle{L_l}{\phi}} - \suml{l}{}{\inner{U_i}{c_{\text{t}}U_l}{x,y}L_l} + \inner{U_i}{Q}{x,y}.
\end{equation}
The corresponding projector splitting integrator for absorption, scattering and the source term is summarized in algorithm~\ref{alg:projec-coll}.

\begin{algorithm}[H]
  \caption{Projector splitting algorithm for absorption, scattering and the source term. \label{alg:projec-coll}}
  \begin{algorithmic}
    \STATE \textbf{Input:} $U_i^0$, $S_{ij}^0$, $V_j^0$ (such that $f^0(x,y,\phi) = \suml{i,j}{}{U_i^0(x,y)S_{ij}^0 V_j^0(\phi)}$) 
    \STATE \textbf{Output:} $U_i^1$, $S_{ij}^3$, $V_j^1$ (such that $f^1(x,y,\phi) = \suml{i,j}{}{U_i^1(x,y)S_{ij}^3 V_j^1(\phi)}$)
  \end{algorithmic}
  \vspace{0.5em}
  \begin{algorithmic}[1]
    \STATE Solve equation \eqref{eq::coll_evolK} with initial value $\suml{i,j}{}{U_i^0 S_{ij}^0}$ up to time $\Delta t$ to obtain $K_j^1$.
    \STATE Perform a QR decomposition of $K_j^1$ to get $U_i^1$ and $S_{ij}^1$.
    \STATE Solve equation \eqref{eq::coll_evolS} with initial value $S_{ij}^1$ up to time $\Delta t$ to obtain $S_{ij}^2$.
    \STATE Solve equation \eqref{eq::coll_evolL} with initial value $\suml{j}{}{S_{ij}^2 V_j^0}$ up to time $\Delta t$ to obtain $L_i^1$.
    \STATE Perform a QR decomposition of $L_i^1$ to get $V_j^1$ and $S_{ij}^3$.
  \end{algorithmic}
\end{algorithm}

\subsection{Augmentation and adaptive rank strategy}\label{sec::adaptive_rank}

As mentioned in the beginning of the section, it is not sufficient to simply run the derived algorithm with the obtained boundary data. We also have to ensure that the span of $V$ is sufficiently rich to accurately represent the boundary data from neighboring domains. To add the relevant information (contained in the $V$ basis of neighboring domains) is the task of the augmentation step described in this section. Augmentation is now a classic concept for dynamical low-rank methods. It is done in the augmented BUG integrator (as the classic variant without augmentation can be unstable; see, e.g., \cite{Einkemmer2023}). In this context it has also been used to get higher order schemes \cite{Ceruti2024} and to enforce conservation \cite{Einkemmer2023c,Baumann2023,Scalone2024} or boundary conditions \cite{Kusch2021}. Augmentation is most commonly encountered for the BUG type integrators, but has also been done for the projector splitting integrator, see, e.g., \cite{adapt_DLR_nonlin_boltz}. A crucial goal in the present context is to augment as few basis functions as necessary in order to keep computational and memory cost as low as possible.

Although, in principle, our augmentation procedure can be done in the continuous case, it is significantly more straightforward to describe it in the language of linear algebra. Thus, from now on we assume that we have already performed a spatial and velocity/angle discretization. Then our low-rank approximation $f=\suml{i,j}{}{U_i S_{ij} V_j}$ becomes
\[ F = U S V^T, \quad \text{where} \quad U \in \mathbb{R}^{n_x \cdot n_y \times r}, S \in \mathbb{R}^{r \times r}, V \in \mathbb{R}^{n_{\phi} \times r}.\]
The boundary data for $x$ advection are then given by
\[ F_L = K_L V_L^T, \quad \text{where} \quad K_L \in \mathbb{R}^{n_y \times r_L}, V_L \in \mathbb{R}^{n_{\phi} \times r_L} \]
and
\[ F_R = K_R V_R^T, \quad \text{where} \quad K_R \in \mathbb{R}^{n_y \times r_R}, V_R \in \mathbb{R}^{n_{\phi} \times r_R}. \]
Note that $r_L$ is the rank used on the left subdomain (i.e.~on  $\Omega_L$) and $r_R$ is the rank used on the right subdomain (i.e.~on $\Omega_R$).
More details on the specific discretization that we use in the numerical results can be found in section \ref{sec::num_results}.

In the following we will describe the augmentation procedure that we perform before Step 1 (i.e.~the $x$ advection step). A similar procedure is performed before step 2. It is clear (see figure \ref{fig::inflow_outflow}) that only information coming from $\partial\Omega_L$ and $\partial\Omega_R$ is relevant. The augmentation procedure that we use is given in algorithm \ref{alg:augmentation}. The main idea is that we keep information from $V_L$ and $V_R$ as long as neglecting it would introduce an error larger than a user specified tolerance $\text{tol}$ in the approximation space (i.e.~the span of $\hat{V}$). The output of the algorithm are new low-rank factors $U^1$, $S^1$ and $V^1$ that, in general, have a higher rank, i.e. the resulting rank $r_t$ satisfies $r_t \geq r$. We can then proceed with the $x$ advection step according to algorithm \ref{alg:projec-advx} with the new low-rank factors. This, in particular, ensures that any boundary information can be represented at least to accuracy $\text{tol}$.

\begin{algorithm}[H]
\begin{algorithmic}
  \STATE \textbf{Input:} $U^0$, $S^0$, $V^0$ (such that $F^0 = U^0S^0V^{0,T}$), $S_L$, $V_L$, $S_R$, $V_R$, tolerance $\text{tol}$
  \STATE \textbf{Output:} $U^1$, $S^1$, $V^1$ (such that $F^1 = U^1S^1V^{1,T}$)
\end{algorithmic}
\vspace{0.5em}
\begin{algorithmic}[1]
  \STATE Set $L^0=[V^0S^{0,T}, V_L S_L^T, V_R S_R^T]$.
  \STATE Perform an SVD to obtain $L^0 = \tilde{V} \tilde{\Sigma} \tilde{U}^T$, where $\tilde{V} \in \mathbb{R}^{n_{\phi} \times m}$, $\tilde{\Sigma} \in \mathbb{R}^{m \times m}$ and $\tilde{U} \in \mathbb{R}^{(r+r_L+r_R) \times m}$.
  \STATE According to some specified tolerance $\text{tol}$, calculate a truncated rank $r_t$: Choose $r_t$ as the minimal number such that 
  \[
    \left(\suml{j=r_t+1}{m}{\sigma_j^2}\right)^{\frac{1}{2}} \leq \text{tol}, 
  \]
  where $\sigma_j$ are the singular values in $\tilde{\Sigma}$. Additionally, impose $r_t \geq r$. 
  \STATE Truncate $\tilde{V}$ to obtain $\hat{V}=\tilde{V}(:,1:r_t) \in \mathbb{R}^{n_{\phi} \times r_t}$.
  \STATE Extend our basis $U^0$ with a random matrix $U_{\text{rand}} \in \mathbb{R}^{n_x \cdot n_{y} \times (r_t-r)}$ and extend $S^0$, such that $\hat{U}\hat{S}\hat{V}^T = U^0S^0V^{0,T} = F^0$ is satisfied, i.e.:
  \begin{equation*}
      \Hat{U} = [U^0, U_{\text{rand}}] \in \mathbb{R}^{n_{x} \cdot n_y \times r_t}, 
      \hspace{0.15cm} \Hat{S} = 
      \begin{bmatrix}
      S^0V^{0,T}\hat{V} \\
      0
      \end{bmatrix} \in \mathbb{R}^{r_t \times r_t}.
  \end{equation*}
  \STATE Orthonormalize the new basis by performing a QR decomposition $\Hat{U} = U^1 S_U$ and $\Hat{V} = V^1 S_V$. 
  Also, set $S^1 = S_U \Hat{S} S_V^T$.
\end{algorithmic}
\caption{Augmentation algorithm. \label{alg:augmentation}}
\end{algorithm}

Let us note that there are two important differences compared to \cite{adapt_DLR_nonlin_boltz}. First, we obtain our augmented basis directly from the low-rank representation of the neighboring subdomains. That is, we do not need to form $F_b$ (see step 1). This makes the algorithm computational more efficient. Second, we proceed with the computation by truncating $L^0$, which includes information of the $V$ basis for both the current as well as the neighboring domains. This implies that in situations where the basis used for the neighboring domains is similar to the current one, $r_t$ is only slightly larger than $r$. We will investigate this in more detail when discussing the numerical results.

Note that since the rank is increased in this step we have to truncate it at some point. Thus, after the advection in $x$ and $y$ we add a standard truncation step (see algorithm \ref{alg:truncate}) that reduces the rank such that the error is below the user specified tolerance $\text{tol}$. This automatically results in a rank adaptive scheme, where the rank can increase and decrease depending on the properties of the solution. 

\begin{algorithm}[H]
\caption{Truncation algorithm. \label{alg:truncate}}
\begin{algorithmic}
  \STATE \textbf{Input:} $U^0$, $S^0$, $V^0$ (such that $F^0 = U^0 S^0 V^{0,T}$), tolerance $\text{tol}$
  \STATE \textbf{Output:} $U^1$, $S^1$, $V^1$ (such that $F^1 = U^1 S^1 V^{1,T}$)
\end{algorithmic}
\vspace{0.5em}
\begin{algorithmic}[1]
  \STATE Perform an SVD to obtain $S^0 = \tilde{U}\tilde{\Sigma}\tilde{V}^T$, where $\tilde{U}, \tilde{V} \in \mathbb{R}^{r \times r}$ are orthonormal and $\tilde{\Sigma} \in \mathbb{R}^{r \times r}$ is diagonal with descending singular values.

  \STATE According to some specified tolerance $\text{tol}$, calculate a truncated rank $r^{\prime}$: 
  Choose $r^{\prime}$ as the minimal number such that 
  \[
    \left( \sum_{j=r'+1}^{r} \sigma_j^2 \right)^{1/2} \le \text{tol},
  \]
  where $\sigma_j$ are the singular values in $\Tilde{\Sigma}$. \\ Then, set $\Bar{U} = \Tilde{U}(:, 1:r^{\prime})$, $\Bar{V} = \Tilde{V}(:, 1:r^{\prime})$ and $S^1 = \Tilde{\Sigma}(1:r^{\prime}, 1:r^{\prime})$.

  \STATE Update the basis as:
  \[
    U^1 = U^0 \bar{U} \in \mathbb{R}^{n_x \cdot n_y \times r'}, \quad
    V^1 = V^0 \bar{V} \in \mathbb{R}^{n_\phi \times r'}.
  \]

\end{algorithmic}
\end{algorithm}

\subsection{Complete domain decomposition dynamical low-rank algorithm}

Putting all the building blocks developed in the previous sections together we arrive at our domain decomposition low-rank algorithm that is detailed in algorithm \ref{alg::ddlr_algorithm}.

\begin{algorithm}[H]
\caption{Domain decomposition low-rank algorithm. \label{alg::ddlr_algorithm}}
\begin{algorithmic}
  \STATE \textbf{Input:} Low rank representation on each domain, {\color{black}tolerance $\text{tol}$}
  \STATE \textbf{Output:} Low rank representation on each domain
\end{algorithmic}
\vspace{0.5em}
\begin{algorithmic}[1]
  \STATE Add basis functions according to some {\color{black}tolerance $\text{tol}$} and algorithm \ref{alg:augmentation}.
  \STATE Run algorithm \ref{alg:projec-advx} to solve $x$ advection.
  \STATE Remove basis functions according to some {\color{black}tolerance $\text{tol}$} and algorithm \ref{alg:truncate}.
  \vspace{1em}
  \STATE Add basis functions for inflowing values in $y$ direction according to some {\color{black}tolerance $\text{tol}$}, similar to algorithm \ref{alg:augmentation}.
  \STATE Run algorithm \ref{alg:projec-advy} to solve $y$ advection.
  \STATE Remove basis functions according to some {\color{black}tolerance $\text{tol}$} and algorithm \ref{alg:truncate}.
  \vspace{1em}
  \STATE Run algorithm \ref{alg:projec-coll} to solve absorption, scattering and the source term.
\end{algorithmic}
\vspace{0.5em}
Note: This procedure is performed on each domain and repeated for every time step.
\end{algorithm}

{
\color{black}
\subsection{Complexity analysis}

Finally, we want to compare our new domain decomposition low-rank algorithm \ref{alg::ddlr_algorithm} to the classic low-rank approach from section \ref{sec::dlra}. For this, we consider the complexity of both algorithms, which is determined by storage cost, computational and memory transaction cost and communication cost. To carry out the complexity analysis in the most general way, we will not consider the specific $2\times 1$ dimensional case, but rather a $d_x$ dimensional spatial domain and a $d_{\phi}$ dimensional angular domain. We let $n$ be the number of grid points in each dimension and $r$ the rank at a certain time step. Note that $r$ changes during a time step due to our rank adaptivity, but for simplicity of our analysis we will always consider $r$ instead of $r_t$ and $r^{\prime}$ from algorithms \ref{alg:augmentation} and \ref{alg:truncate}. Also, the rank for domain decomposition is typically smaller than the one for simulations on one domain, which can be seen in the numerical examples. Additionally, we will assume a splitting into $s^{d_x}$ subdomains of equal size for the domain decomposition algorithm.

\subsubsection{Storage cost}

Let us start by looking at the storage cost. For the classic low-rank approach, we need to save the low-rank factors $U$, $S$ and $V$. The overall storage cost of the classic DLRA is then given by $\mathcal{O}(n^{d_x}  r + r^2 + r  n^{d_{\phi}})$. Here the most expensive part results from the $U$ matrix, which is of size $n^{d_x} \times r$, because for the radiative transfer equation usually $d_x > d_{\phi}$. Assuming the rank $r$ is sufficiently small, i.e., $r \ll n$, this is significantly more efficient than for the full solver, where the storage cost would be $\mathcal{O}(n^{d_x}n^{d_{\phi}})$. For the domain decomposition algorithm, the storage requirement on each node (i.e. for each subdomain) is given by $\mathcal{O}(\left(\frac{n}{s}\right)^{d_x}  r + r^2 + r  n^{d_{\phi}})$.

\subsubsection{Computational and memory transaction cost}

Here we consider the cost for computation and memory transactions on a single node. As we mentioned above, for domain decomposition it is sufficient to consider the cost for a simulation on one node, because of the parallelism. Note that, on a modern CPU computing architecture, computational cost is roughly 40-80 times cheaper than memory access. Still, we consider both. The memory transaction cost is similar to the storage cost. For the classic low-rank approach we get $\mathcal{O}(n^{d_x}  r + r^2 + r  n^{d_{\phi}})$, while for the domain decomposition algorithm we get $\mathcal{O}(\left(\frac{n}{s}\right)^{d_x}  r + r^2 + r  n^{d_{\phi}})$. The computational cost for the main steps of the classic low-rank approach are listed in table \ref{tab::classic_dlra_cost}. For the domain decomposition algorithm, we have the same steps for each subdomain, but only have $\left(\frac{n}{s}\right)^{d_x}$ spatial points in each subdomain, resulting in a lower computational cost for those steps compared to the classic approach. Thus, in table \ref{tab::ddlr_cost}, we omit those steps and only show the additional steps and the total cost (including the omitted steps) of one time step for the domain decomposition low-rank algorithm. When comparing the total computational cost with the total memory transaction cost, we get for both the classic and the domain decomposition low-rank algorithm:
\begin{equation*}
  \frac{\text{computational cost}}{\text{memory transaction cost}} = O(r).
\end{equation*}
This means, depending on the rank and the architecture of the system, the algorithm can be either memory bound or compute bound. 

\begin{table}[H]

  \color{black}

  \centering
  \caption{Computational cost of classic low-rank algorithm from section \ref{sec::dlra}.}
  \label{tab::classic_dlra_cost}

  \renewcommand{\arraystretch}{1.2}

  \begin{tabular}{|l|c|}
  \hline
  \textbf{Operations} & \textbf{Floating point operations (flops)} \\

  \hline\hline

  Preparation of $K$
  & $O(n^{d_x}  r^2)$ \\
  \hline
  Preparation of $L$
  & $O(n^{d_{\phi}}  r^2)$ \\
  \hline
  Preparation of $\innersingle{V}{\phi}$
  & $O(n^{d_{\phi}}  r)$ \\
  \hline
  Preparation of $\innersingle{L}{\phi}$
  & $O(n^{d_{\phi}}  r)$ \\
  \hline
  Preparation of $\inner{V}{vV}{\phi}$
  & $O(n^{d_{\phi}}  r^2)$ \\
  \hline
  Preparation of $\inner{U}{\nabla_{x,y}U}{x,y}$
  & $O(n^{d_x}  r^2)$ \\
  \hline
  Total cost of preparation steps
  & $O(n^{d_x}  r^2 + n^{d_{\phi}}  r^2)$ \\

  \hline\hline

  Evaluation of right-hand side for $K$ step
  & $O(n^{d_x}  r^2)$ \\
  \hline
  Solving $K$ step with RK4
  & $O(n^{d_x}  r)$ \\
  \hline
  QR decomposition of $K$
  & $O(n^{d_x}  r^2)$ \\
  \hline
  Total cost of $K$ step
  & $O(n^{d_x}  r^2)$ \\

  \hline\hline

  Evaluation of right-hand side for $S$ step
  & $O(r^3)$ \\
  \hline
  Solving $S$ step with RK4
  & $O(r^2)$ \\
  \hline
  Total cost of $S$ step
  & $O(r^3)$ \\

  \hline\hline

  Evaluation of right-hand side for $L$ step
  & $O(n^{d_{\phi}}  r^2)$ \\
  \hline
  Solving $L$ step with RK4
  & $O(n^{d_{\phi}}  r)$ \\
  \hline
  QR decomposition of $L$
  & $O(n^{d_{\phi}}  r^2)$ \\
  \hline
  Total cost of $L$ step
  & $O(n^{d_{\phi}}  r^2)$ \\

  \hline\hline

  Total cost of one time step
  & $O(n^{d_x} r^2 + n^{d_{\phi}} r^2 + r^3)$ \\
  \hline

  \end{tabular}
\end{table}

\begin{table}[H]

  \color{black}

  \centering
  \caption{Computational cost of domain decomposition low-rank algorithm \ref{alg::ddlr_algorithm}.}
  \label{tab::ddlr_cost}

  \renewcommand{\arraystretch}{1.2}
  \begin{tabular}{|l|c|}
  \hline
  \textbf{Operations} & \textbf{Floating point operations (flops)} \\

  \hline\hline

  Boundary data preparation using eq \ref{eq::generate_Kboundary}
  & $O\left(\left(\frac{n}{s}\right)^{d_x-1}  r^2 + n^{d_{\phi}}  r^2\right)$ \\

  \hline\hline

  Preparation of $L^0$
  & $O(n^{d_{\phi}}  r^2)$ \\
  \hline
  Reduced SVD of $L^0$
  & $O(n^{d_{\phi}}  r^2)$ \\
  \hline
  QR decomposition of $\hat{U}$
  & $O\left(\left(\frac{n}{s}\right)^{d_x}  r^2\right)$ \\
  \hline
  QR decomposition of $\hat{V}$
  & $O(n^{d_{\phi}}  r^2)$ \\
  \hline
  Calculation of $S^1$ 
  & $O(r^3)$ \\
  \hline
  Total cost of augmentation algorithm \ref{alg:augmentation}
  & $O\left(n^{d_{\phi}}  r^2 + \left(\frac{n}{s}\right)^{d_x} r^2 + r^3\right)$ \\

  \hline\hline

  Reduced SVD of $S^0$
  & $O(r^3)$ \\
  \hline
  Calculation of $U^1$ and $V^1$
  & $O\left(\left(\frac{n}{s}\right)^{d_x} r^2 + n^{d_{\phi}} r^2\right)$ \\
  \hline
  Total cost of truncation algorithm \ref{alg:truncate}
  & $O\left(\left(\frac{n}{s}\right)^{d_x}  r^2 + n^{d_{\phi}}  r^2 + r^3\right)$ \\

  \hline\hline

  Total cost of one time step
  & $O\left(n^{d_{\phi}}  r^2 + \left(\frac{n}{s}\right)^{d_x}  r^2 + r^3\right)$ \\
  \hline

  \end{tabular}
\end{table}

\subsubsection{Communication cost}

Now we consider the additional cost we get, when running a domain decomposition simulation on multiple nodes. To run the simulation on a subdomain, we need to communicate data for all the advection steps. In each advection step we need to exchange the boundary data in $K$ from the 2 neighboring subdomains in the advection direction. The size of the data that we need to exchange is $\mathcal{O}(\left(\frac{n}{s}\right)^{d_x-1}  r)$. This we need to do for all $d_x$ advection steps, which still gives us a total communication cost of $\mathcal{O}(\left(\frac{n}{s}\right)^{d_x-1}  r)$. For relatively large $n$, this is significantly smaller than the memory transaction cost of the domain decomposition algorithm:
\begin{equation*}
  \frac{\text{memory transaction cost}}{\text{communication cost}} = O\left(\frac{n}{s}\right).
\end{equation*}
Thus, as long as the number of grid points is sufficiently large, the communication cost is negligible compared to the memory transaction cost. This is a significant advantage of our domain decomposition low-rank algorithm, as it allows us to run simulations on multiple nodes without incurring a significant communication overhead.

}

\section{Numerical results}\label{sec::num_results}

In this section we consider several well known test cases of the radiative transfer equation in order to demonstrate the effectiveness of the proposed method. We compare the proposed algorithm \ref{alg::ddlr_algorithm} with a standard dynamical low-rank approach that does not make use of domain decomposition. {\color{black} This standard approach is based on the evolution equations \eqref{eq::Kstep}, \eqref{eq::Sstep} and \eqref{eq::Lstep} from section \ref{sec::dlra} and is implemented in a way that allows for rank adaptivity, which lets us compare the two methods directly. Note that in this standard low-rank approach we do not use the ideas from section \ref{sec::algorithm}, meaning we do not apply operator splitting to our radiative transfer equation. We also compare the results to a reference solution, which is computed with the same dynamical low-rank algorithm without domain decomposition. The only difference is, that with the reference solution we use a very small tolerance for rank adaptivity, resulting in a very high simulation rank and thus an accurate solution.}

We first study a classic lattice and Hohlraum problem and show that our algorithm is more efficient than the classic low-rank approach. We then study a point source {\color{black} and a line source} problem, where the rank on some part of the domain needs to be large in order to obtain an accurate solution. {\color{black} Those examples} show that the low-rank domain decomposition algorithm can work efficiently in situations where classic low-rank methods require large memory and computational effort. The source code necessary to reproduce the numerical experiments can be found in \cite{Brunner2025_code}.

All simulations conducted here use the classic Runge--Kutta method of order four (RK4) to perform the time integration. For the spatial discretization {\color{black} a finite difference upwind scheme of first order is used}. This is, in particular, necessary for the problems described here as the discontinuous interface between cells can otherwise lead to unphysical oscillations. For more details we refer to appendix \ref{sec::appendix_upwind}.

\subsection{Lattice test problem}\label{sec::lattice}

\begin{figure}[H]
    \centering
    \begin{subfigure}{0.48\textwidth}
        \centering
        \includegraphics[width=\textwidth]{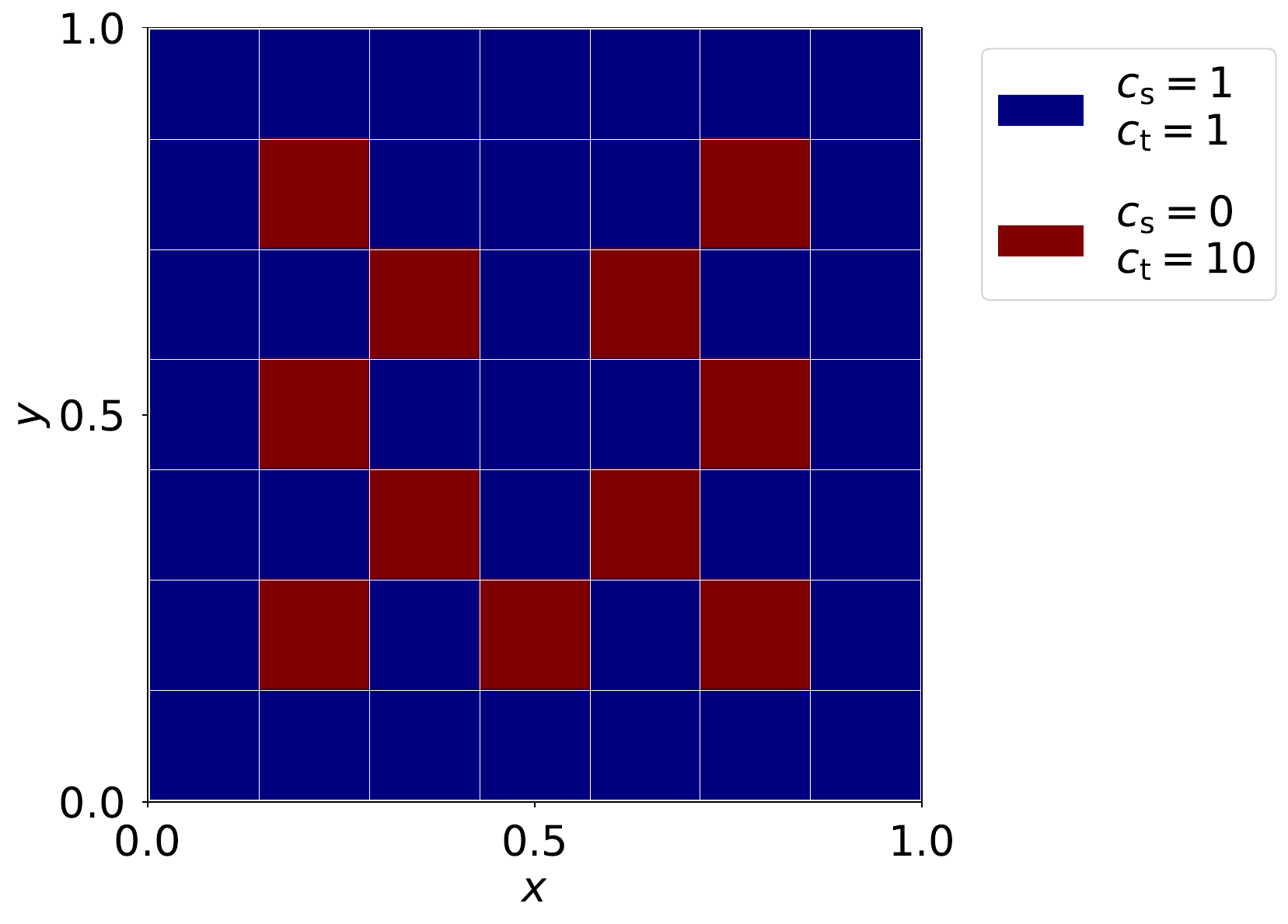}
    \end{subfigure}
    \begin{subfigure}{0.51\textwidth}
        \centering
        \includegraphics[width=\textwidth]{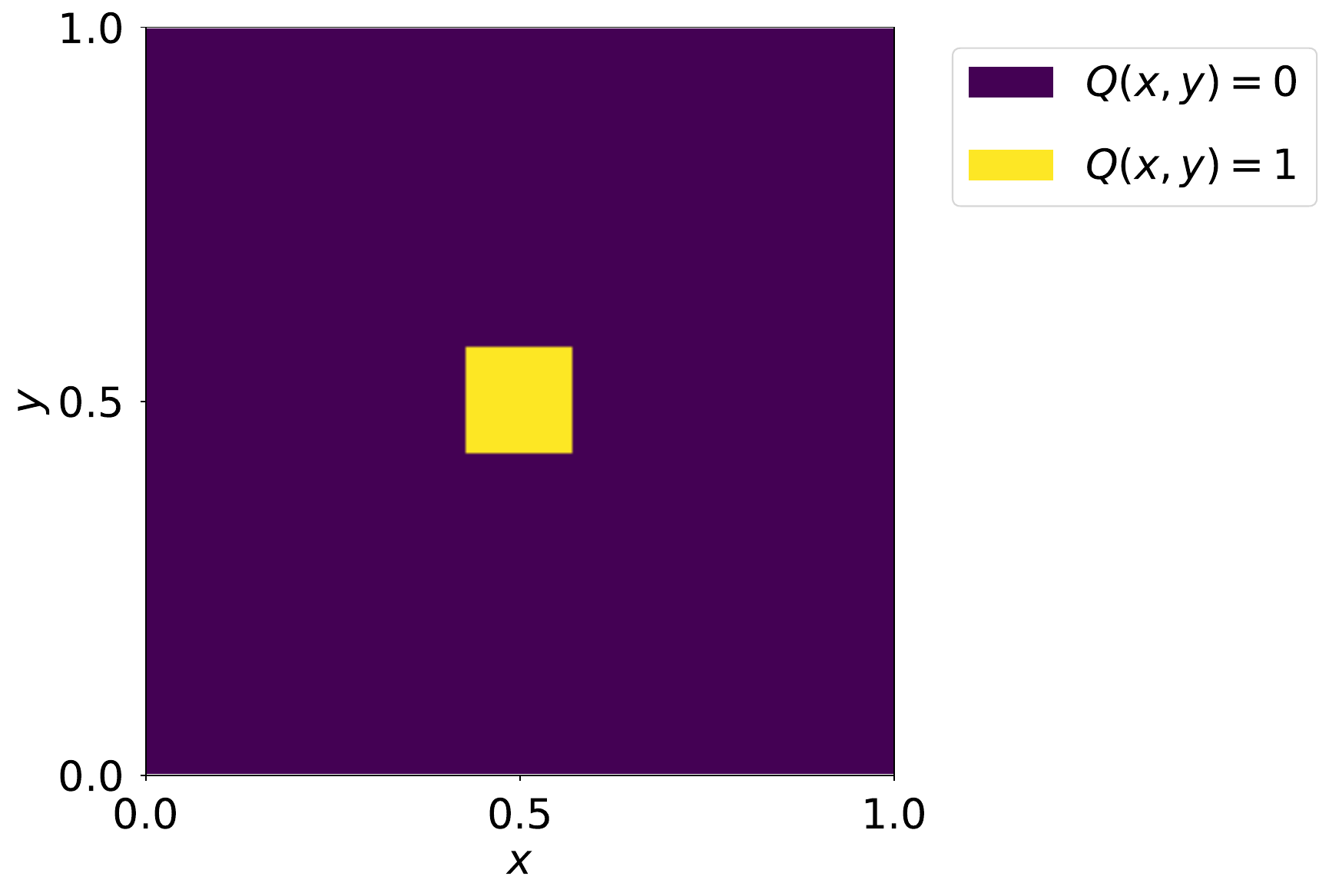}
    \end{subfigure}
    \caption{Setup of the lattice test problem.}
    \label{fig::lattice_setup}
\end{figure}

Our first problem is the lattice test case that has often been considered in the literature 
(see, e.g., \cite{Brunner_testproblem,McClarren_testproblem,Peng_lowrank_testproblem,Parallel_Gianluca_testproblem}). It is based on the blocks of a nuclear reactor core that has both highly absorbing and scattering regions.
Our setup can be seen in figure \ref{fig::lattice_setup}. The full domain is one centimeter wide, i.e. $x\in [0,1]$ and $y \in [0,1]$. For the angle we have $\phi \in [0,2\pi]$. The spatial domain is then split up into a grid containing $7\times 7$ blocks. We choose our unit system such that $c_{\text{adv}}=1$. The choice of the other coefficients can be seen in figure \ref{fig::lattice_setup}. The red blocks describe a highly absorbing material with $c_{\text{s}}=0$ and $c_{\text{t}}=10$, while the blue blocks describe a material that scatters but has no absorption ($c_{\text{s}}=1$ and $c_{\text{t}}=1$). The right of figure \ref{fig::lattice_setup} describes the source $Q(x,y)$. For the yellow block we have $Q(x,y)=1$, while in the remainder of the domain we set $Q(x,y)=0$. The initial condition is {\color{black}$f(t=0,x,y,\phi)=0$}.

\begin{figure}[H]
    \centering
    \begin{subfigure}{\textwidth}
        \centering
        \begin{minipage}{0.12\textwidth}
          \raggedleft
          \caption{$t=0.3$:}
        \end{minipage}
        \begin{minipage}{0.424\textwidth}
          \includegraphics[width=\textwidth]{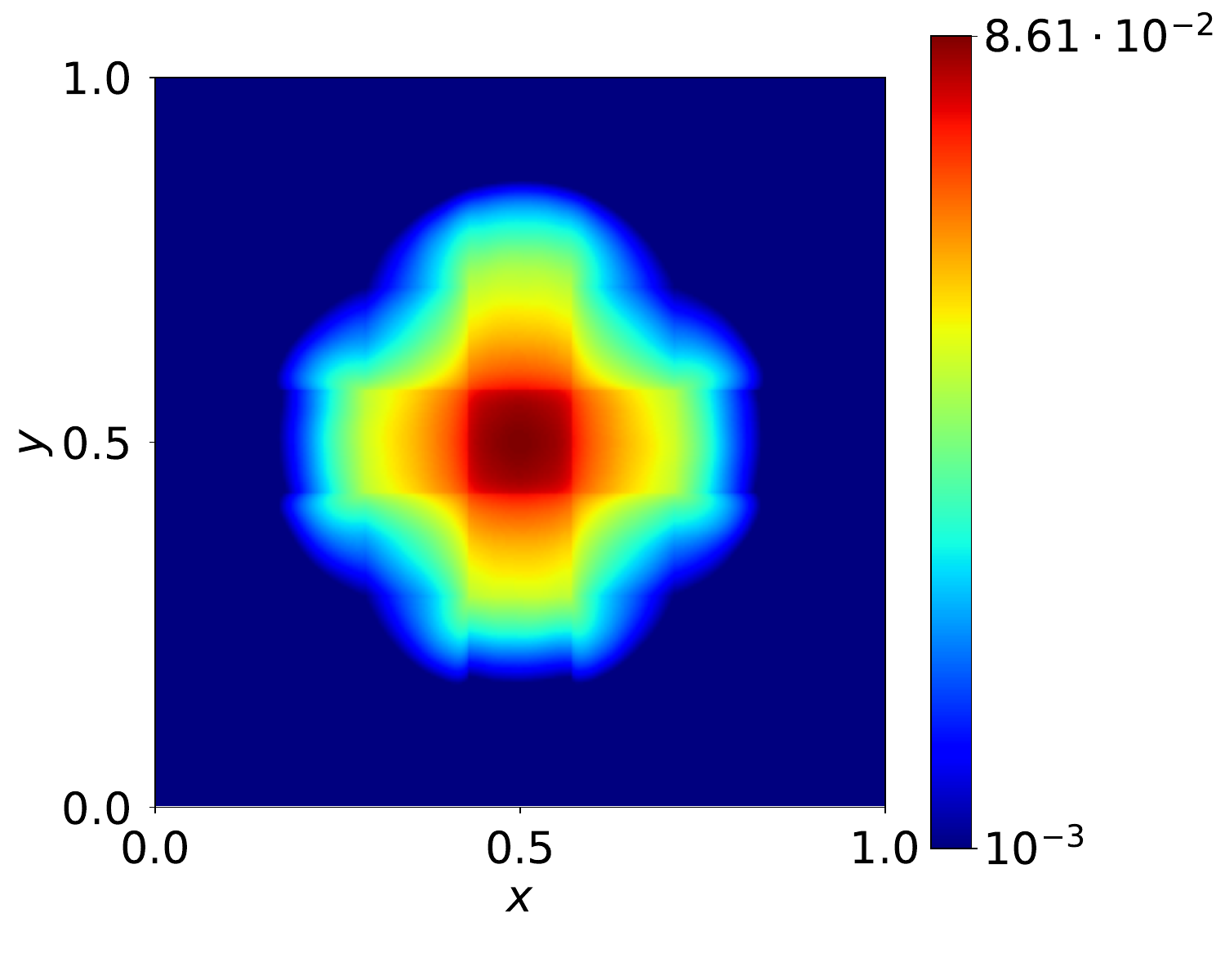}
        \end{minipage}
        \begin{minipage}{0.424\textwidth}
            \centering
            \includegraphics[width=\textwidth]{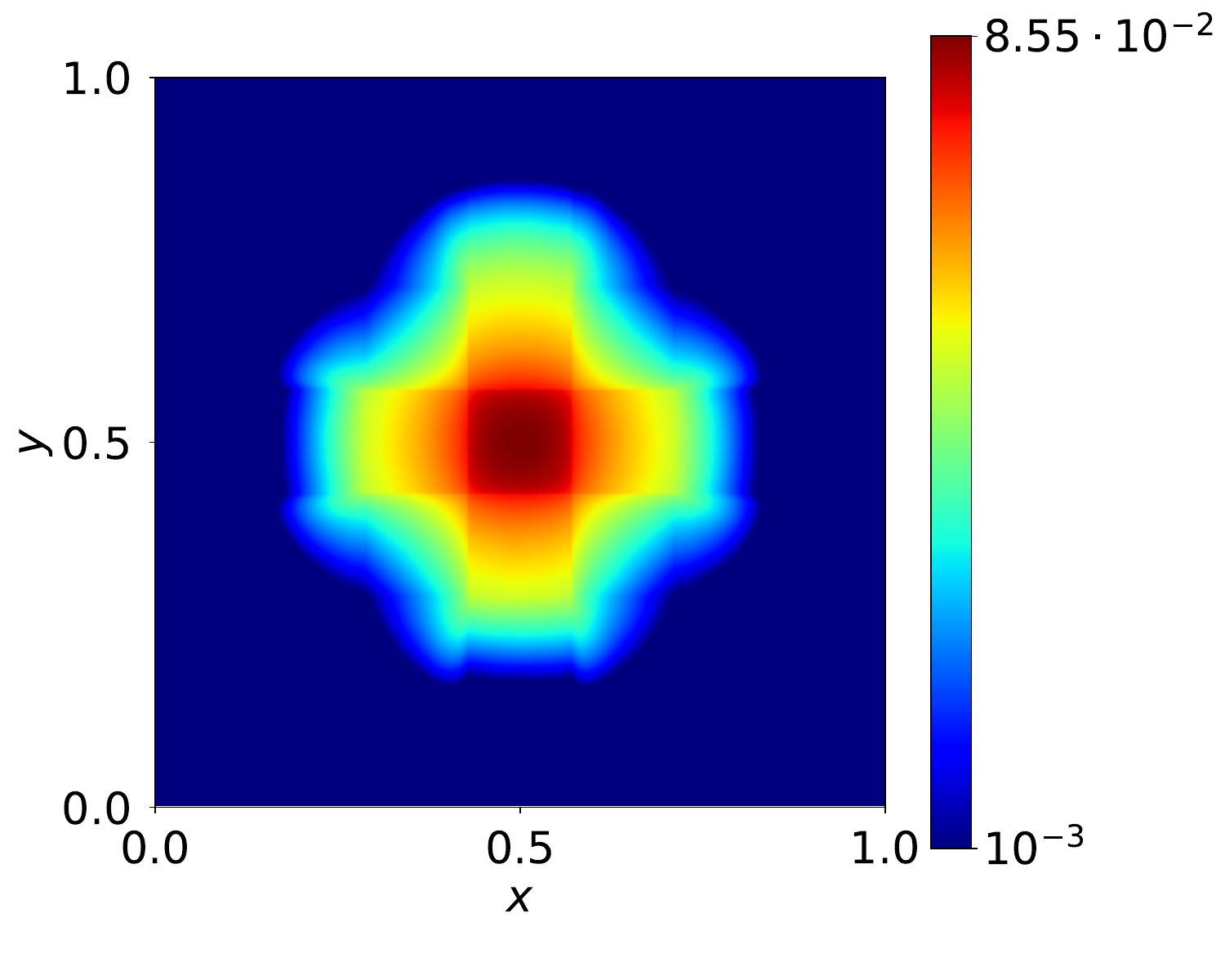}
        \end{minipage}
    \end{subfigure}
    
    \begin{subfigure}{\textwidth}
        \centering
        \begin{minipage}{0.12\textwidth}
          \raggedleft
          \caption{$t=0.5$:}
        \end{minipage}
        \begin{minipage}{0.424\textwidth}
          \includegraphics[width=\textwidth]{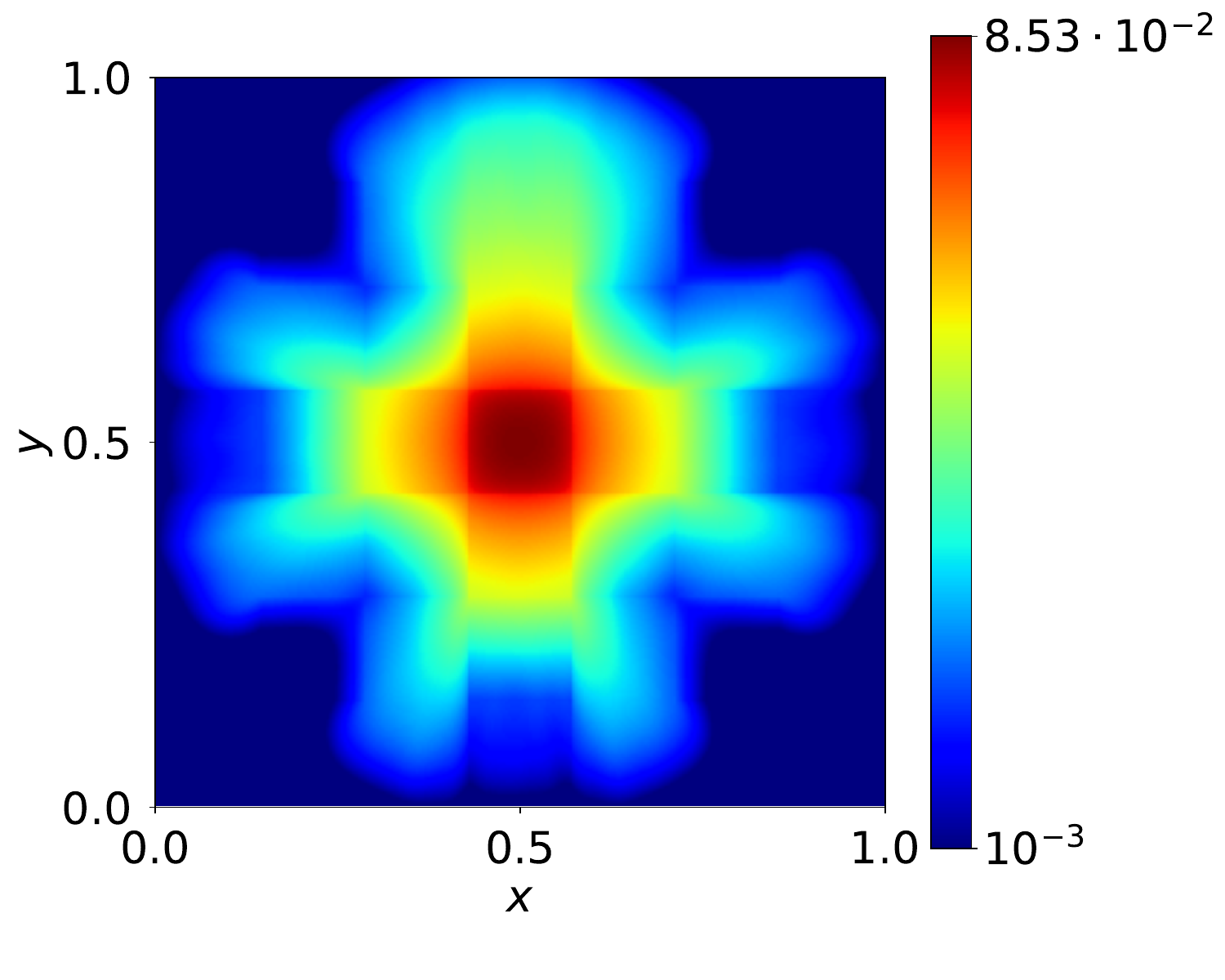}
        \end{minipage}
        \begin{minipage}{0.424\textwidth}
            \centering
            \includegraphics[width=\textwidth]{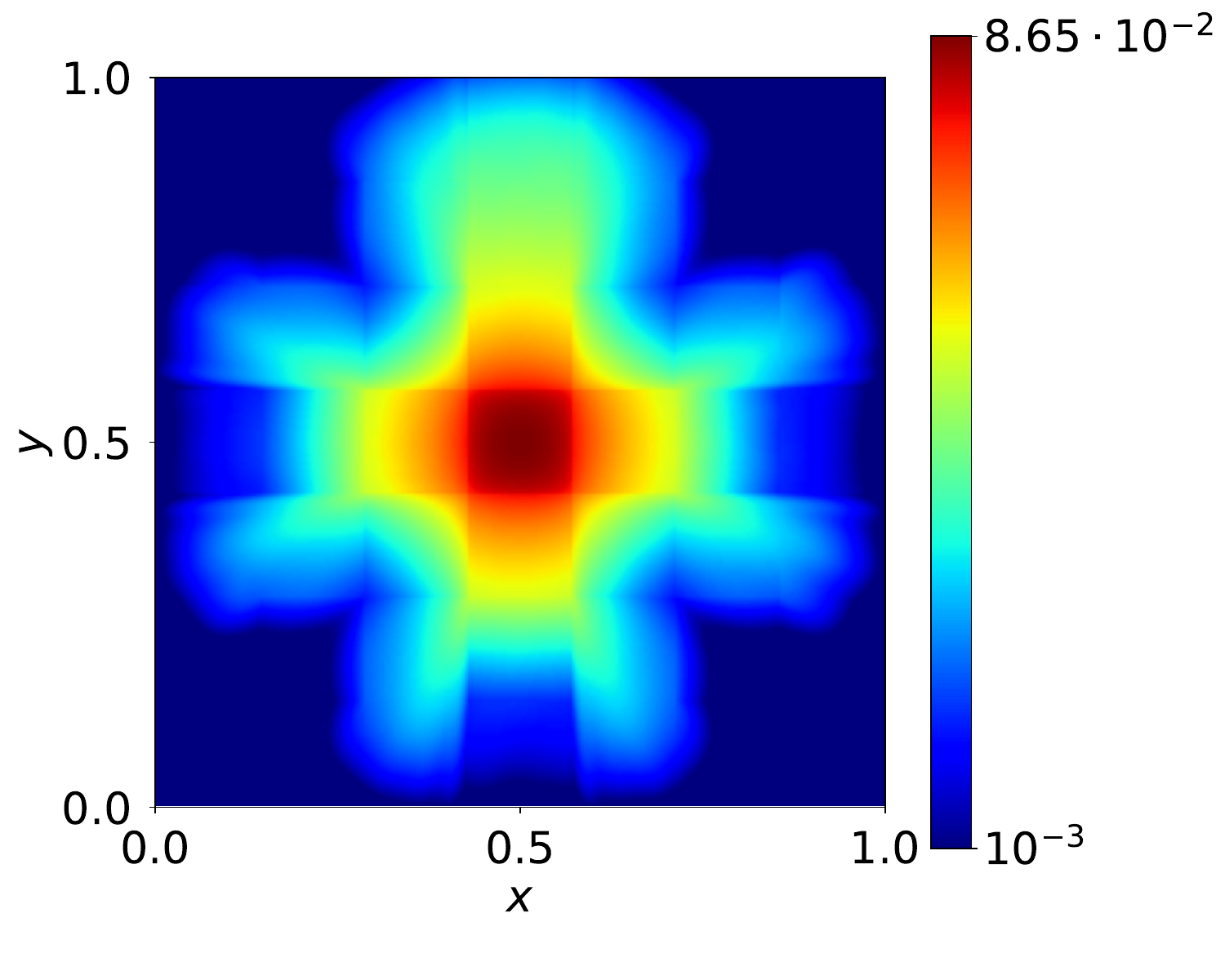}
        \end{minipage}
    \end{subfigure}
    
    \begin{subfigure}{\textwidth}
        \centering
        \begin{minipage}{0.12\textwidth}
          \raggedleft
          \caption{$t=0.7$:}
        \end{minipage}
        \begin{minipage}{0.424\textwidth}
          \includegraphics[width=\textwidth]{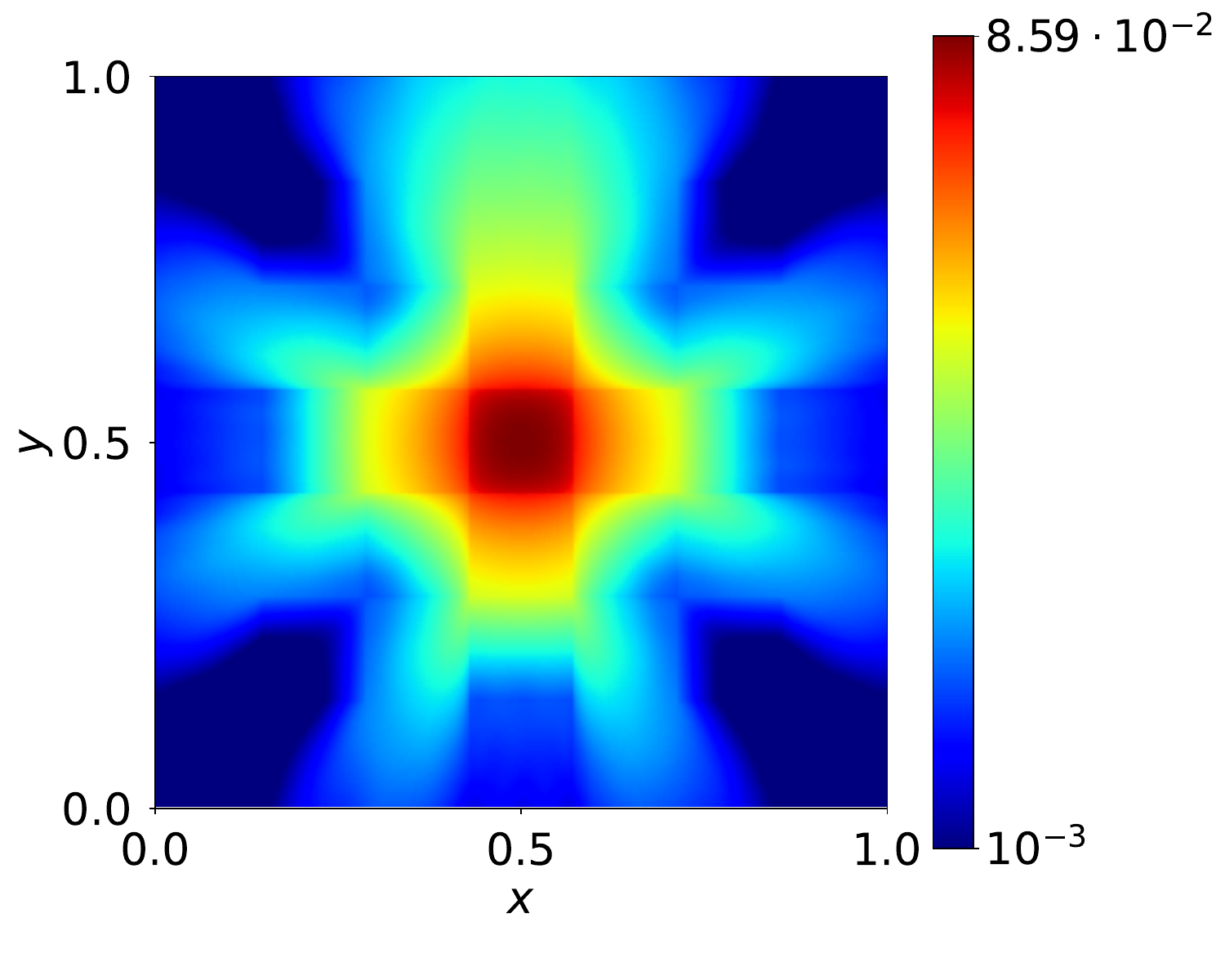}
        \end{minipage}
        \begin{minipage}{0.424\textwidth}
            \centering
            \includegraphics[width=\textwidth]{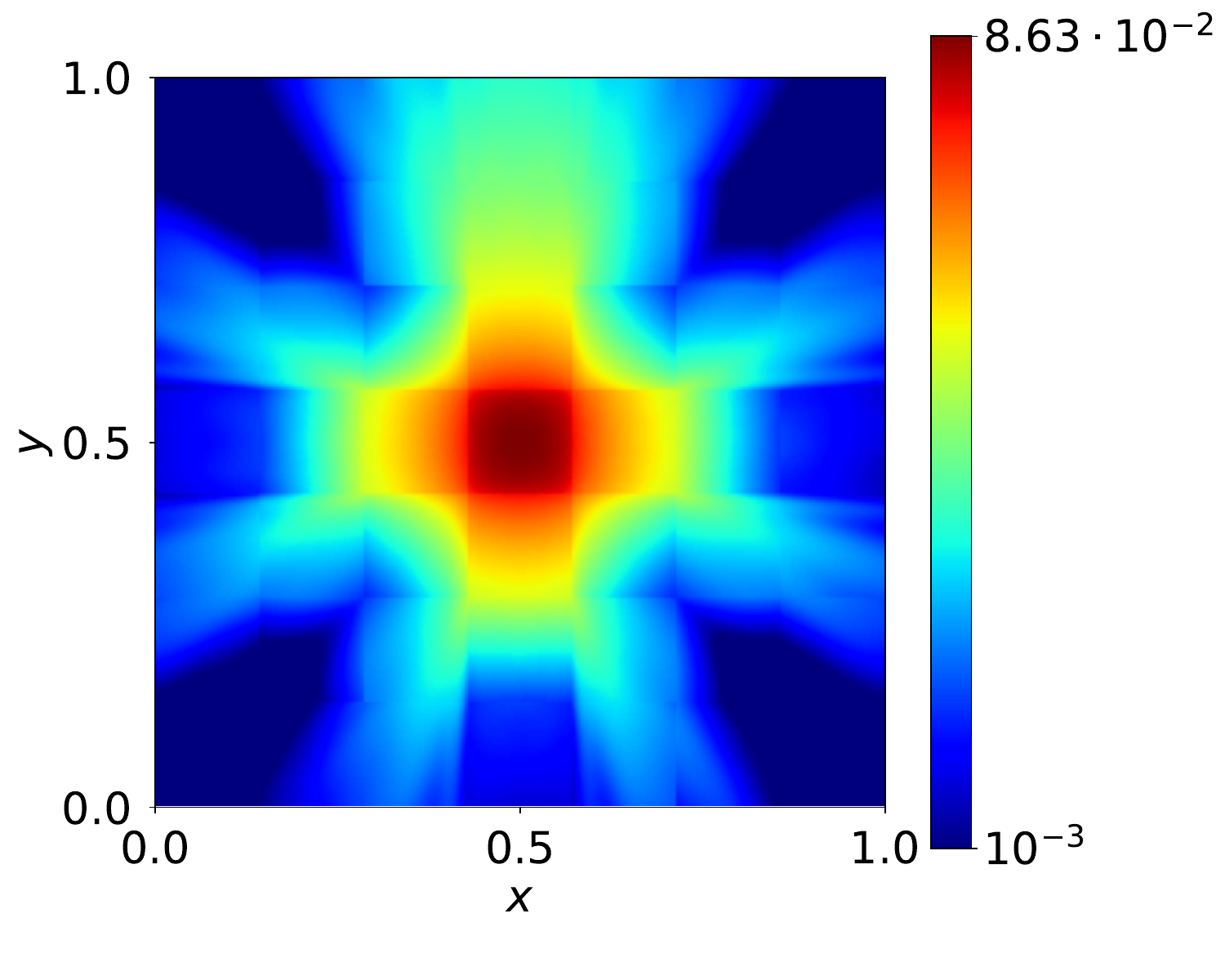}
        \end{minipage}
    \end{subfigure}
    \caption{We show $\rho (t,x,y)$ obtained by solving the lattice test problem with a classic dynamical low-rank algorithm (i.e.~with a single low-rank representation for the entire domain) on the left and with the domain decomposition low-rank algorithm \ref{alg::ddlr_algorithm} on the right. For the latter, $7 \times 7$ subdomains, each containing $36$ grid points in both the $x$ and $y$ directions, have been used.}
    \label{fig::lattice_rho_plots}
\end{figure}

We compare our domain decomposition approach to a classic dynamical low-rank algorithm (i.e.~without domain decomposition). The grid is discretized equidistantly, such that we have 252 grid points in $x$, $y$ and $\phi$. The time step is chosen to be $\Delta t = 0.5 \cdot \Delta x$. On the boundary we have an outflow boundary condition. The density $\rho$ for times $t=0.3$, $t=0.5$ and $t=0.7$ is shown on the left-hand side of figure \ref{fig::lattice_rho_plots}. As expected, the radiation induced by the source term  spreads outwards, moving around the absorbing blocks. A tolerance {\color{black} $\text{tol}= 5.5\cdot10^{-5}$}, which is used to perform rank adaptivity, is used {\color{black}for the classic low-rank} simulations. 

We compare this with our proposed domain decomposition approach which is shown on the right of figure \ref{fig::lattice_rho_plots}. The plots agree well. {\color{black}The tolerance used for the rank-adaptivity of the domain decomposition low-rank simulation is $
\text{tol} = 10^{-5}$. This yields a similar error for both methods.} Comparing to a reference solution $f_{\text{ref}}$ that has been computed with the classic low-rank method without domain decomposition and $\text{tol} = 10^{-10}$, we get
\begin{equation*}
  \norm{f-f_{\text{ref}}} := \frac{\norm{f-f_{\text{ref}}}_{L^2}}{\norm{1}_{L^2}} \approx {\color{black}2.07\cdot 10^{-4}}.
\end{equation*}
The error of the low-rank method without domain decomposition is approximately {\color{black} $2.04 \cdot 10^{-4}$} (compared to the same reference solution). {\color{black}Note that we have chosen the tolerance in the simulations such that we get roughly the same overall error. This required the choice of different tolerances for the two methods.} Note that for the classic dynamical low-rank scheme we use a different approach to achieve rank adaptivity (as there is no inflow in this case), which can be found in \cite{Hochbruck_rankadaptivity}.

\begin{figure}[H]
    \centering
    \includegraphics[width=0.44\textwidth]{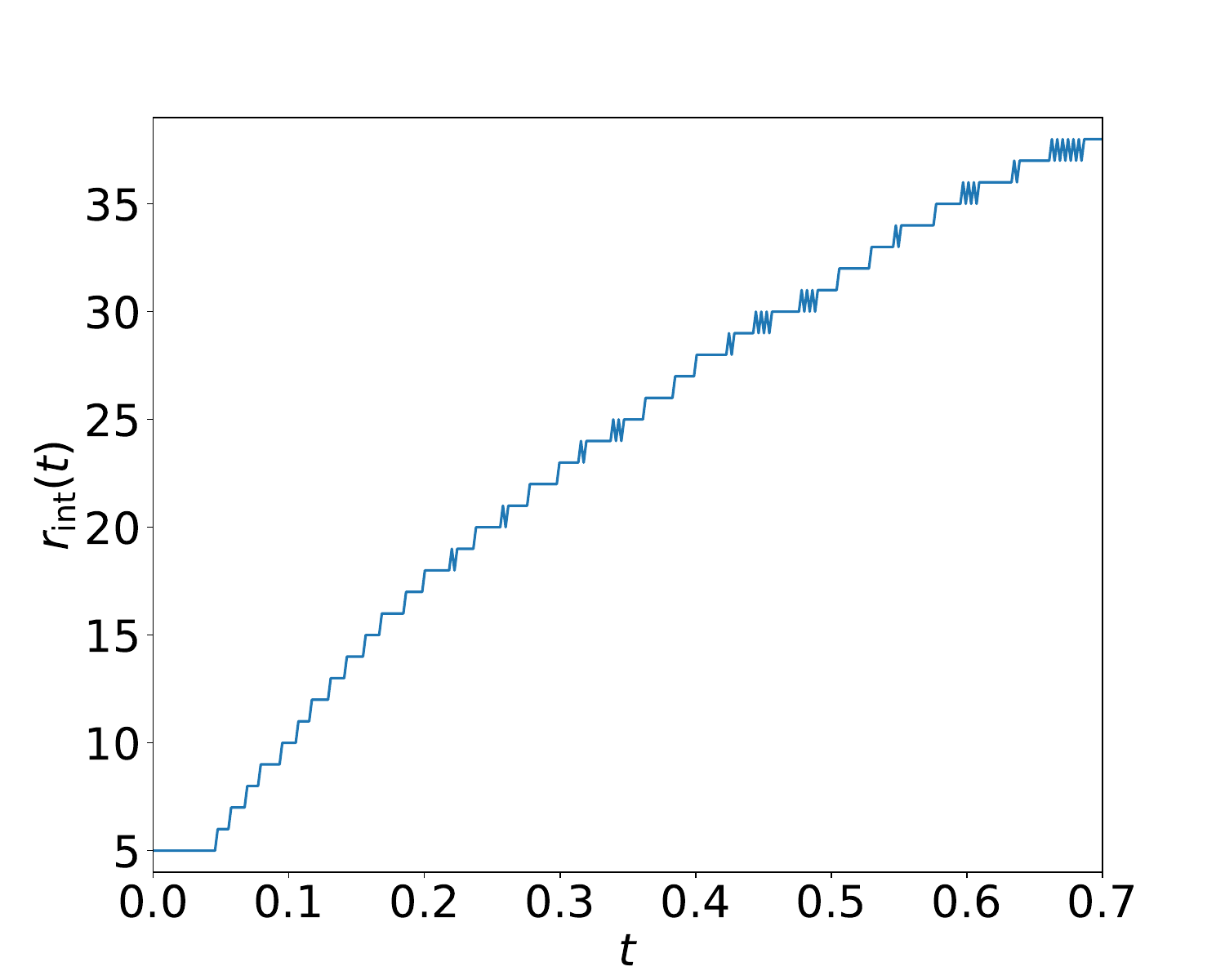}
    \caption{We show the time evolution of the rank for the classic low-rank simulation of the lattice test problem.  As an adaptive rank scheme we used the scheme described in \cite{Hochbruck_rankadaptivity} with {\color{black} $\text{tol}=5.5\cdot10^{-5}$}.}
    \label{fig::lattice_1domain_ranks}
\end{figure}

\begin{figure}[H]
    \centering
    \begin{subfigure}{0.44\textwidth}
        \centering
        \includegraphics[width=\textwidth]{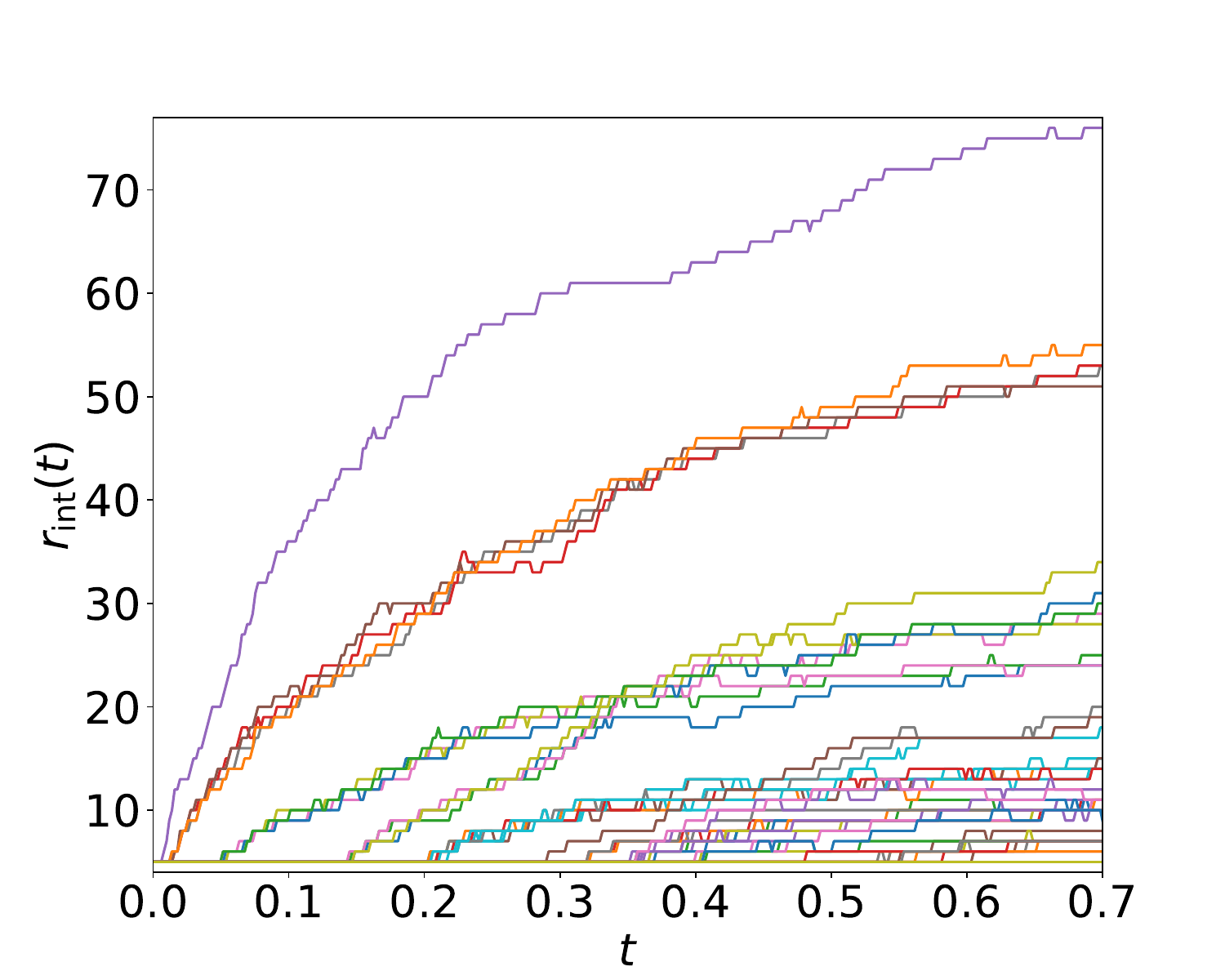}
        \caption{Intermediate ranks}
    \end{subfigure}
    \begin{subfigure}{0.46\textwidth}
        \centering
        \includegraphics[width=\textwidth]{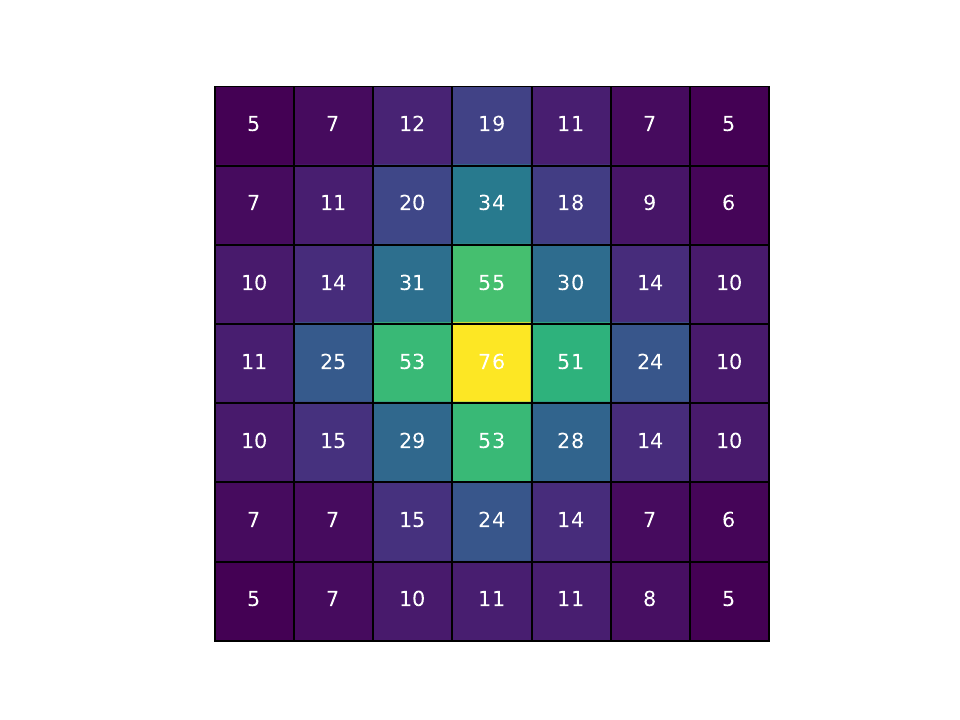}
        \caption{Intermediate ranks at time $t=0.7$}
    \end{subfigure}
    \begin{subfigure}{0.44\textwidth}
        \centering
        \includegraphics[width=\textwidth]{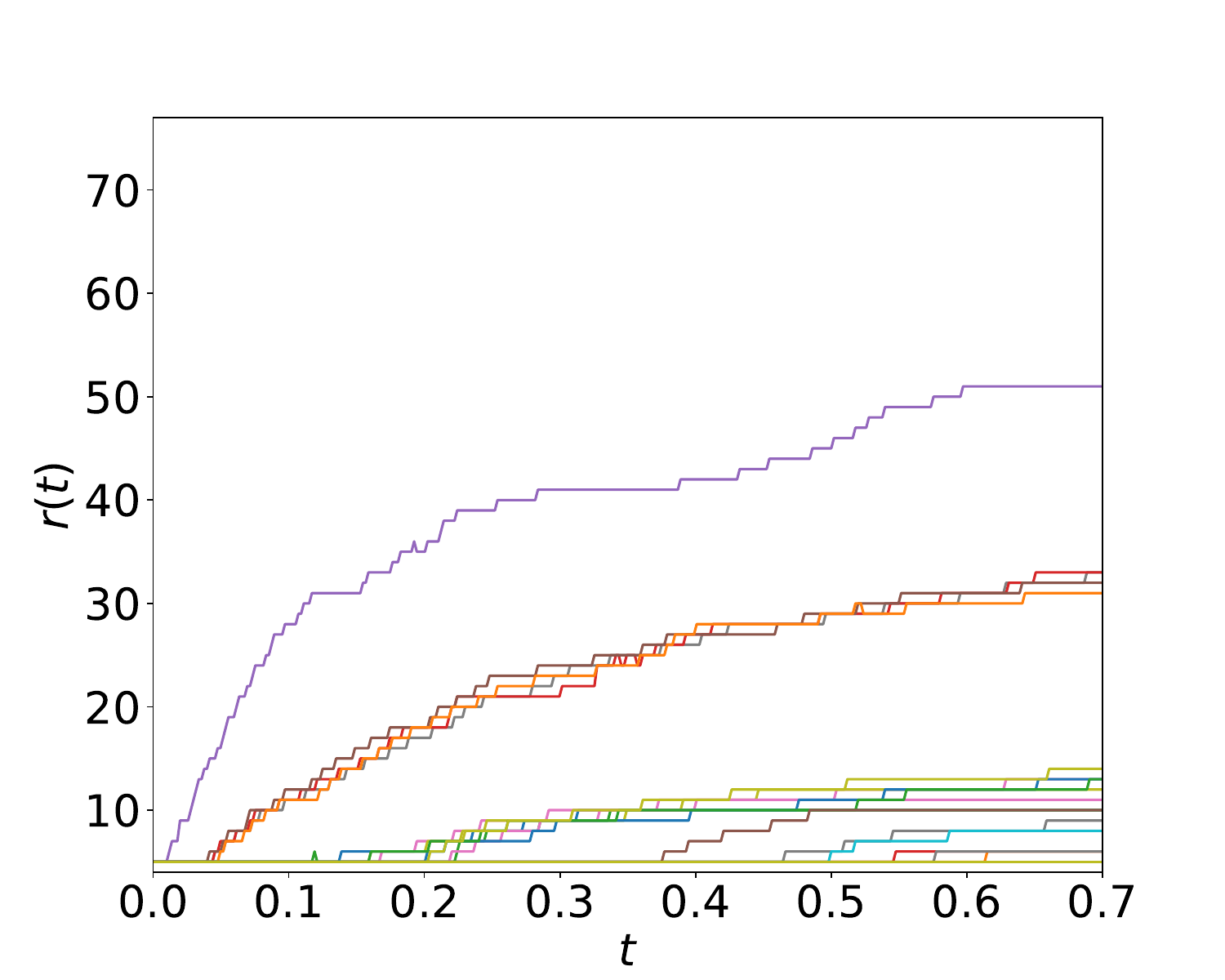}
        \caption{Ranks}
    \end{subfigure}
    \begin{subfigure}{0.46\textwidth}
        \centering
        \includegraphics[width=\textwidth]{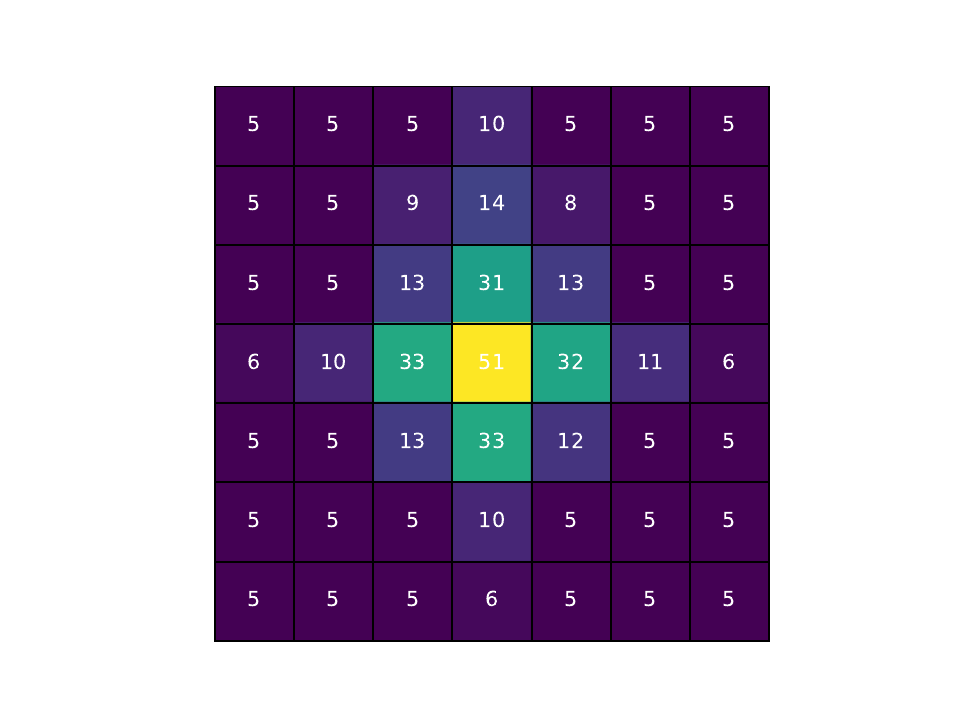}
        \caption{Ranks at time $t=0.7$}
    \end{subfigure}
    \caption{We show the time evolution of the intermediate ranks and ranks after the time step for the domain decomposition dynamical low-rank simulation of the lattice test problem. The adaptive rank scheme as proposed in algorithms \ref{alg:augmentation} and \ref{alg:truncate} with {\color{black} $\text{tol}=10^{-5}$} was used.}
    \label{fig::lattice_domaindecomp_ranks}
\end{figure}

We now compare the ranks that are required to obtain the described results. For the low-rank method without domain decomposition the rank as a function of time is given in figure \ref{fig::lattice_1domain_ranks}. As expected, the rank increases steadily over time. The highest rank required is {\color{black} $38$} at the end of our simulation. For the domain decomposition approach the rank is different in each block. The ranks in each block as a function of time is shown in figure \ref{fig::lattice_domaindecomp_ranks}. Here we have to distinguish between the rank required to store our solution (bottom of the plot; this would be the approximation that is written to disk) and the highest rank that occurs during the computation, also called intermediate rank (top of the plot; this is the rank that determines the amount of memory required to run the simulation). 

\begin{figure}[H]
    \centering
    \includegraphics[width=0.45\textwidth]{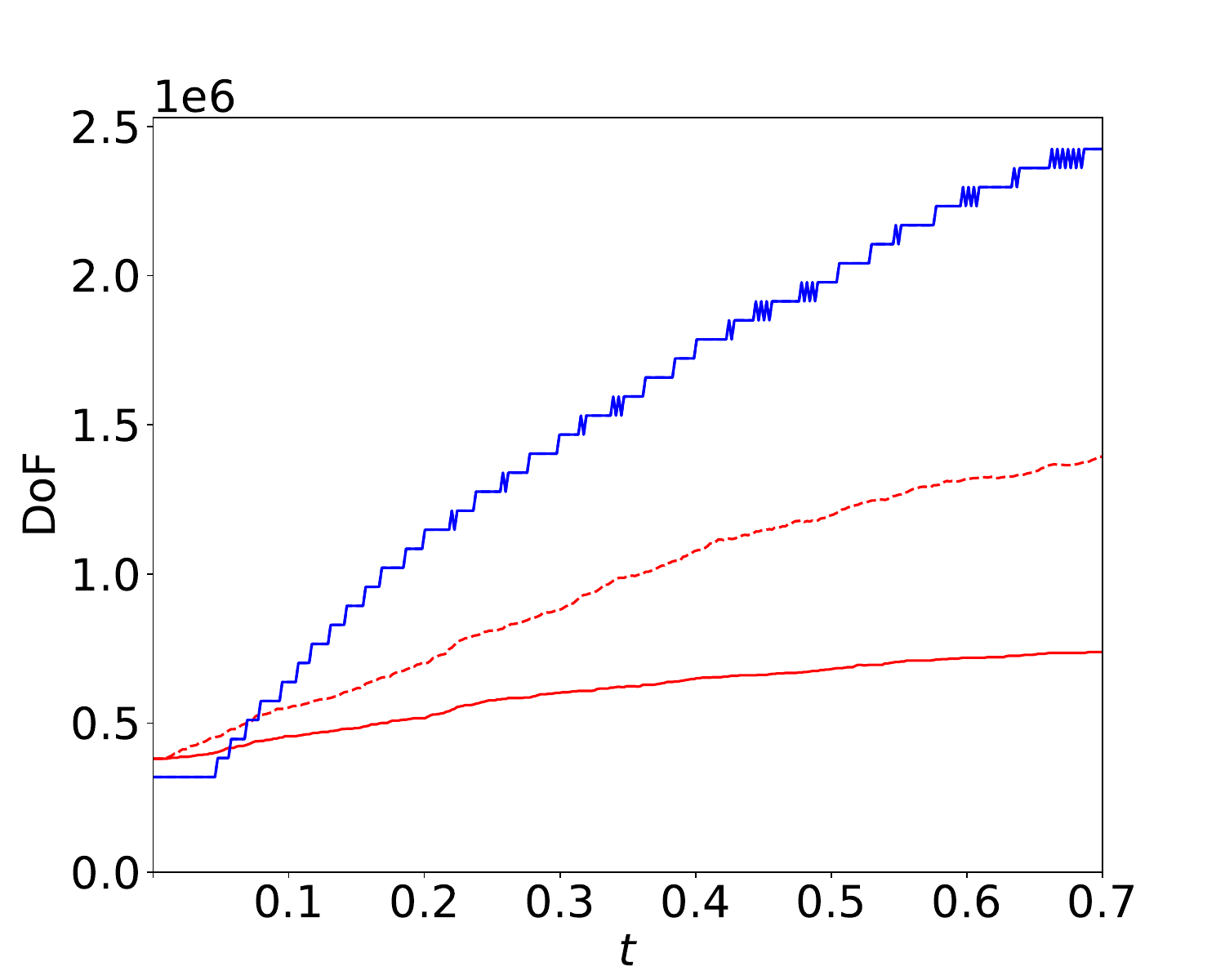}
    \includegraphics[width=0.45\textwidth]{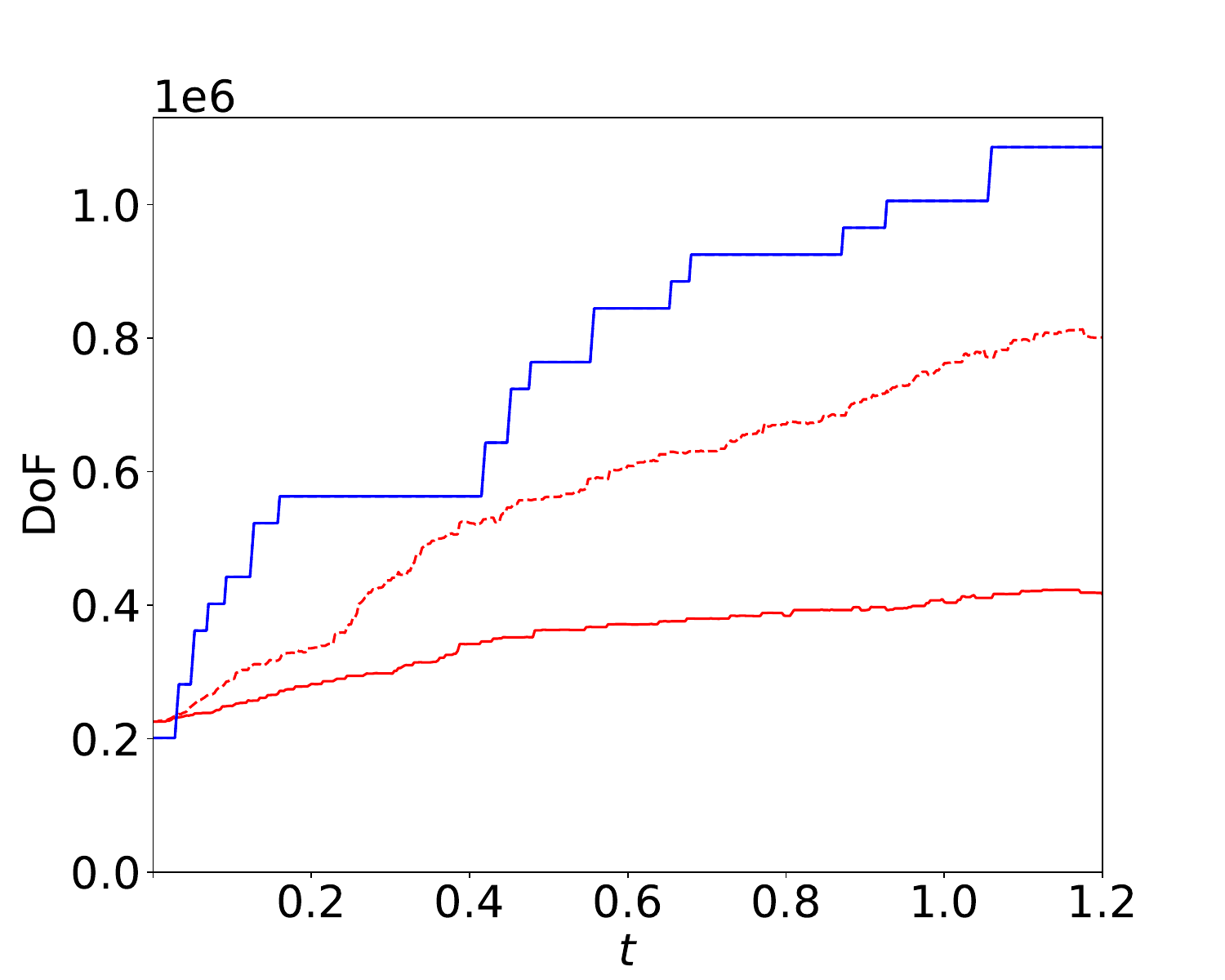}
    \caption{We show degrees of freedom (DoF) as a function of time for the low-rank simulation of the lattice test problem (left) and the Hohlraum test problem (right) for the classic low-rank method (blue) and the domain decomposition low-rank method (algorithm~\ref{alg::ddlr_algorithm}) (red).}
    \label{fig::lattice_dof_comparison}
\end{figure}

We note that only in {\color{black} the central blocks} the rank required by the domain decomposition low-rank method is as large as what is required for the classic low-rank method. In most other blocks the rank required is significantly lower. This is because the low-rank representation only needs to resolve a more localized region of space. Moreover, in some blocks the solution is naturally low-rank (e.g.~due to strong absorption or because not much is happening in that block at early times of the simulation). Thus, the domain decomposition approach reduces the number of degrees of freedom required to run the simulation significantly. This is shown in figure \ref{fig::lattice_dof_comparison}, where we see  roughly a factor of {\color{black} $1.7$} reduction when comparing maximum memory usage (i.e.~we compare with the intermediate rank, given by the dotted line) and roughly a factor of {\color{black} four} less memory to store the solution (i.e.~we compare with the rank obtained after each time step, given by the full line).

\subsection{Hohlraum test problem}\label{sec::hohlraum}

\begin{figure}[H]
    \centering
    \includegraphics[width=0.48\textwidth]{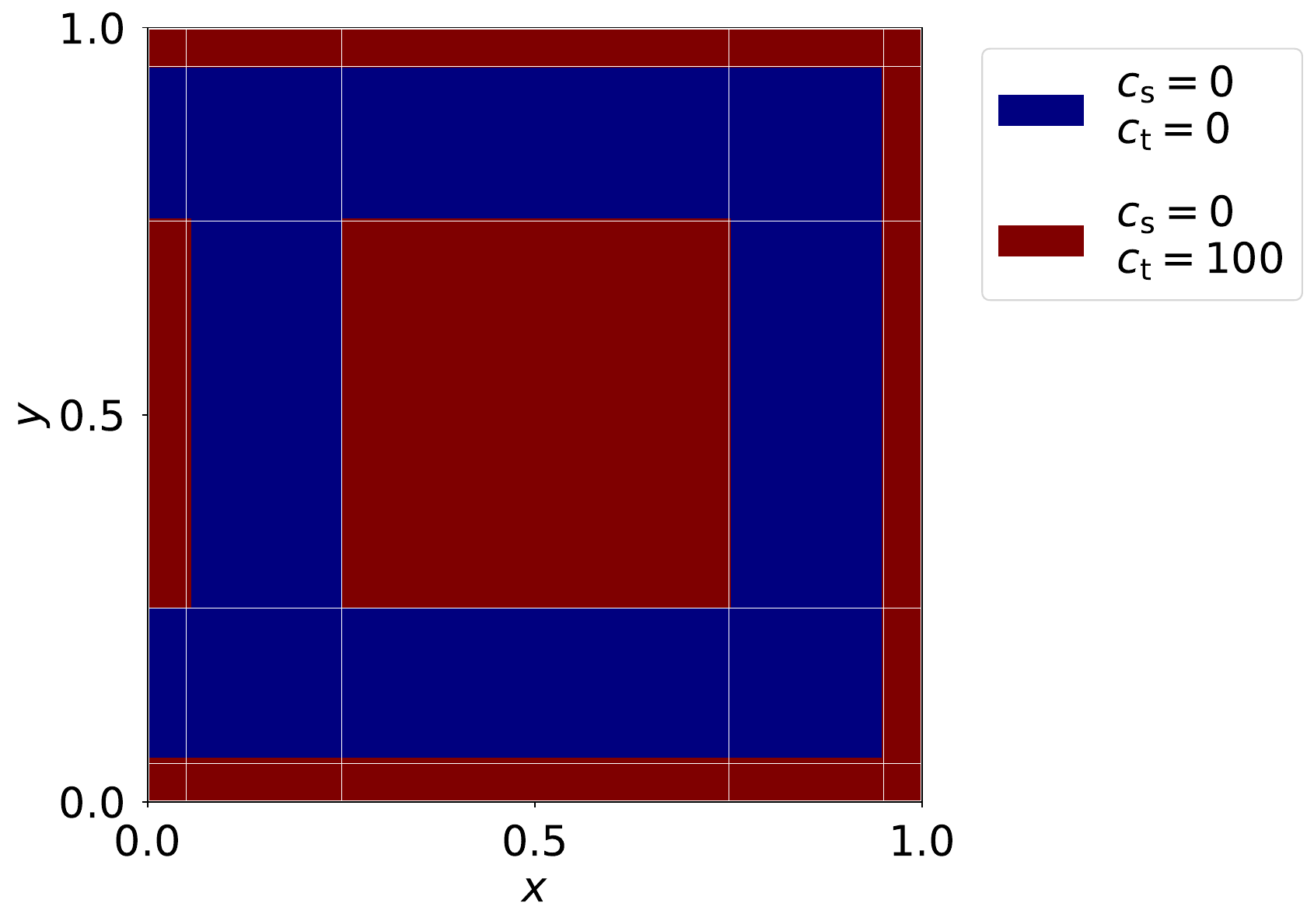}
    \caption{Setup of the Hohlraum and point source test problem.}
    \label{fig::hohlraum_setup}
\end{figure}

The second test case we consider is a hohlraum simulation (see, e.g., \cite{Peng_lowrank_testproblem,Brunner_testproblem,McClarren_testproblem}). 
The motivation for this test case comes from inertial confinement fusion, where the hohlraum holds a small capsule of fusion fuel. In our case, we look at a 
simplified 2-dimensional version of a hohlraum. As can be seen in figure \ref{fig::hohlraum_setup}, we have a rectangular domain 1 cm wide in the $x$ and $y$ directions.
The angular domain is again $\phi \in [0,2\pi]$. The advection coefficient is equal to 1 on the whole domain. The blue blocks describe a vacuum, where $c_{\text{t}}=0$ and $c_{\text{s}}=0$, while the red blocks describe a strongly absorbing material with $c_{\text{t}}=100$ and $c_{\text{s}}=0$.

In this problem there is no source, i.e.~$Q(x,y)=0$. Radiation enters via an inflow boundary condition given by
\begin{equation*}
  f(t,x=0,y,\phi) = 1, \text{ for all } t\in \mathbb{R}^+, y\in [0,1] \text{ and } \phi\in \left[0,\frac{\pi}{2}\right]\cup \left[\frac{3\pi}{2},2\pi\right].
\end{equation*}
At all other boundaries in the domain we have either a zero-inflow condition or an outflow boundary. The initial condition is again given by  {\color{black}$f(t=0,x,y,\phi)=0$}.

\begin{figure}[H]
    \centering
    \begin{subfigure}{\textwidth}
        \centering
        \begin{minipage}{0.12\textwidth}
          \raggedleft
          \caption{$t=0.2$:}
        \end{minipage}
        \begin{minipage}{0.424\textwidth}
          \includegraphics[width=\textwidth]{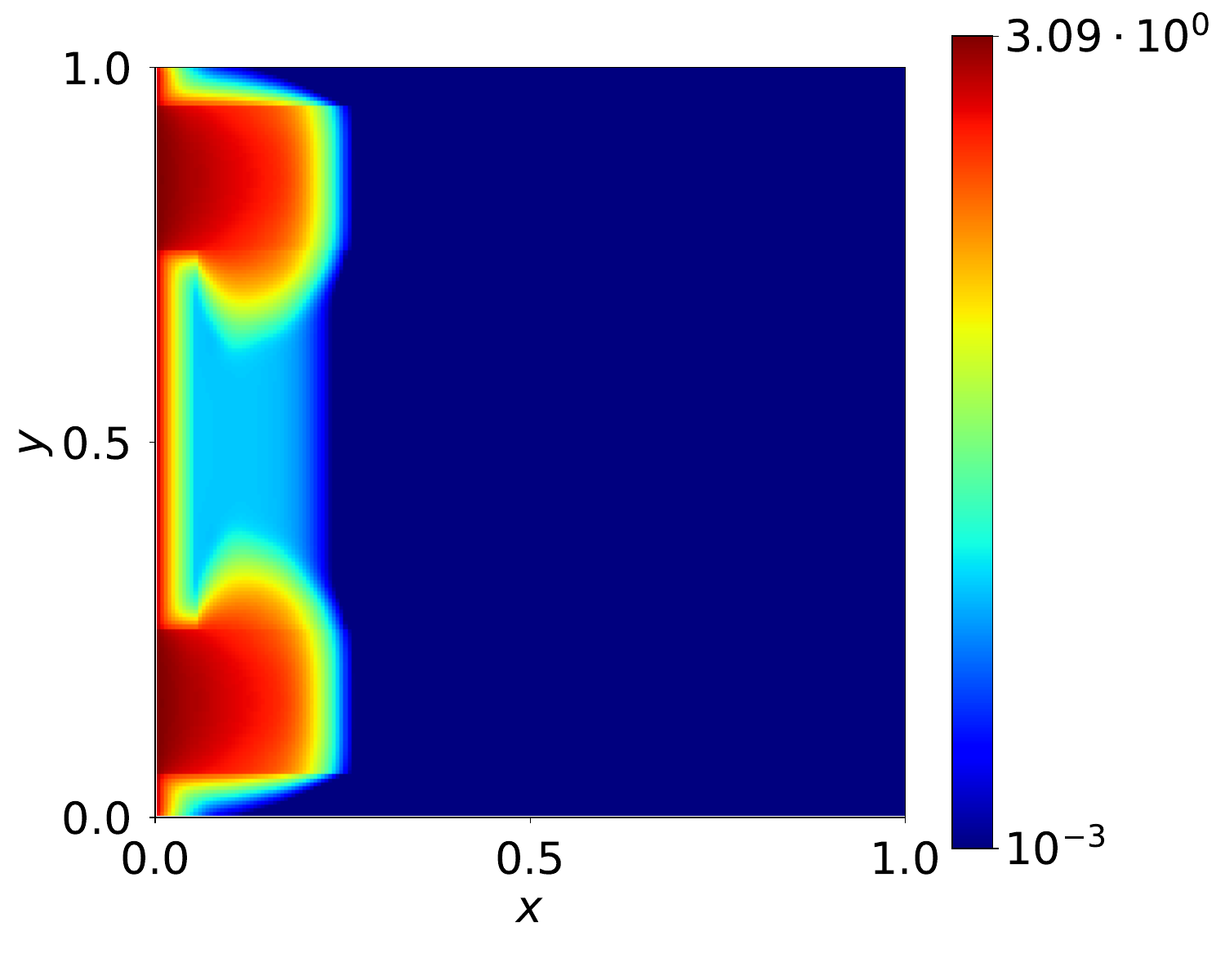}
        \end{minipage}
        \begin{minipage}{0.424\textwidth}
            \centering
            \includegraphics[width=\textwidth]{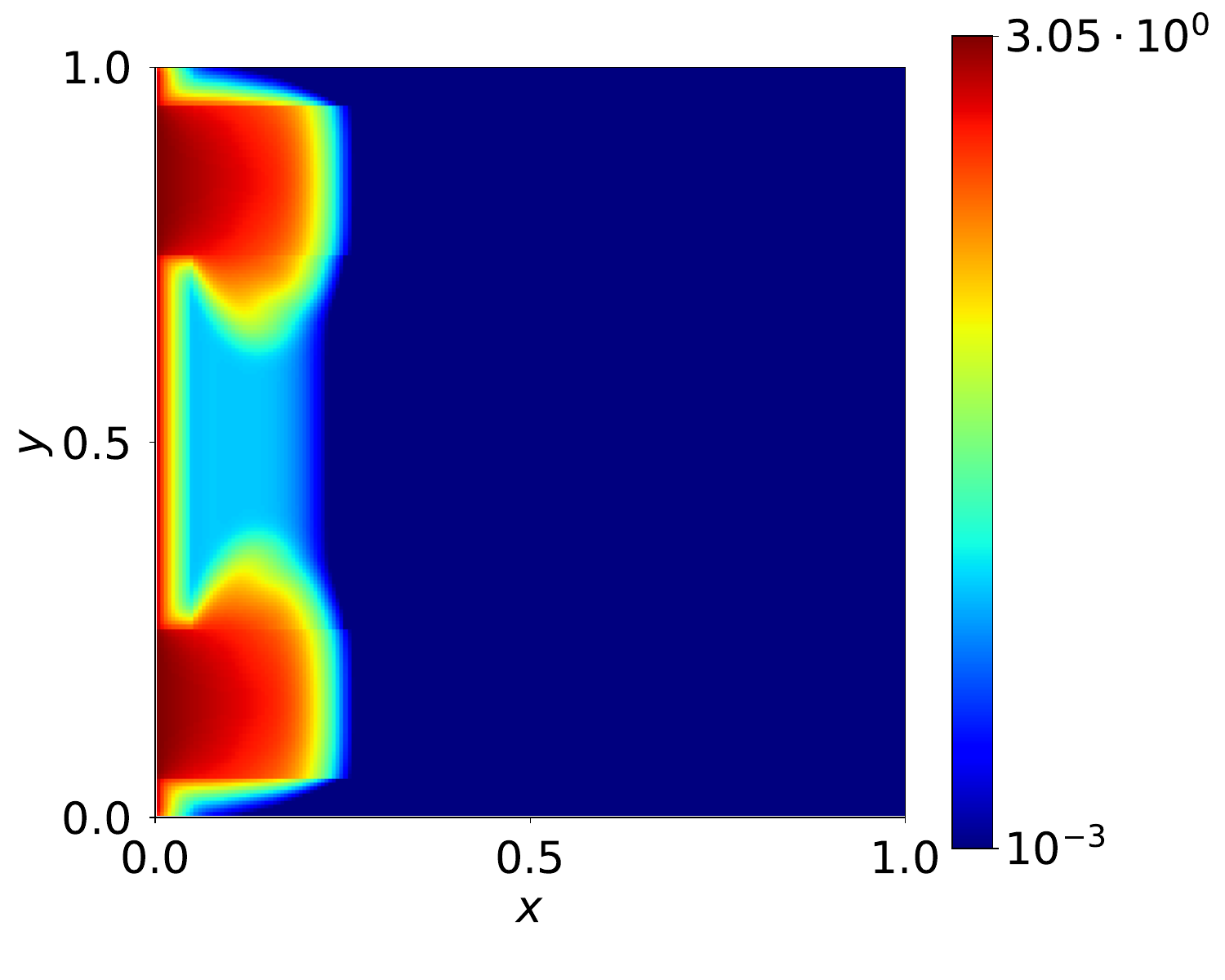}
        \end{minipage}
    \end{subfigure}

    \begin{subfigure}{\textwidth}
        \centering
        \begin{minipage}{0.12\textwidth}
          \raggedleft
          \caption{$t=0.4$:}
        \end{minipage}
        \begin{minipage}{0.424\textwidth}
          \includegraphics[width=\textwidth]{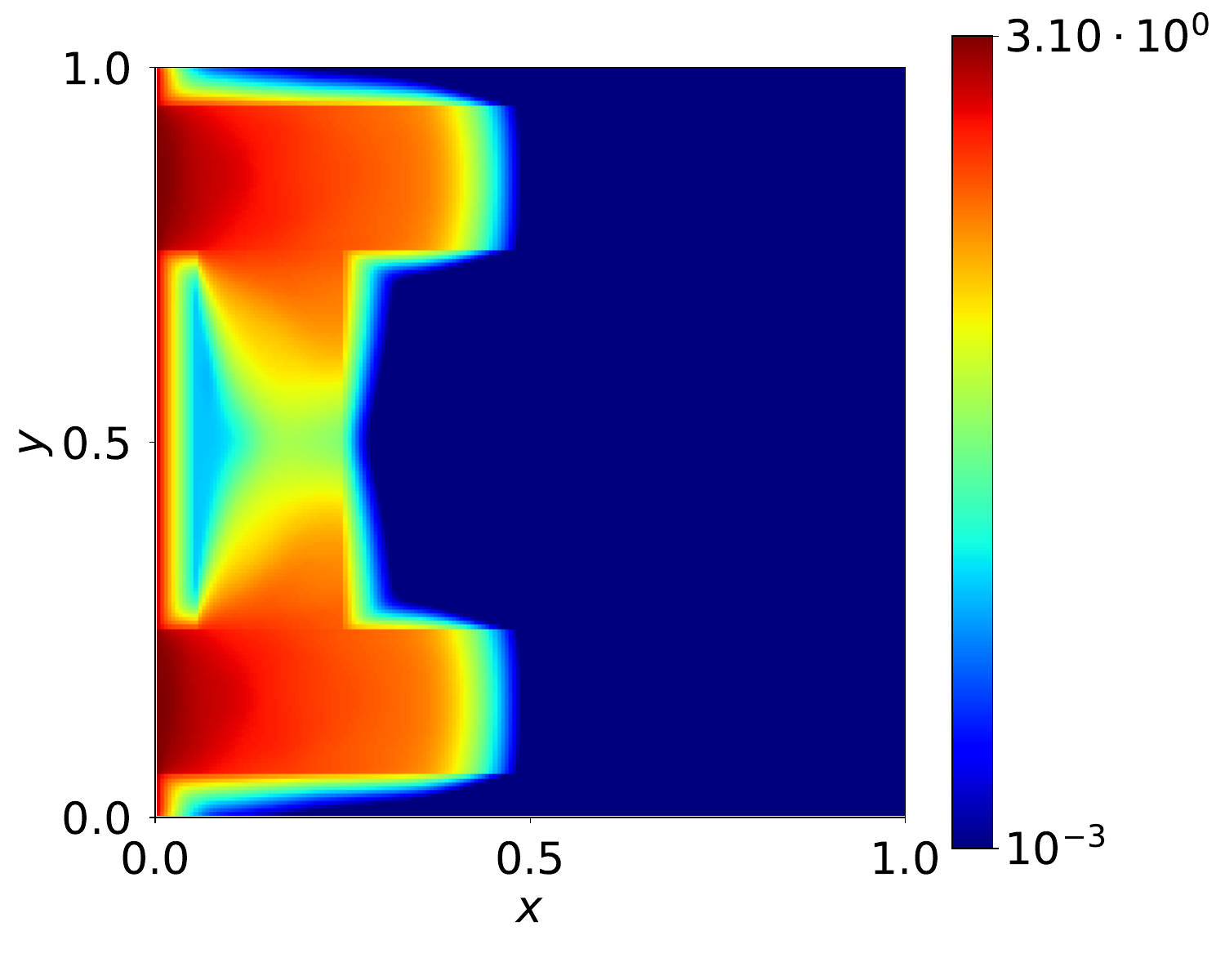}
        \end{minipage}
        \begin{minipage}{0.424\textwidth}
            \centering
            \includegraphics[width=\textwidth]{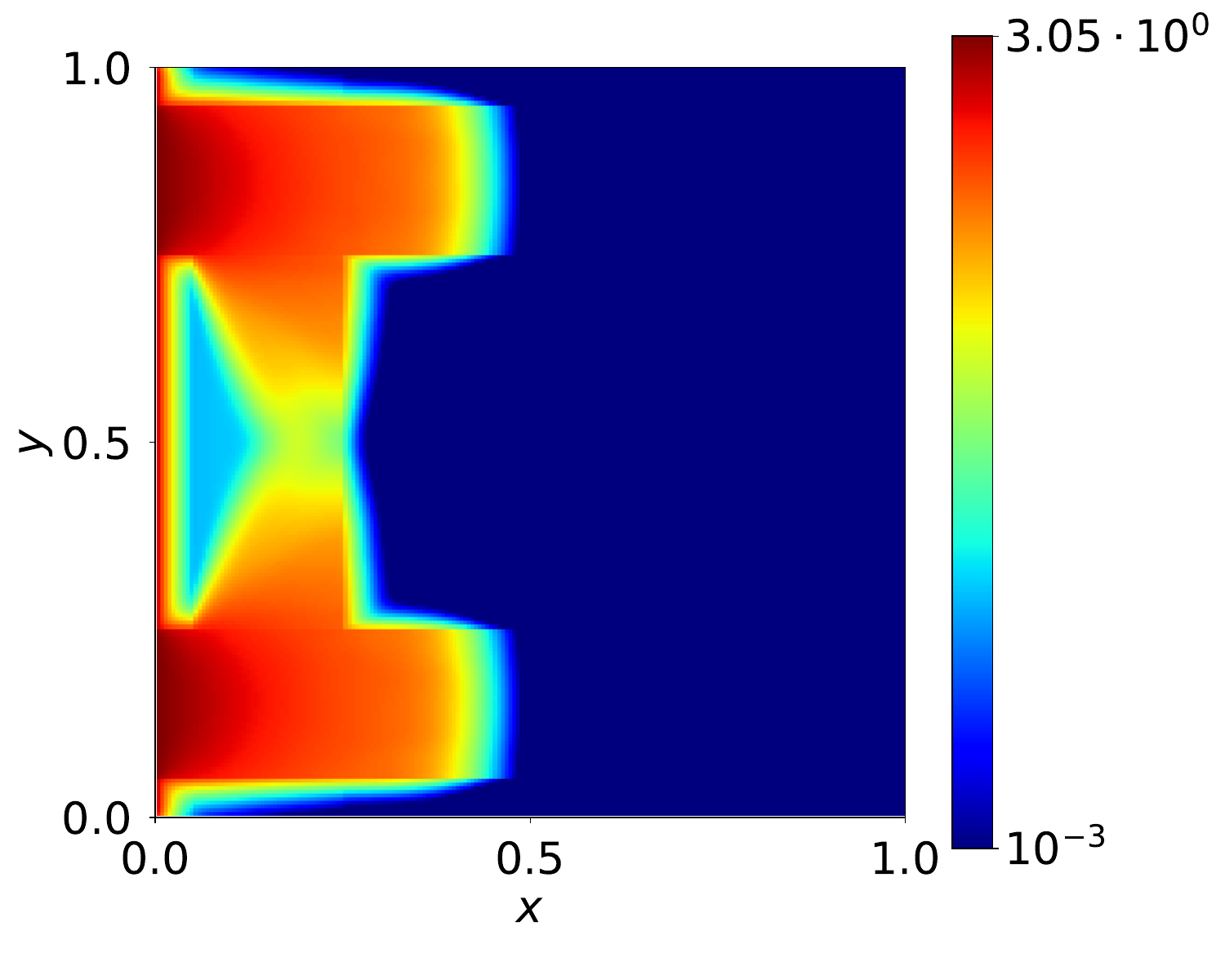}
        \end{minipage}
    \end{subfigure}

    \begin{subfigure}{\textwidth}
        \centering
        \begin{minipage}{0.12\textwidth}
          \raggedleft
          \caption{$t=1.2$:}
        \end{minipage}
        \begin{minipage}{0.424\textwidth}
          \includegraphics[width=\textwidth]{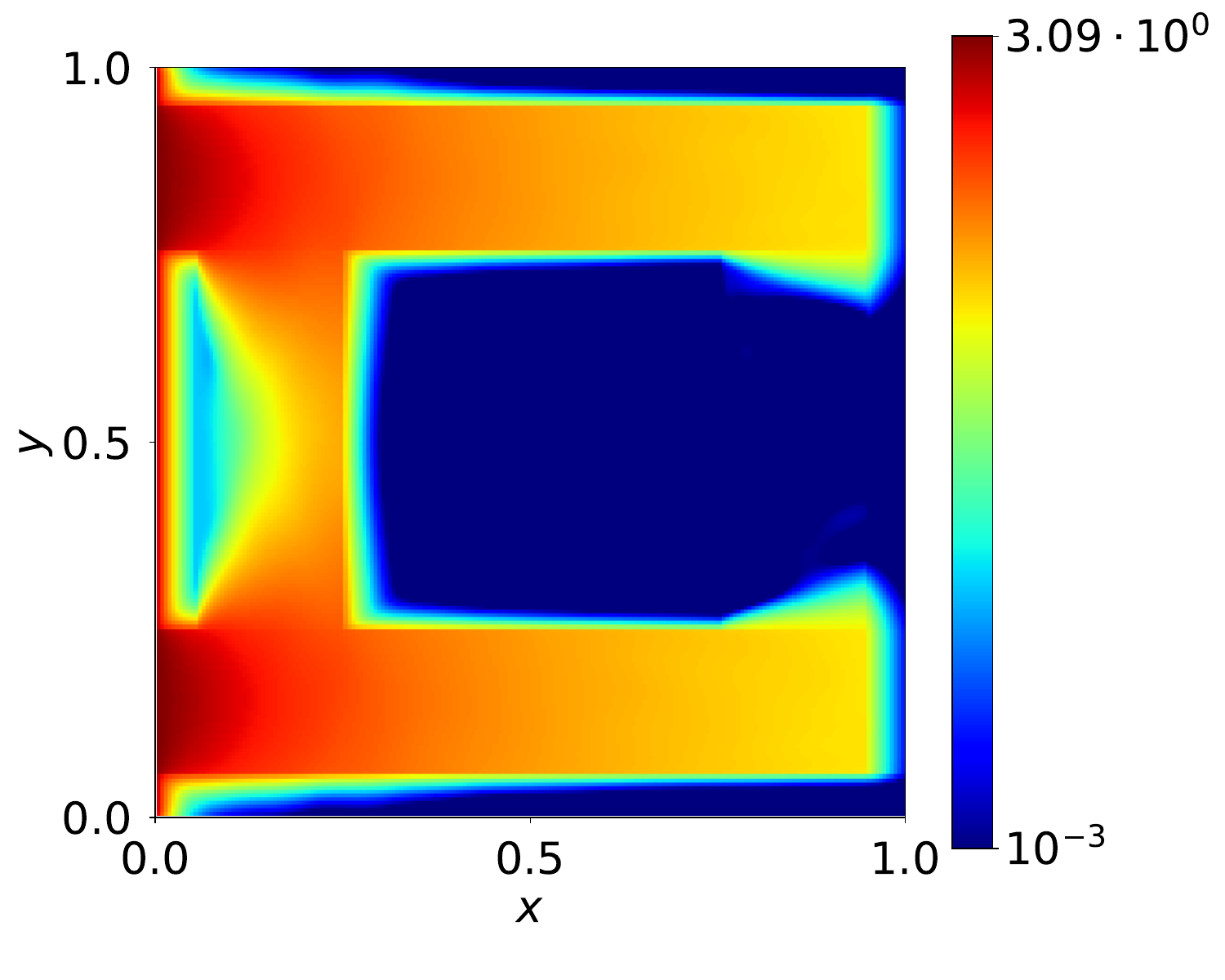}
        \end{minipage}
        \begin{minipage}{0.424\textwidth}
            \centering
            \includegraphics[width=\textwidth]{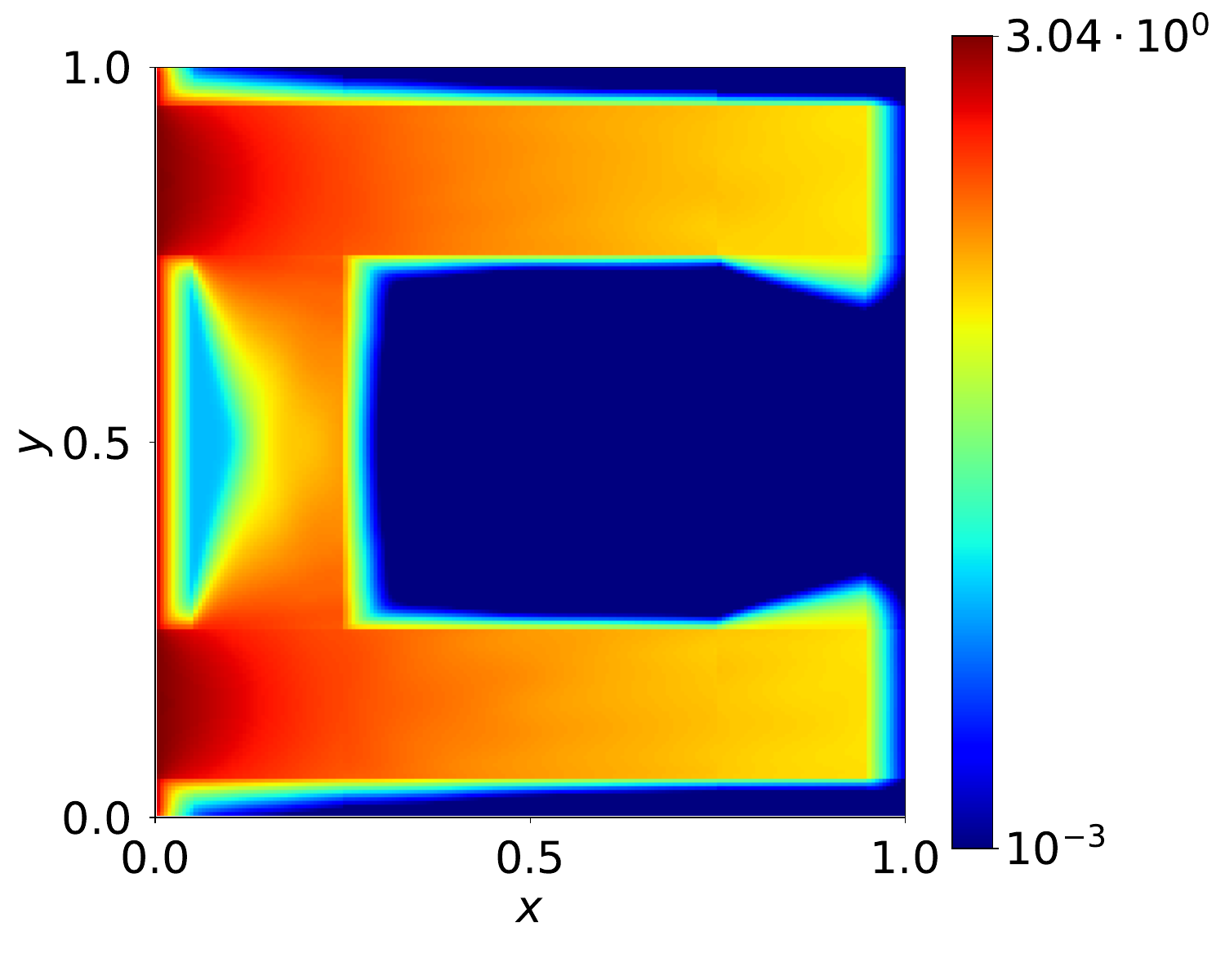}
        \end{minipage}
    \end{subfigure}
    \caption{We show $\rho (t,x,y)$ obtained by solving the Hohlraum test problem with the classic dynamical low-rank algorithm (left) and with the domain decomposition low-rank algorithm \ref{alg::ddlr_algorithm} (right). For the latter, the grid is divided into $5 \times 5$ subdomains.}
    \label{fig::hohlraum_rho_plots}
\end{figure}

We again compare the low-rank methods with and without domain decomposition. An equidistant grid with 200 grid points in $x$, $y$ and $\phi$ is used, and the time step size is set to $\Delta t = 0.5 \cdot \Delta x$. The results for $\rho(t,x,y)$ at times $t=0.2$, $t=0.4$ and $t=1.2$ are shown in figure~\ref{fig::hohlraum_rho_plots} (left without domain decomposition and right with domain decomposition). We observe that the radiation enters the domain on the left boundary and moves through the vacuum. The absorbing blocks prevent radiation to move through the material and instead force it around. Similar results are observed for the low-rank method with and without domain decomposition.

\begin{figure}[H]
    \centering
    \begin{subfigure}{0.44\textwidth}
        \centering
        \includegraphics[width=\textwidth]{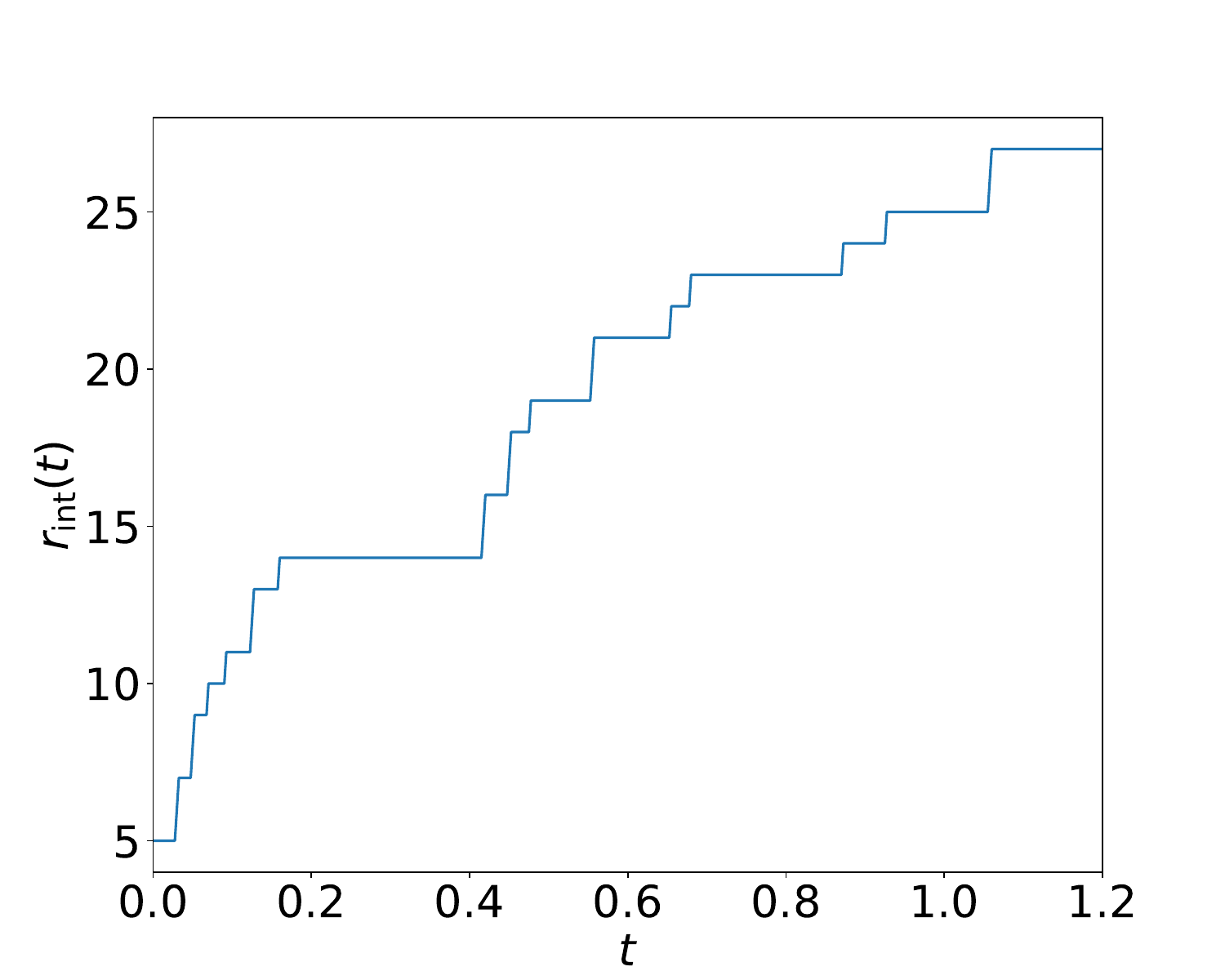}
        \caption{Intermediate ranks}
    \end{subfigure}
    \begin{subfigure}{0.44\textwidth}
        \centering
        \includegraphics[width=\textwidth]{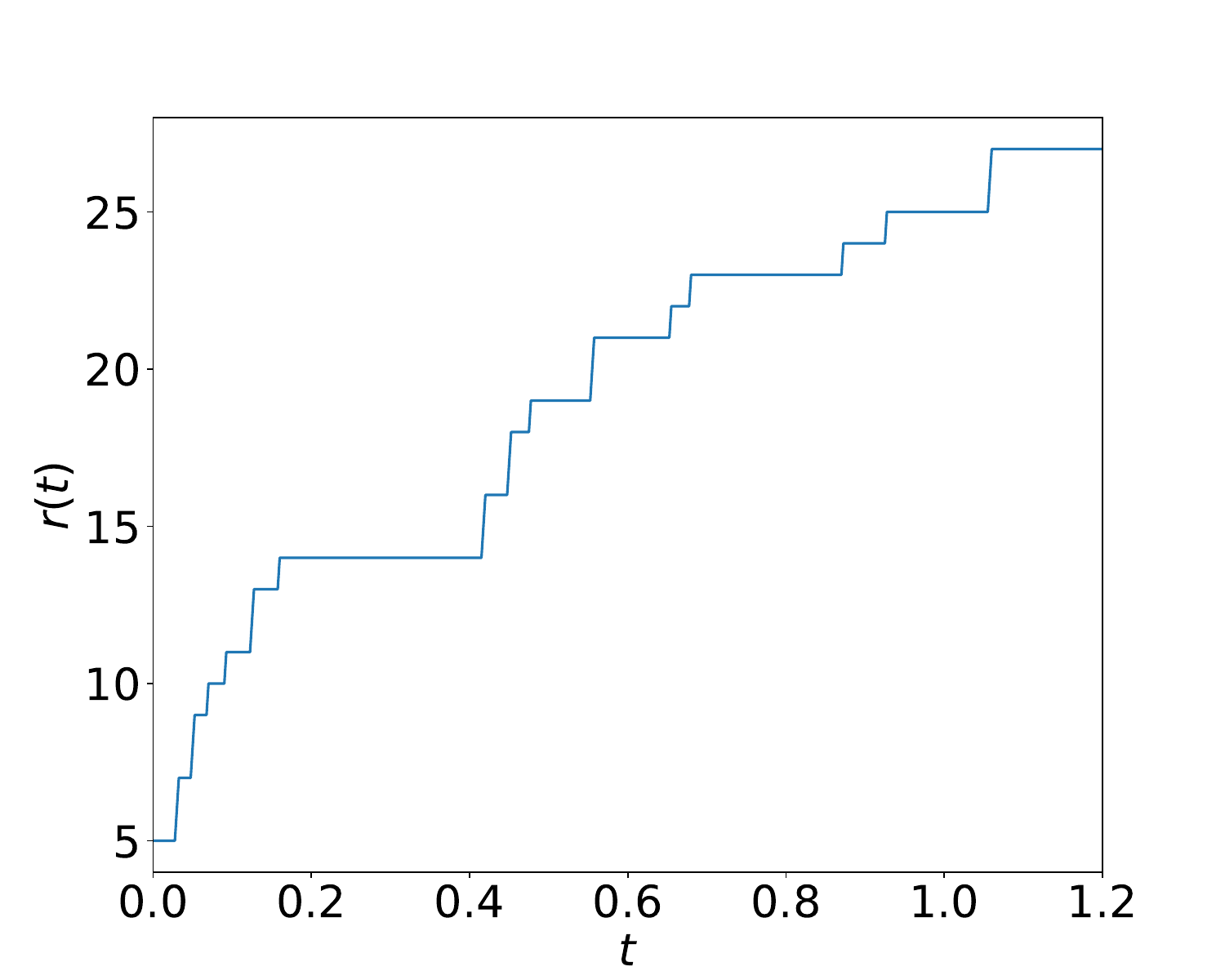}
        \caption{Ranks}
    \end{subfigure}
    \caption{We show the time evolution of the intermediate rank and rank after the time step for the classic low-rank simulation of the Hohlraum test problem.  As an adaptive rank scheme we used the scheme as proposed in algorithms \ref{alg:augmentation} and \ref{alg:truncate} with {\color{black} $\text{tol}=1.3\cdot 10^{-4}$}.}
    \label{fig::hohlraum_1domain_ranks}
\end{figure}

\begin{figure}[H]
    \centering
    \begin{subfigure}{0.44\textwidth}
        \centering
        \includegraphics[width=\textwidth]{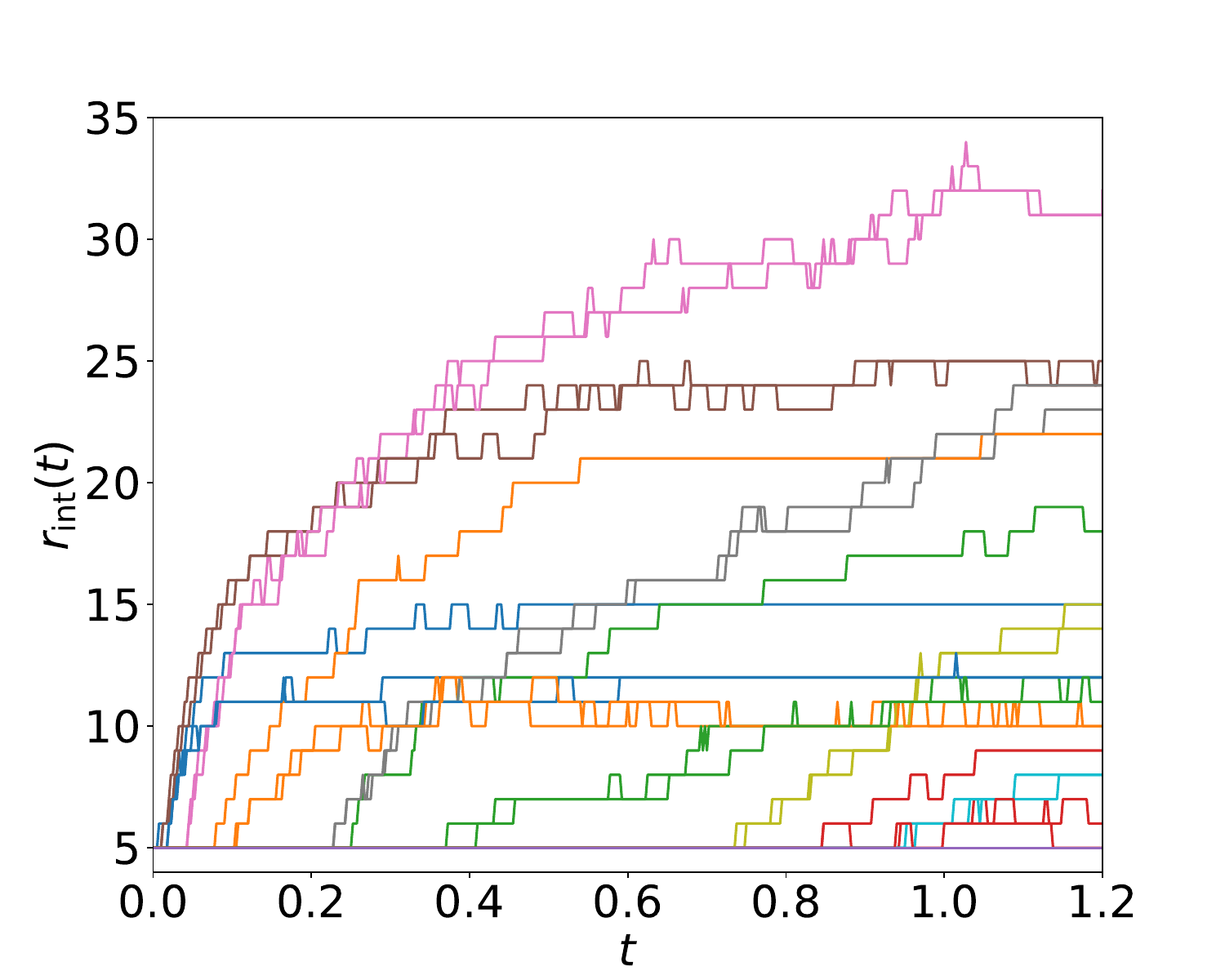}
        \caption{Intermediate ranks}
    \end{subfigure}
    \begin{subfigure}{0.46\textwidth}
        \centering
        \includegraphics[width=\textwidth]{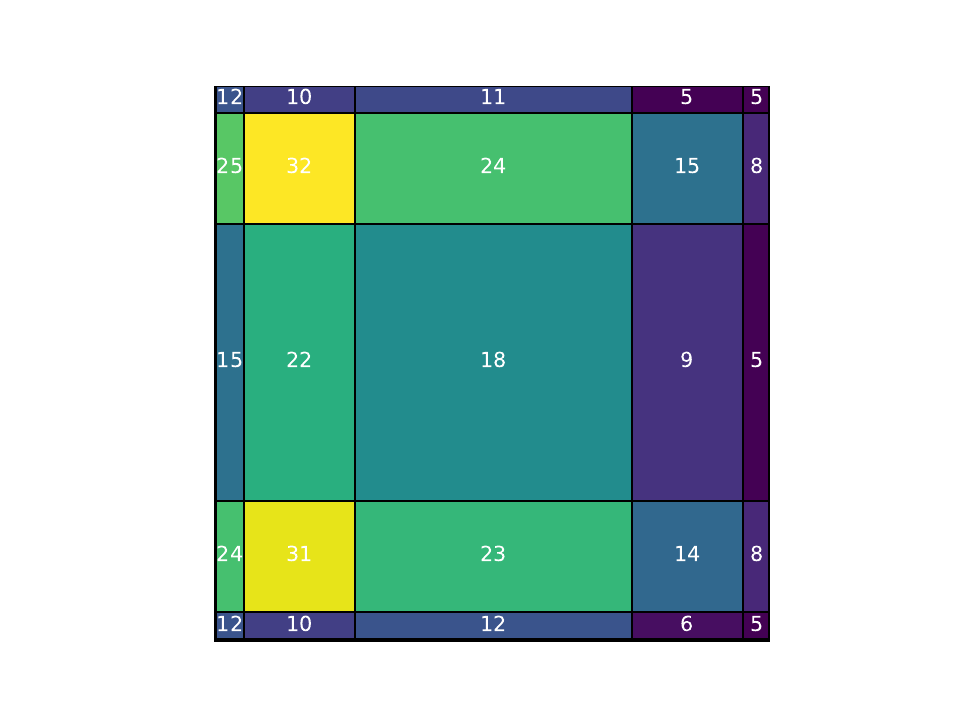}
        \caption{Intermediate ranks at time $t=1.2$}
    \end{subfigure}
    \begin{subfigure}{0.44\textwidth}
        \centering
        \includegraphics[width=\textwidth]{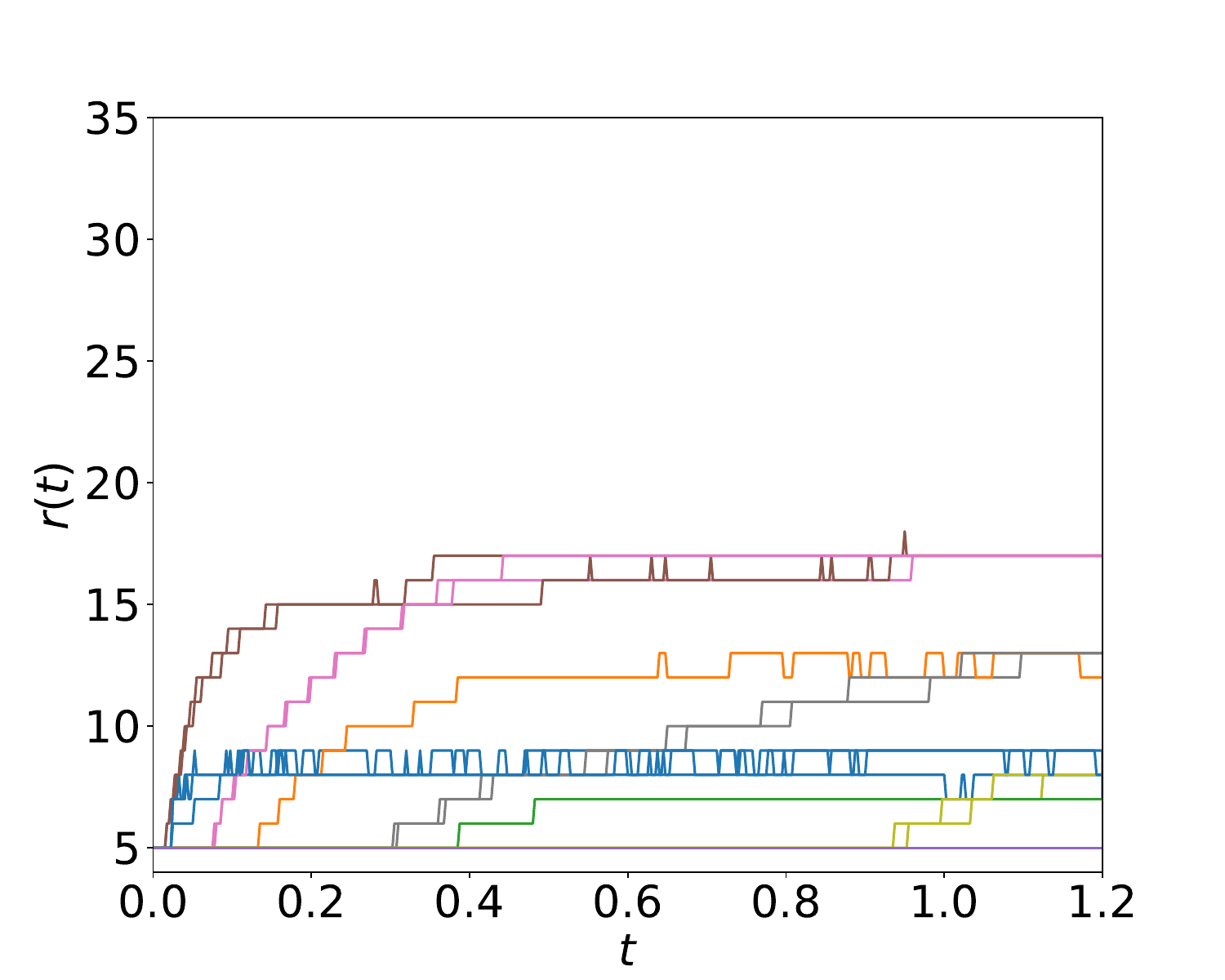}
        \caption{Ranks}
    \end{subfigure}
    \begin{subfigure}{0.46\textwidth}
        \centering
        \includegraphics[width=\textwidth]{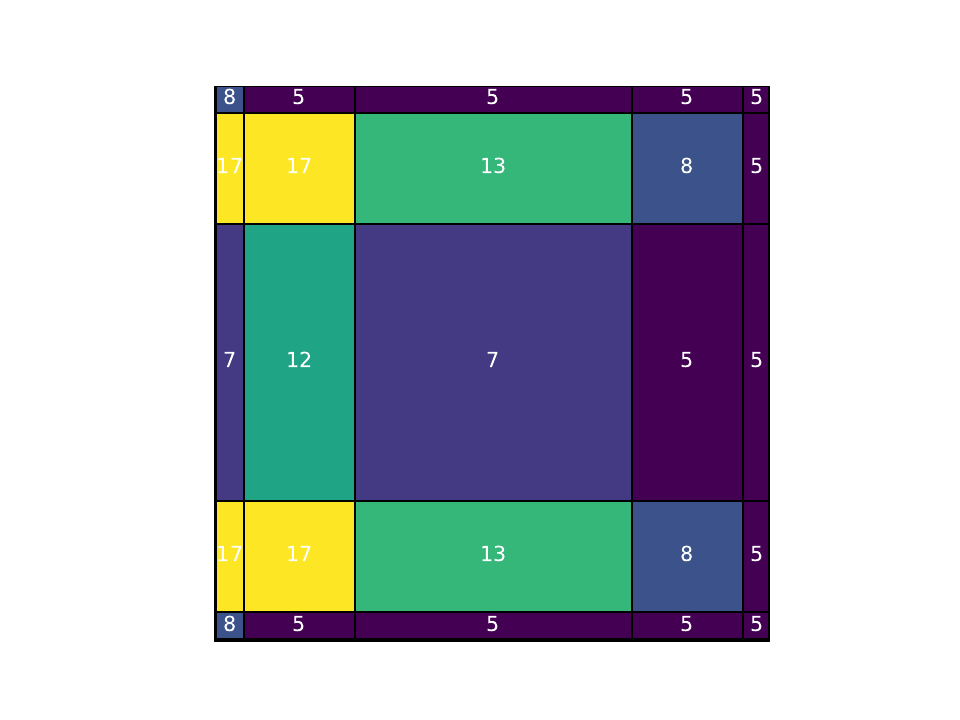}
        \caption{Ranks at time $t=1.2$}
    \end{subfigure}
    \caption{We show the time evolution of the intermediate ranks and ranks after the time step for the domain decomposition dynamical low-rank simulation of the Hohlraum test problem. The adaptive rank scheme as proposed in algorithms \ref{alg:augmentation} and \ref{alg:truncate} with $\text{tol}=10^{-3}$ was used.}
    \label{fig::hohlraum_domaindecomp_ranks}
\end{figure}

We again compare the time evolution of the ranks for the approach without (figure \ref{fig::hohlraum_1domain_ranks}) and with domain decomposition (figure \ref{fig::hohlraum_domaindecomp_ranks}). {\color{black}We note that we have chosen a tolerance of $\text{tol} = 1.3 \cdot 10^{-4}$ for the classic dynamical low-rank approach. This is necessary to avoid visible artifacts in the numerical solution. For the domain decomposition low-rank approach a tolerance of $\text{tol}=10^{-3}$ is sufficient to obtain good results.} We note, however, that the largest intermediate rank of the domain decomposition low-rank algorithm is roughly the same in both cases as is the error compared to a reference solution computed with the classic method and $\text{tol}=10^{-10}$ ($\norm{f-f_{\text{ref}}} \approx 1.3 \cdot 10^{-2}$ for the classic DLRA approach and $\norm{f-f_{\text{ref}}} \approx 2.5 \cdot 10^{-2}$ for the domain decomposition). We believe that the slightly larger error in case of the domain decomposition is actually not a low-rank error, but a time discretization error caused by freezing the boundary condition in time. {\color{black}One could reduce this error by extending our time integration at the boundary to second order (e.g. doing half a step of an evolution equation as a predictor, then a full step as a corrector), which is left for future work.}

Comparing the degrees of freedom required by each method in figure \ref{fig::lattice_dof_comparison} on the right we see that the domain decomposition approach requires roughly {\color{black} 30} \% less degrees of freedom when comparing maximum memory usage (i.e.~we compare with the intermediate rank) and more than a factor of two less memory to store the solution (i.e.~we compare with the rank obtained after the time step). This is despite the fact that both methods have a similar maximal rank. However, for the domain decomposition approach a significantly smaller rank is used in most subdomains.

\subsection{Isotropic point source test problem}\label{sec::iso_pointsource}

The {\color{black} next} test case we consider is a simulation with an isotropic point source. It is well known that an isotropic point source propagating into vacuum can only be approximated well by using a large rank. We use the same setup as in the Hohlraum example from section \ref{sec::hohlraum}, but instead of having a uniform inflow on the entire left boundary, we now place a point source on the left boundary, i.e.
\begin{equation*}
  \begin{aligned}
  f(t,x=0,y,\phi) = \frac{1}{\sqrt{2\pi}\sigma} \exp\left(-\frac{(y-y_0)^2}{2\sigma^2}\right), &\text{ for all } t\in \mathbb{R}^+, y\in [0,1] \\
  &\text{ and } \phi\in \left[0,\frac{\pi}{2}\right]\cup \left[\frac{3\pi}{2},2\pi\right],
  \end{aligned}
\end{equation*}
where $y_0 = 0.85$ and a small $\sigma = 10^{-2}$. To resolve the point source accurately, we require $600$ grid points in both the $x$ and $y$ directions. We still use 200 grid points in $\phi$.

\begin{figure}[H]
    \centering
    \begin{subfigure}{\textwidth}
        \centering
        \begin{minipage}{0.12\textwidth}
          \raggedleft
          \caption{$t=0.3$:}
        \end{minipage}
        \begin{minipage}{0.424\textwidth}
          \includegraphics[width=\textwidth]{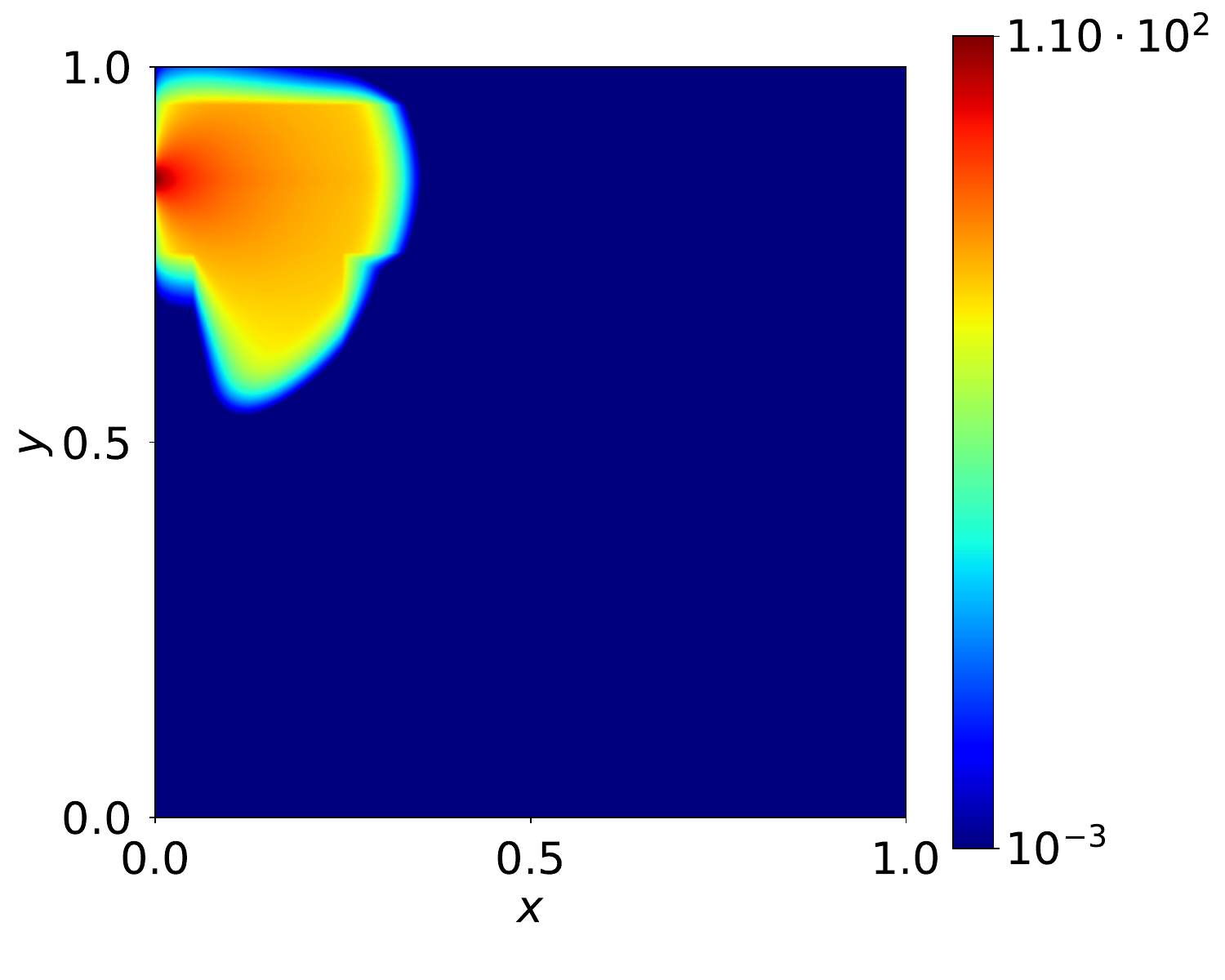}
        \end{minipage}
        \begin{minipage}{0.424\textwidth}
            \centering
            \includegraphics[width=\textwidth]{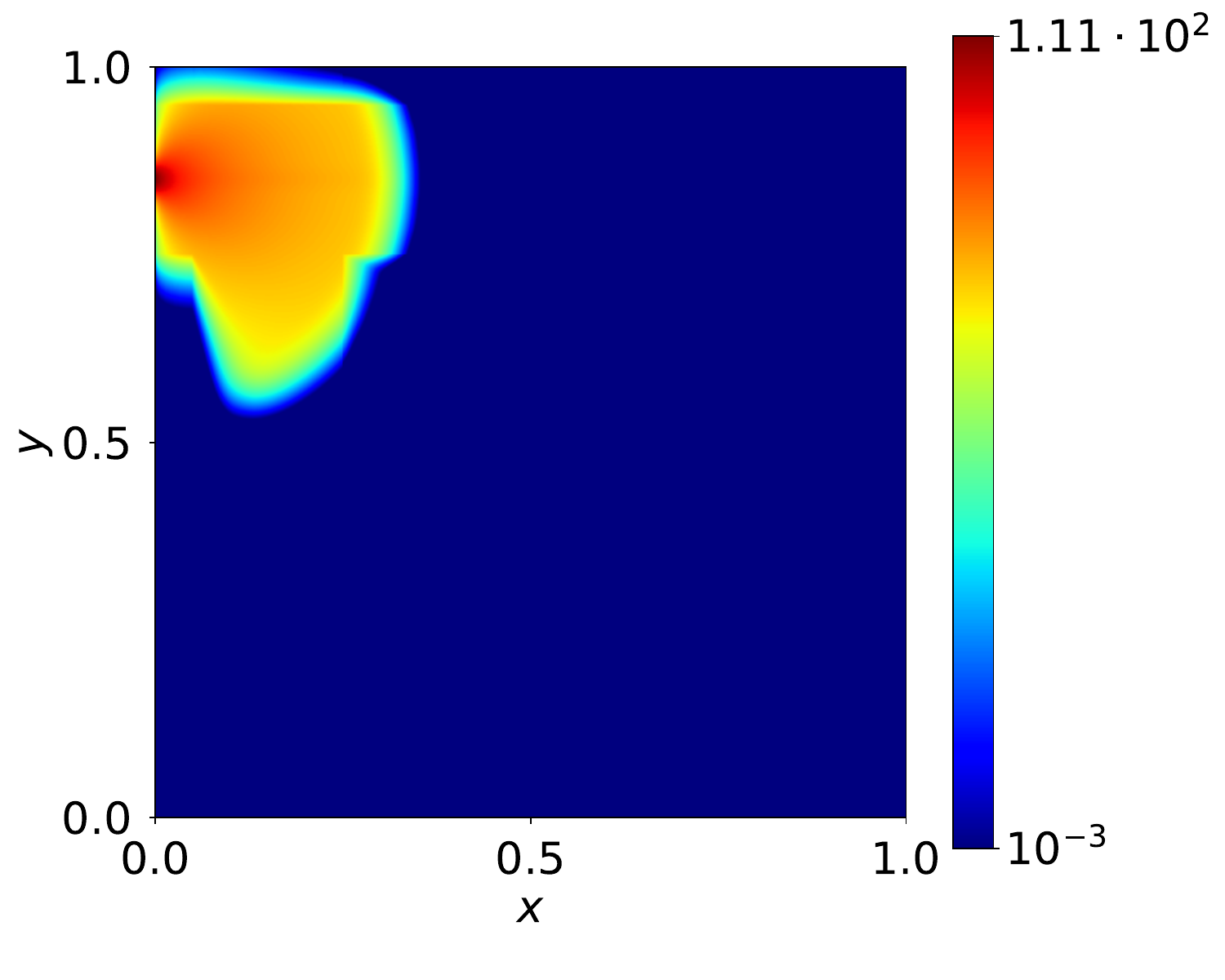}
        \end{minipage}
    \end{subfigure}
    
    \begin{subfigure}{\textwidth}
        \centering
        \begin{minipage}{0.12\textwidth}
          \raggedleft
          \caption{$t=0.7$:}
        \end{minipage}
        \begin{minipage}{0.424\textwidth}
          \includegraphics[width=\textwidth]{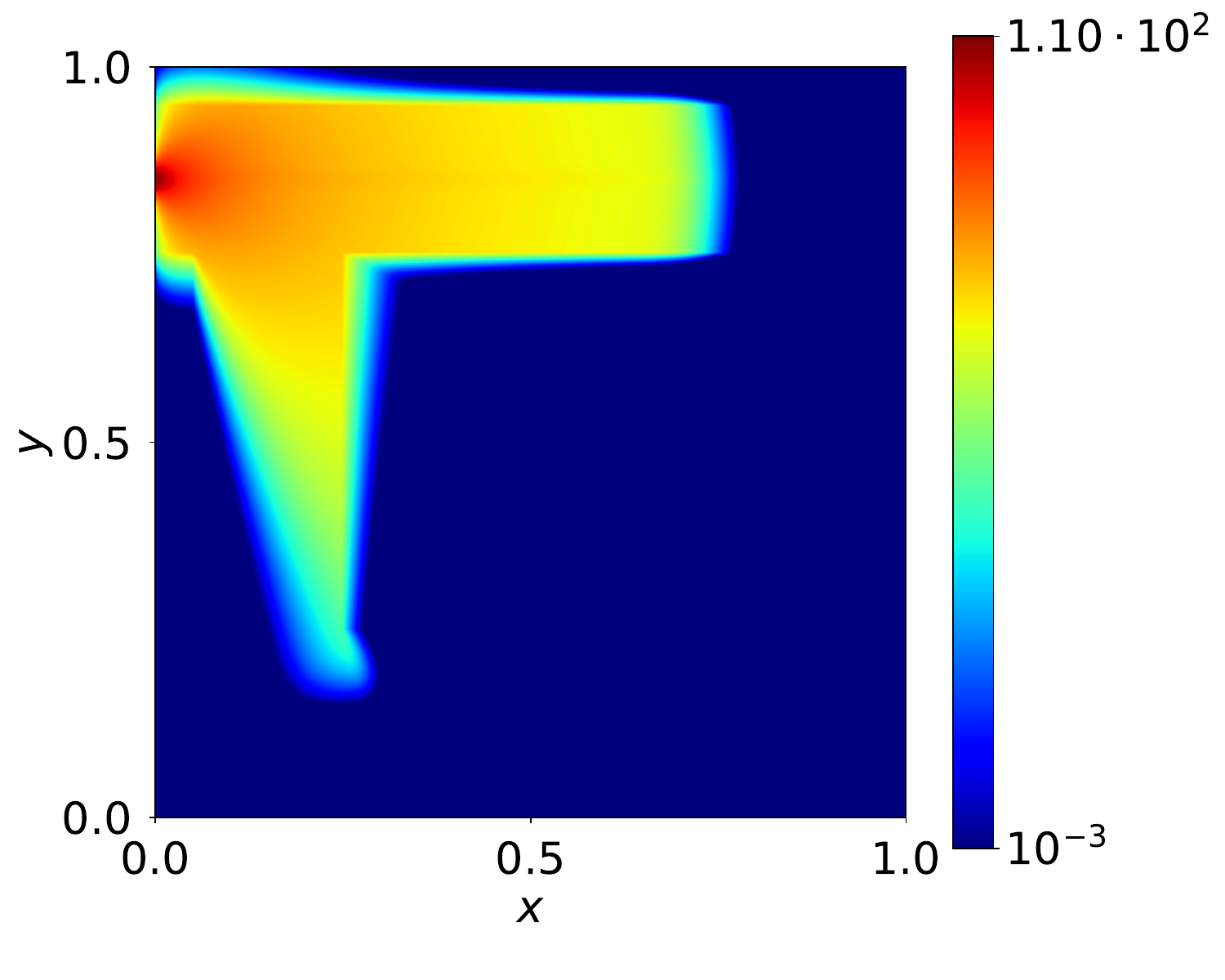}
        \end{minipage}
        \begin{minipage}{0.424\textwidth}
            \centering
            \includegraphics[width=\textwidth]{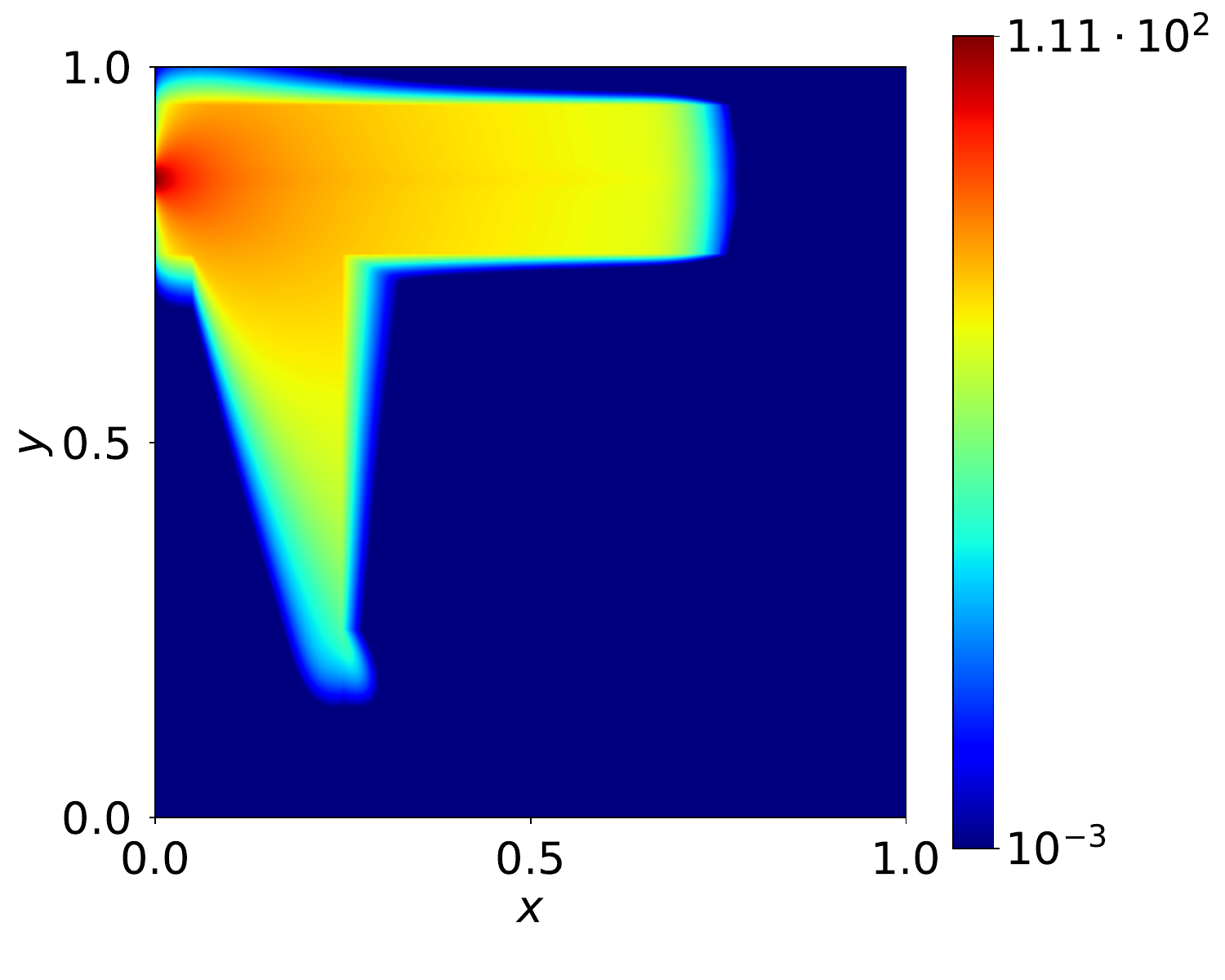}
        \end{minipage}
    \end{subfigure}
    
    \begin{subfigure}{\textwidth}
        \centering
        \begin{minipage}{0.12\textwidth}
          \raggedleft
          \caption{$t=1.0$:}
        \end{minipage}%
        \begin{minipage}{0.424\textwidth}
          \includegraphics[width=\textwidth]{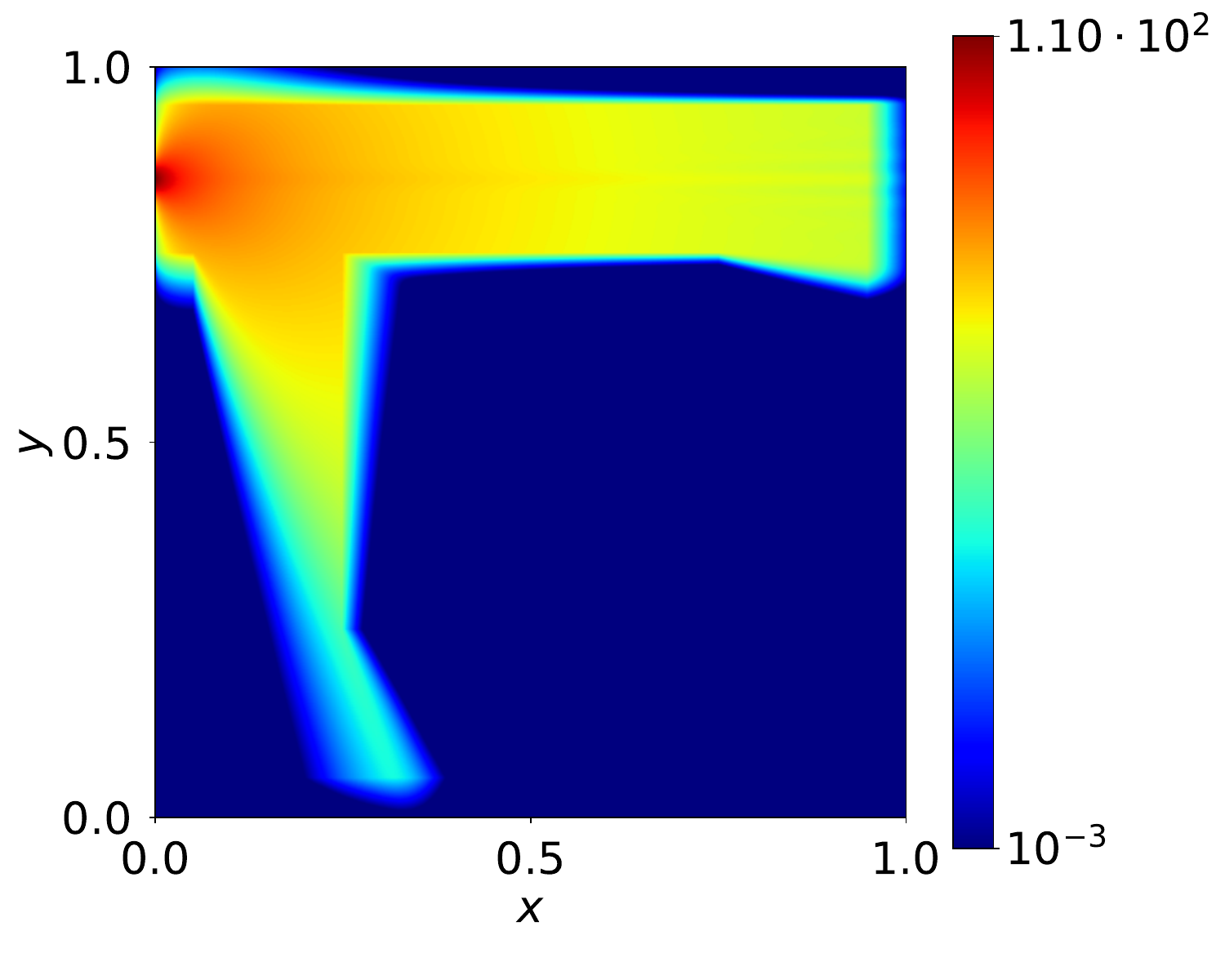}
        \end{minipage}
        \begin{minipage}{0.424\textwidth}
            \centering
            \includegraphics[width=\textwidth]{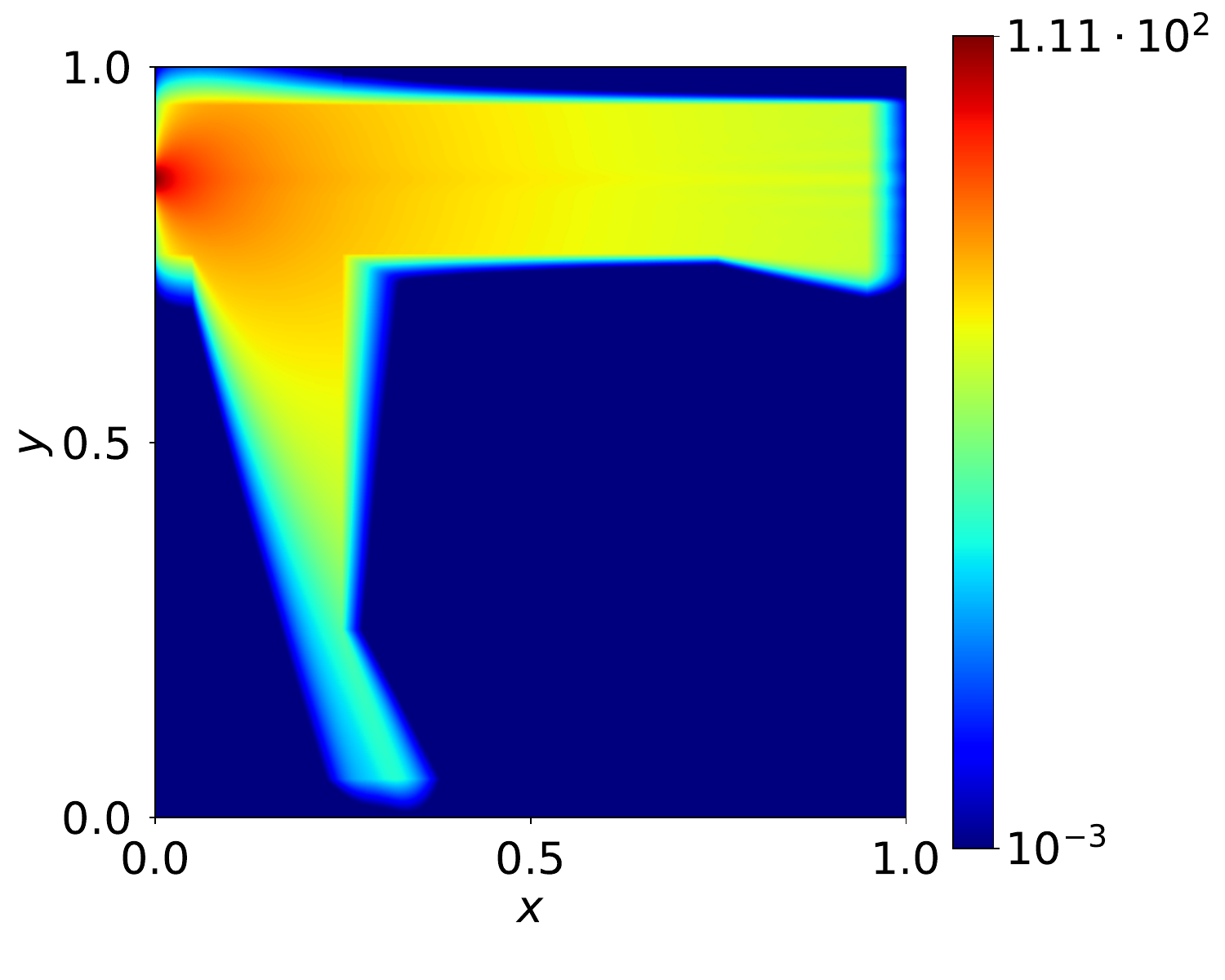}
        \end{minipage}
    \end{subfigure}
    \caption{We show $\rho (t,x,y)$ obtained by solving the point source test problem with the classic dynamical low-rank algorithm (left) and with the domain decomposition low-rank algorithm \ref{alg::ddlr_algorithm} (right). For the latter, the grid is divided into $5 \times 5$ subdomains.}
    \label{fig::pointsource_rho_plots}
\end{figure}

The results for the classic low-rank method and the proposed domain decomposition low-rank algorithm \ref{alg::ddlr_algorithm} for $\rho(t,x,y)$ at times $t=0.3$, $t=0.7$, and $t=1.0$ are shown in figure \ref{fig::pointsource_rho_plots}. As can be clearly observed, the radiation enters our domain on the upper left side of the boundary and propagates into the domain. {\color{black}We chose our tolerances, such that we get identical results for both simulations.}

\begin{figure}[H]
    \centering
    \begin{subfigure}{0.41\textwidth}
        \centering
        \includegraphics[width=\textwidth]{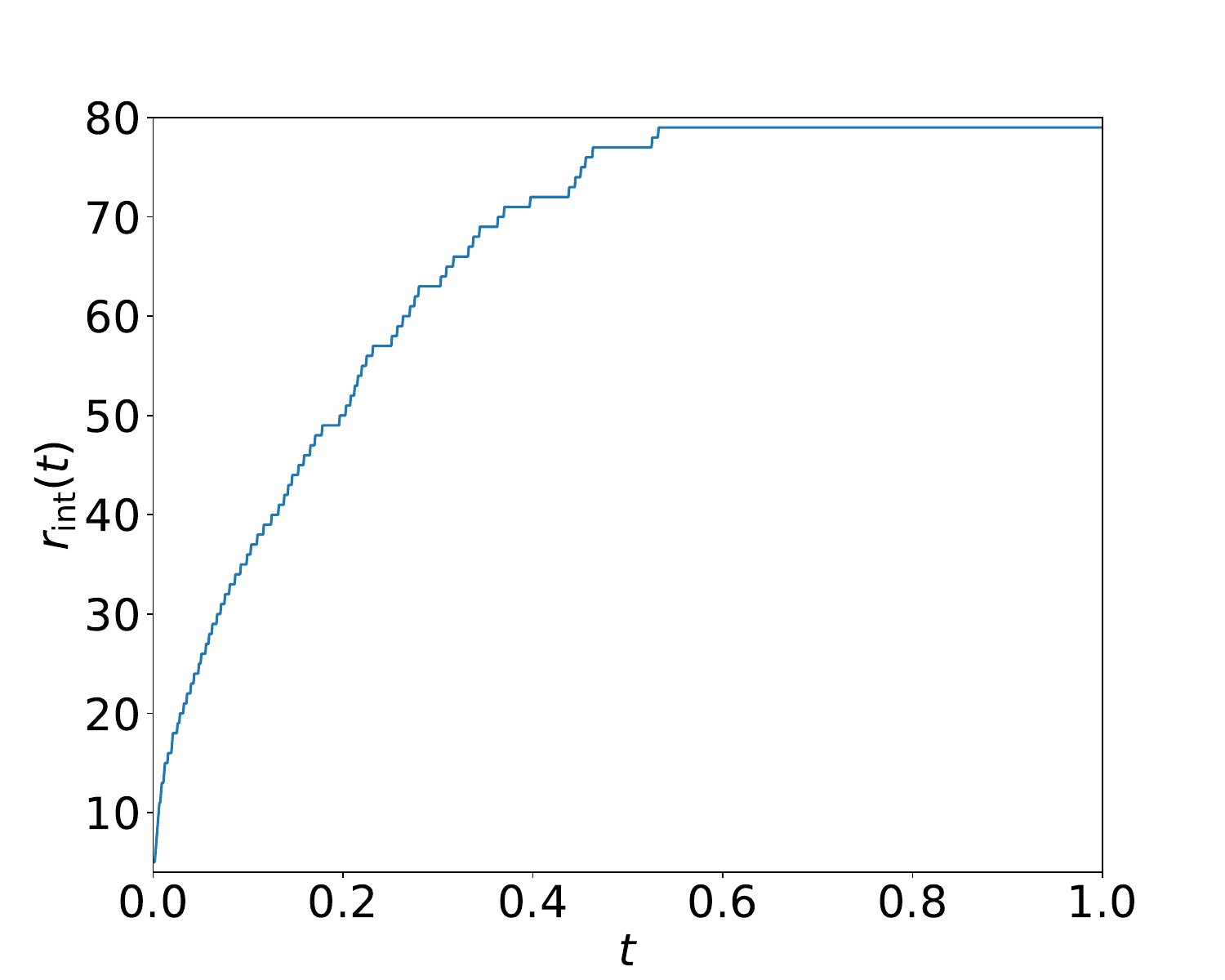}
        \caption{Intermediate ranks}
    \end{subfigure}
    \begin{subfigure}{0.41\textwidth}
        \centering
        \includegraphics[width=\textwidth]{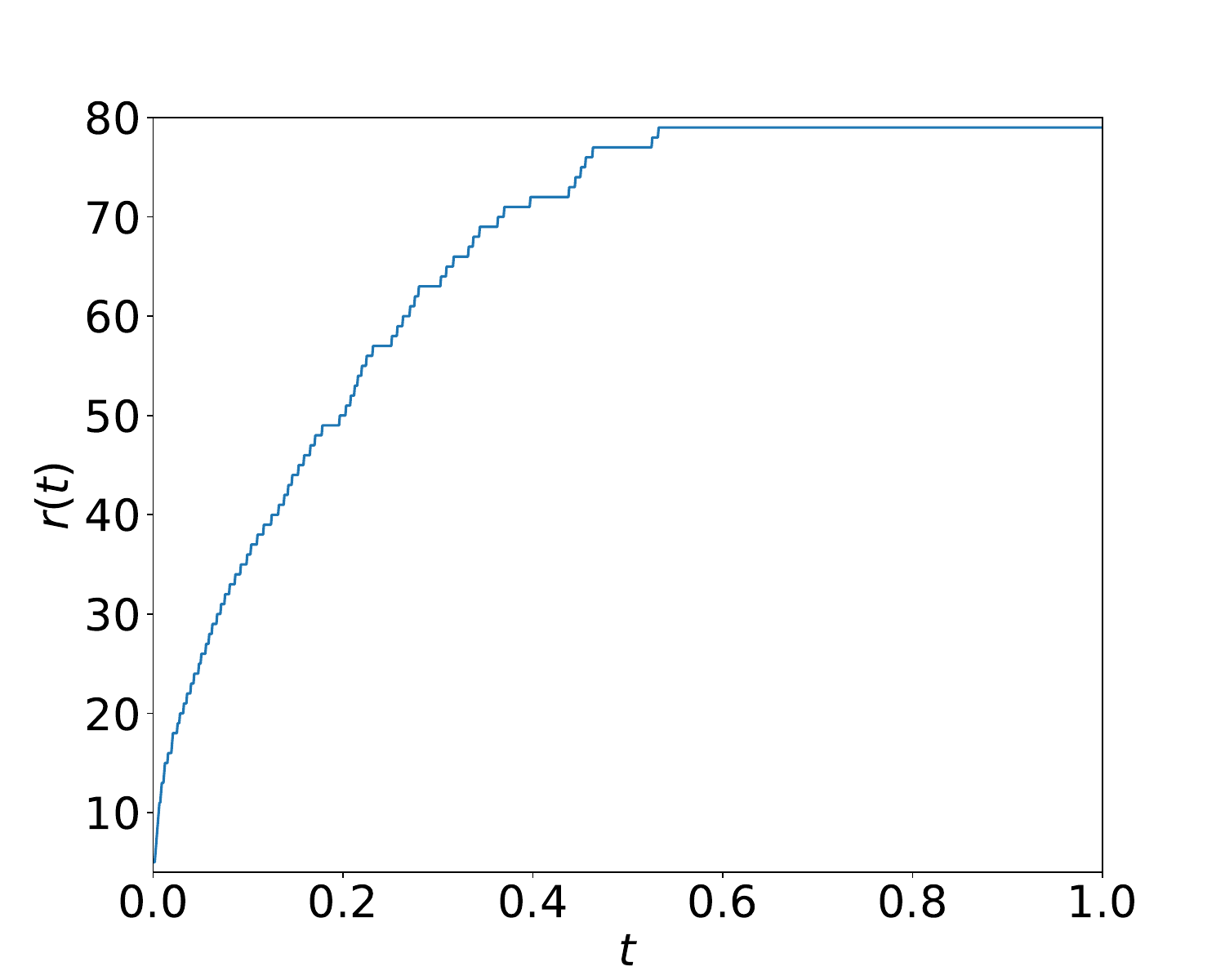}
        \caption{Ranks}
    \end{subfigure}
    \caption{We show the time evolution of the intermediate rank and rank after the time step for the classic low-rank simulation of the point source test problem.  As an adaptive rank scheme we used the scheme as proposed in algorithms \ref{alg:augmentation} and \ref{alg:truncate} with {\color{black} $\text{tol}=1.2\cdot 10^{-5}$}.}
    \label{fig::pointsource_classic_ranks}
\end{figure}

\begin{figure}[H]
    \centering
    \begin{subfigure}{0.44\textwidth}
        \centering
        \includegraphics[width=\textwidth]{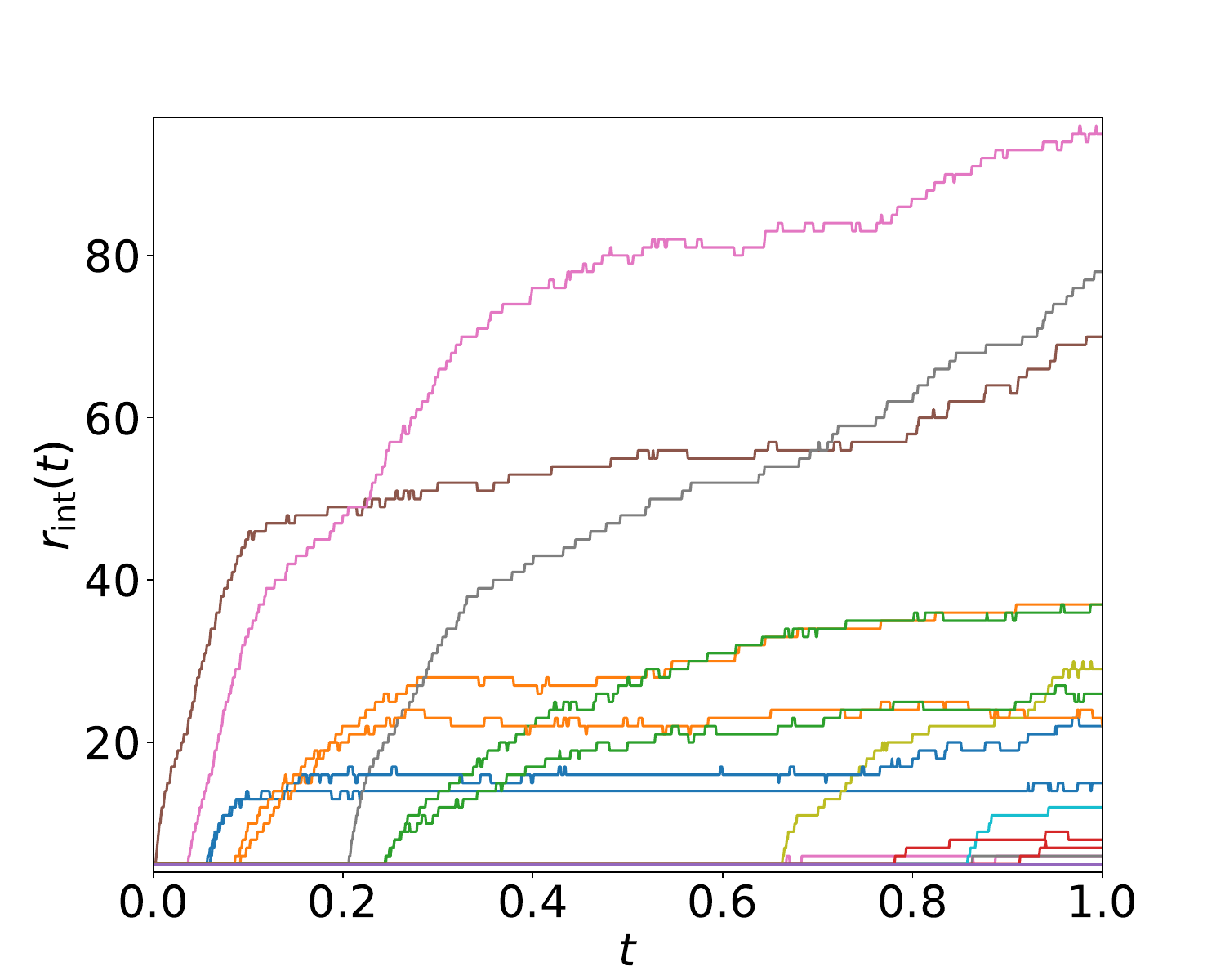}
        \caption{Intermediate ranks}
    \end{subfigure}
    \begin{subfigure}{0.46\textwidth}
        \centering
        \includegraphics[width=\textwidth]{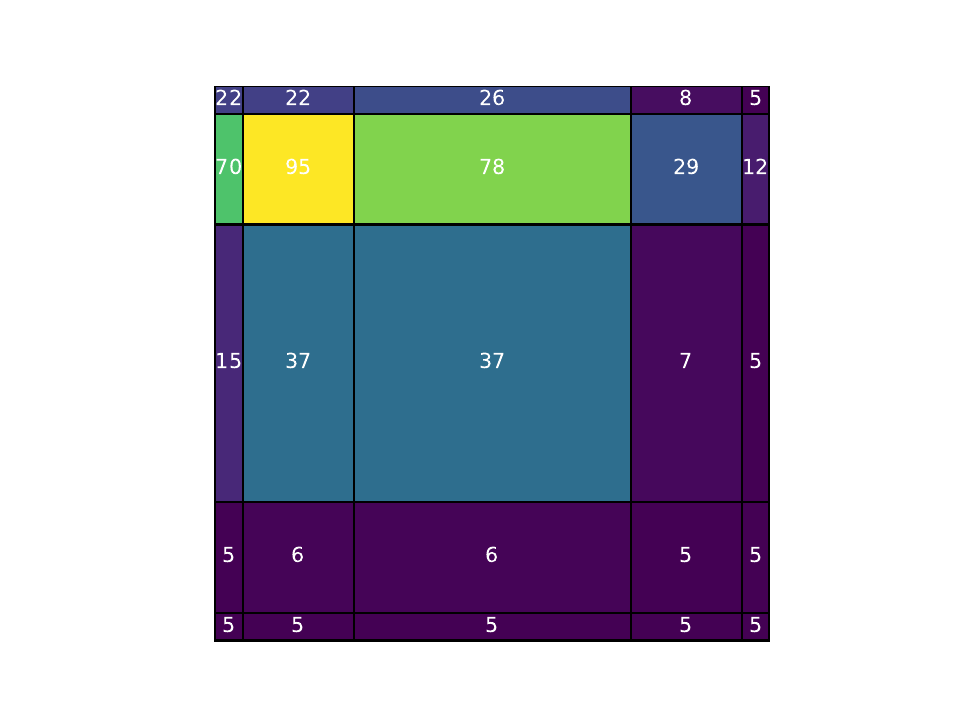}
        \caption{Intermediate ranks at time $t=1.0$}
    \end{subfigure}
    \begin{subfigure}{0.44\textwidth}
        \centering
        \includegraphics[width=\textwidth]{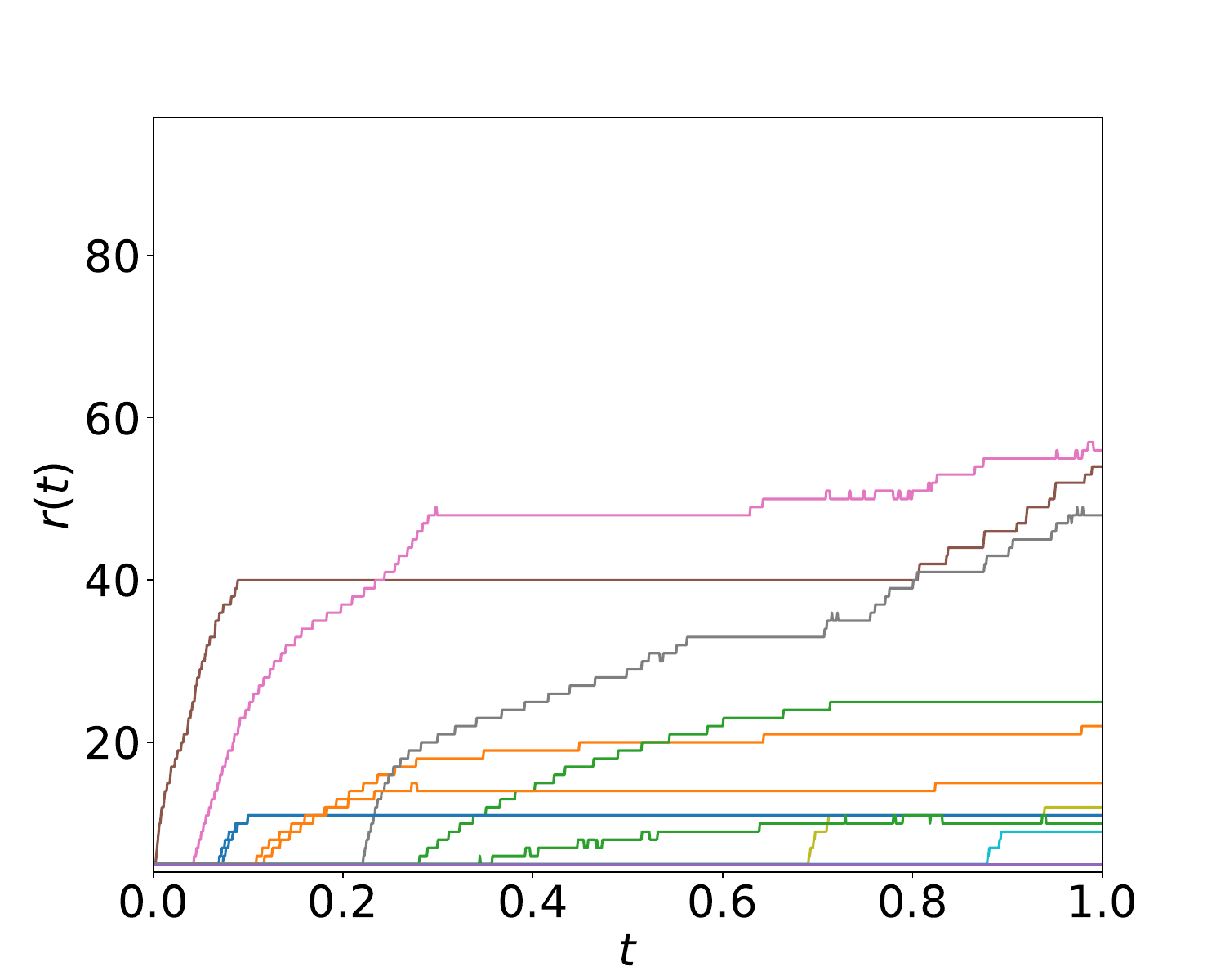}
        \caption{Ranks}
    \end{subfigure}
    \begin{subfigure}{0.46\textwidth}
        \centering
        \includegraphics[width=\textwidth]{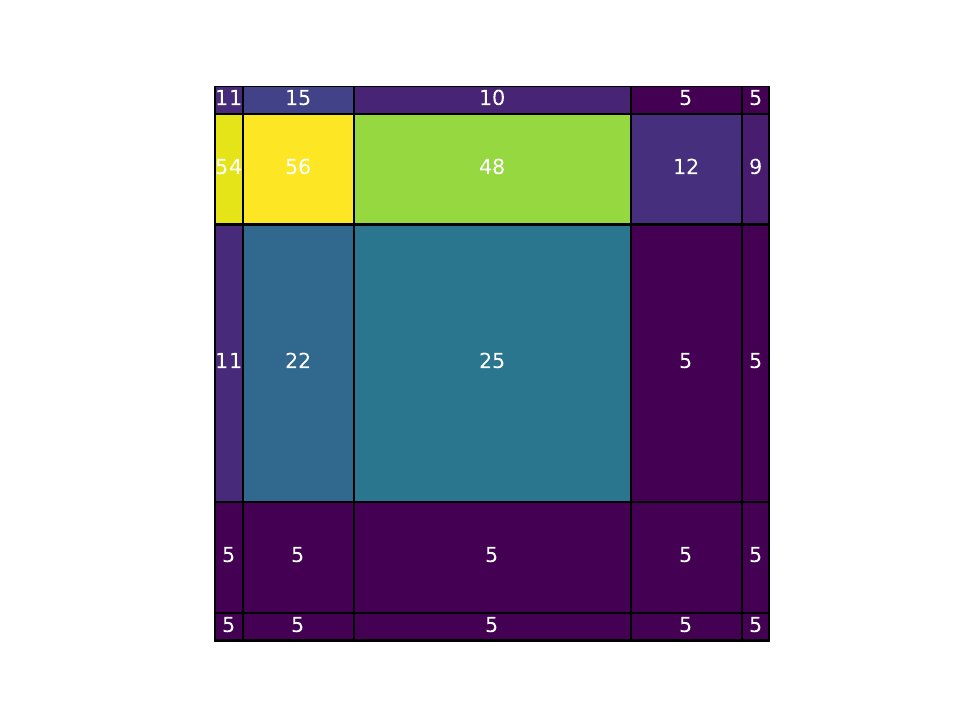}
        \caption{Ranks at time $t=1.0$}
    \end{subfigure}
    \caption{We show the time evolution of the intermediate ranks and ranks after the time step for the domain decomposition dynamical low-rank simulation of the point source test problem. The adaptive rank scheme as proposed in algorithms \ref{alg:augmentation} and \ref{alg:truncate} with {\color{black} $\text{tol}=1.2\cdot 10^{-4}$} was used.}    
    \label{fig::pointsource_domaindecomp_ranks}
\end{figure}

Here, we really want to focus on the ranks required by each method and compare the degrees of freedom. The ranks as a function of time for the two configurations are shown in figures \ref{fig::pointsource_classic_ranks} and \ref{fig::pointsource_domaindecomp_ranks}, respectively. While for classic simulations the rank needs to be high on the full domain, for domain decomposition simulations one can clearly see that even though the rank that is required close to the point source is high, as the radiation spreads into the domain a significantly lower rank can be used.

\begin{figure}[H]
    \centering
    \includegraphics[width=0.45\textwidth]{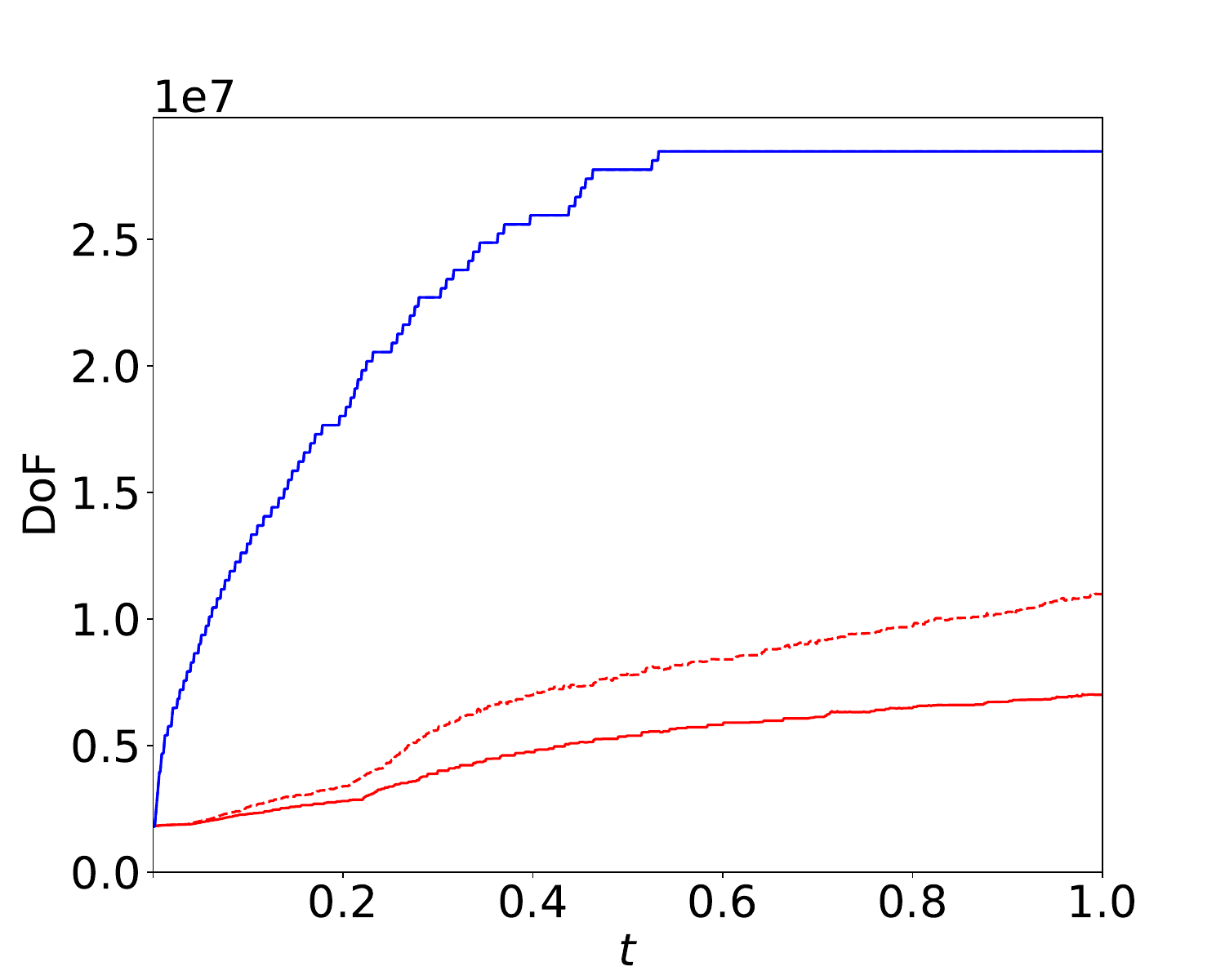}
    \includegraphics[width=0.45\textwidth]{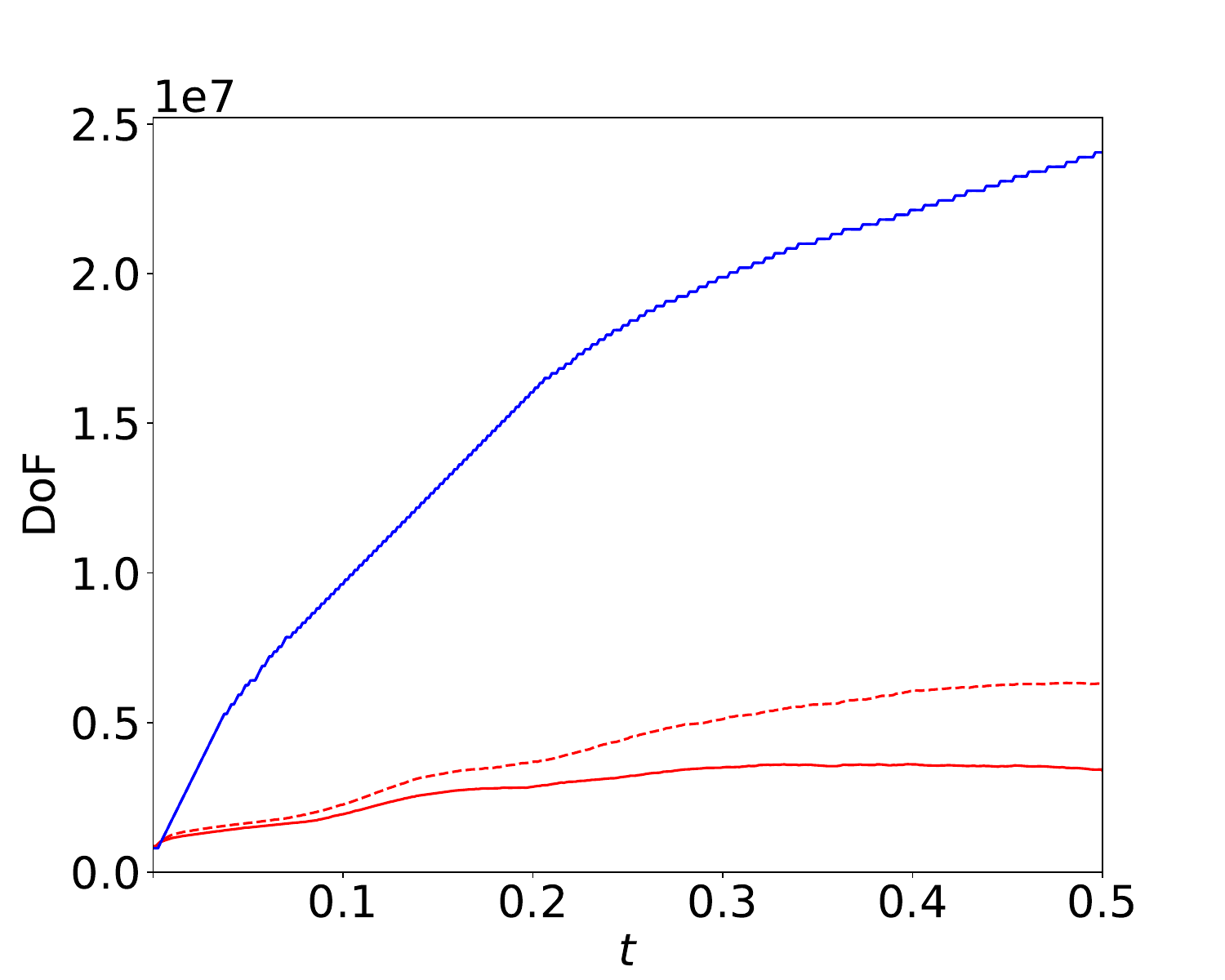}
    \caption{We show degrees of freedom (DoF) as a function of time for the low-rank simulation of the point source test problem (left) {\color{black} and the line source test problem (right)} for the classic low-rank method (blue) and the domain decomposition low-rank method (algorithm~\ref{alg::ddlr_algorithm}) (red).}
    \label{fig::pointsource_dof}
\end{figure}

Thus, the domain decomposition approach here has the clear advantage that the overall degrees of freedom can be kept significantly lower (roughly by a factor of $3$ for the intermediate rank and $5$ for the rank after each time step) compared to a classic low-rank approach, as can be seen in figure \ref{fig::pointsource_dof}.

{
\color{black}
\subsection{Line source test problem}\label{sec::linesource}
Finally, we consider a line source test problem. In this test case, our full domain is a purely scattering domain, and we do not have varying coefficients, i.e., $c_{\text{t}}=c_{\text{s}}=1$ on the whole domain. The initial condition is given by a line source in the middle of the domain, i.e.,
\begin{equation*}
  \begin{aligned}
  f(t=0,x,y,\phi) = \frac{1}{2\pi\cdot 0.001^2} \exp\left(-\frac{(x-0.5)^2+(y-0.5)^2}{2\cdot 0.05^2}\right), &\text{ for all } x\in [0,1], y\in [0,1] \\
  &\text{ and } \phi\in [0,2\pi].
  \end{aligned}
\end{equation*}
We note that in this case we have no inflow condition and no source, and thus again use the adaptive rank algorithm from \cite{Hochbruck_rankadaptivity} for the simulation on a single domain. The simulations are conducted on a grid with 400 grid points in $x$ and $y$ and 200 grid points in $\phi$. For domain decomposition we decompose the problem into $8\times 8$ subdomains. On the outer boundary, we impose outflow boundary conditions. The time step size is again set to $\Delta t = 0.5 \cdot \Delta x$. 

The results using the classic dynamical low-rank algorithm (left) and the domain decomposition low-rank algorithm \ref{alg::ddlr_algorithm} (right) for $\rho(t,x,y)$ at times $t=0.1$, $t=0.3$, and $t=0.5$ are shown in figure \ref{fig::linesource_rho_plots}. As can be clearly observed, the radiation spreads out from the line source into the domain. We again chose our tolerances, such that we get similar results for both simulations. Still, one can see differences in the results, the main difference being that our domain decomposition algorithm seems to capture the circular shape of the solution better.

\begin{figure}[H]
    \centering
    \begin{subfigure}{\textwidth}
        \centering
        \begin{minipage}{0.12\textwidth}
          \raggedleft
          \caption{$t=0.1$:}
        \end{minipage}
        \begin{minipage}{0.424\textwidth}
          \includegraphics[width=\textwidth]{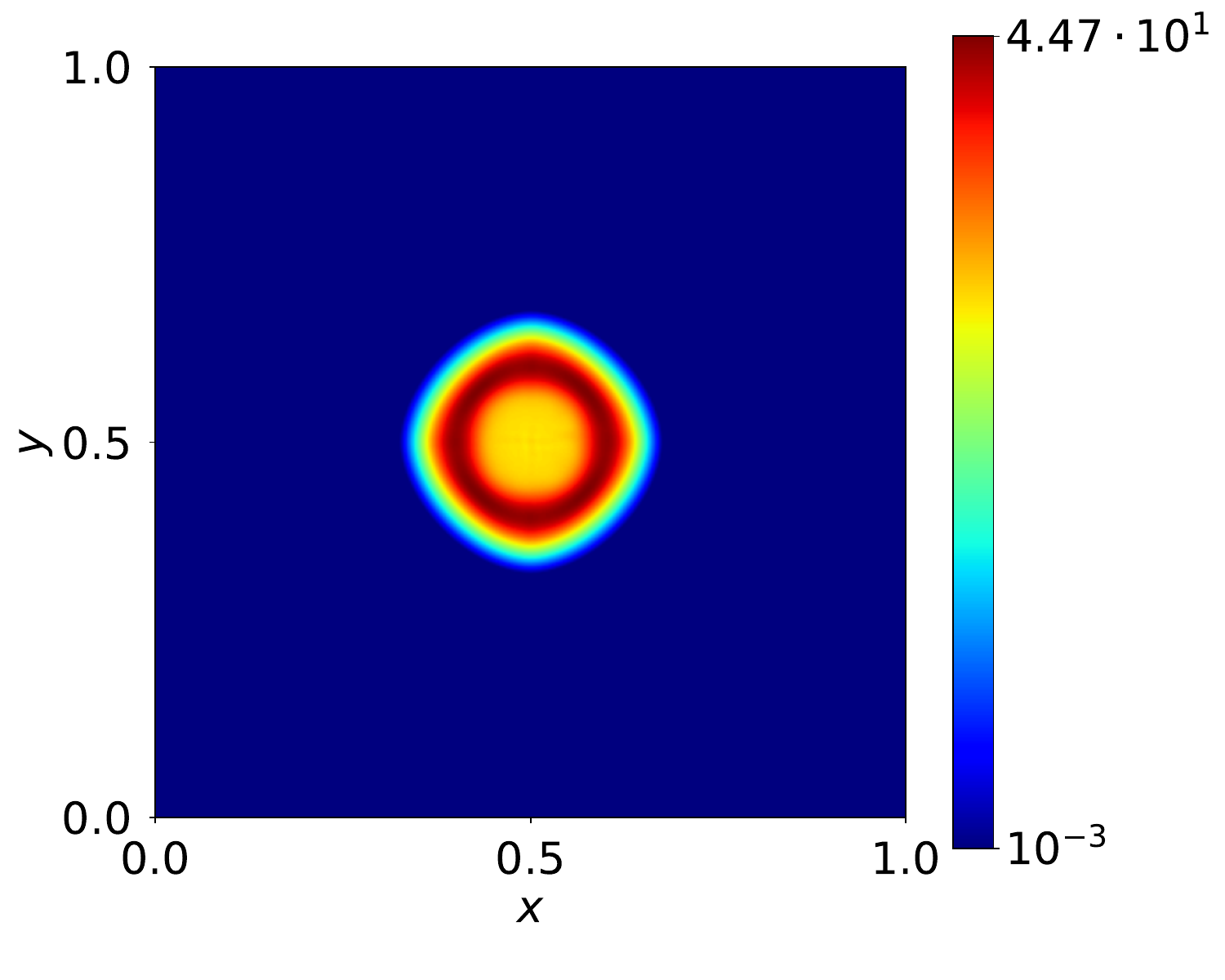}
        \end{minipage}
        \begin{minipage}{0.424\textwidth}
            \centering
            \includegraphics[width=\textwidth]{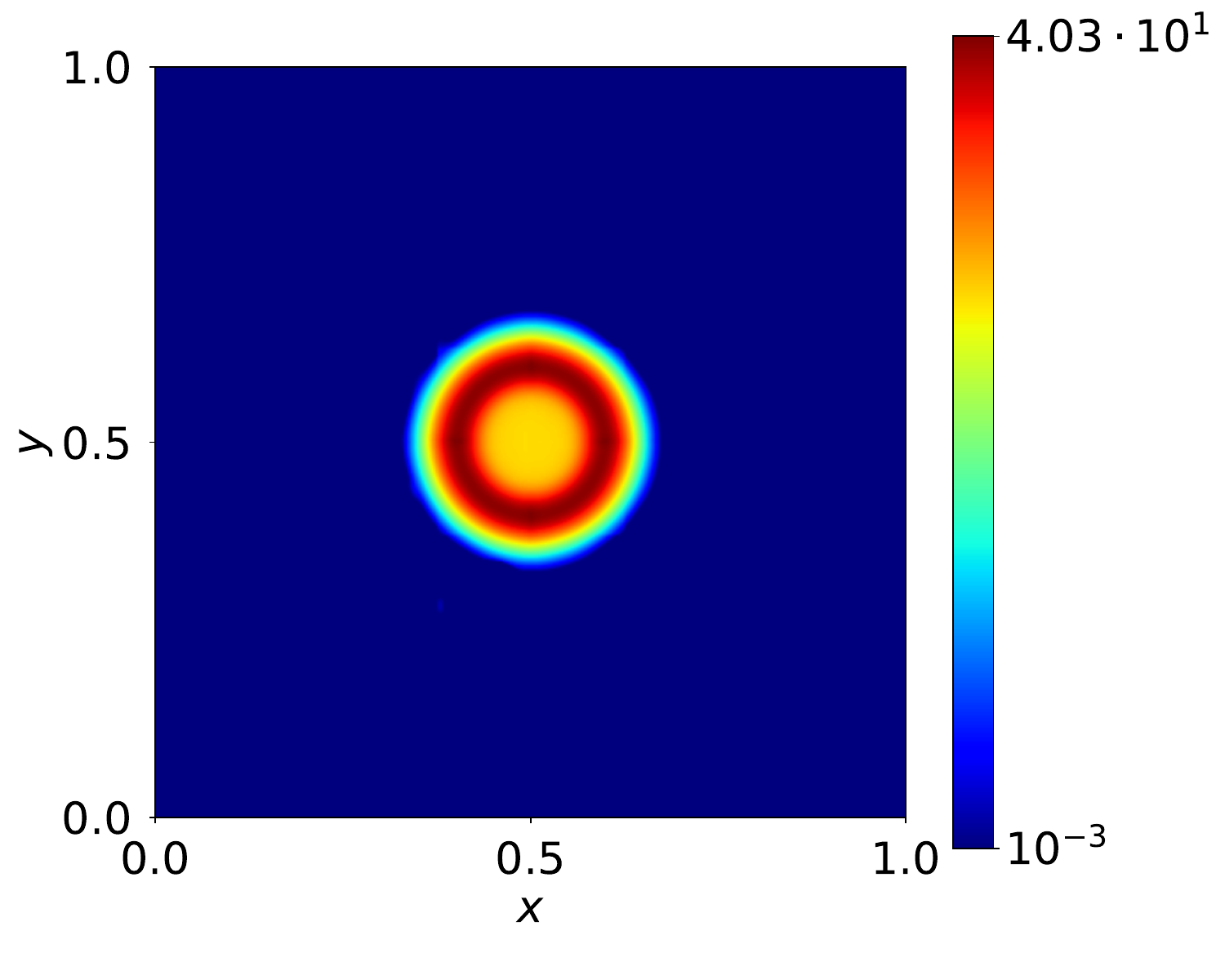}
        \end{minipage}
    \end{subfigure}
    
    \begin{subfigure}{\textwidth}
        \centering
        \begin{minipage}{0.12\textwidth}
          \raggedleft
          \caption{$t=0.3$:}
        \end{minipage}
        \begin{minipage}{0.424\textwidth}
          \includegraphics[width=\textwidth]{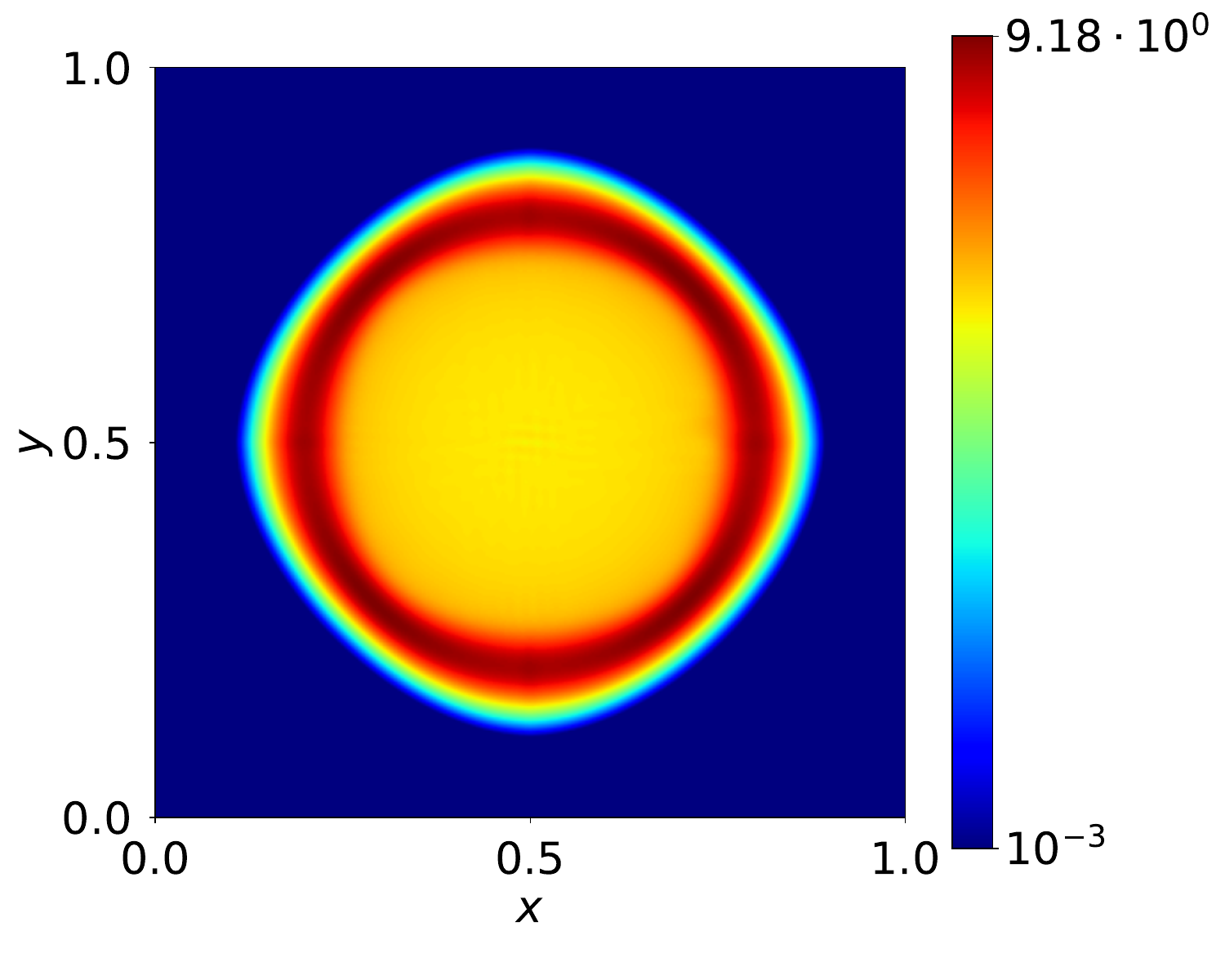}
        \end{minipage}
        \begin{minipage}{0.424\textwidth}
            \centering
            \includegraphics[width=\textwidth]{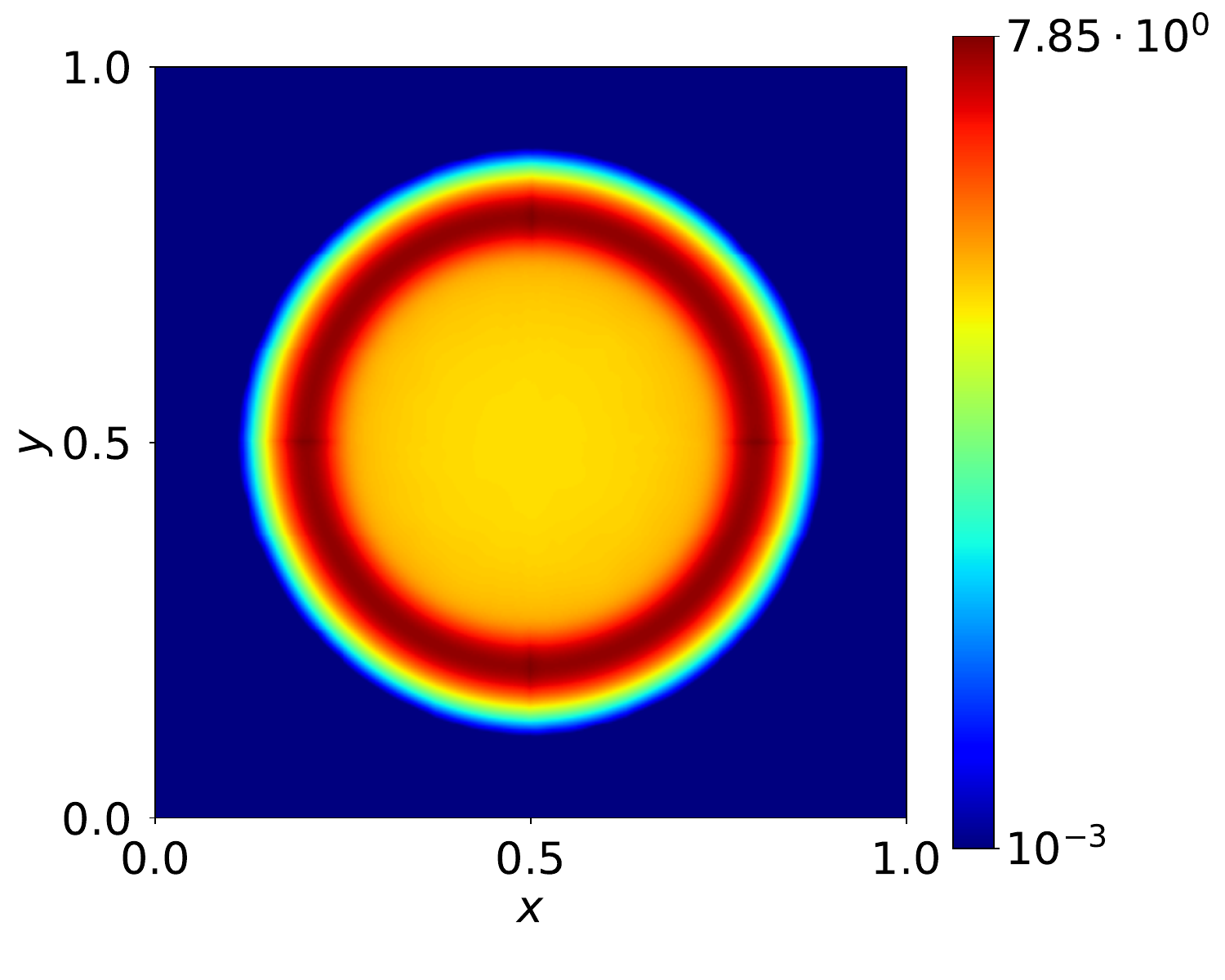}
        \end{minipage}
    \end{subfigure}
    
    \begin{subfigure}{\textwidth}
        \centering
        \begin{minipage}{0.12\textwidth}
          \raggedleft
          \caption{$t=0.5$:}
        \end{minipage}%
        \begin{minipage}{0.424\textwidth}
          \includegraphics[width=\textwidth]{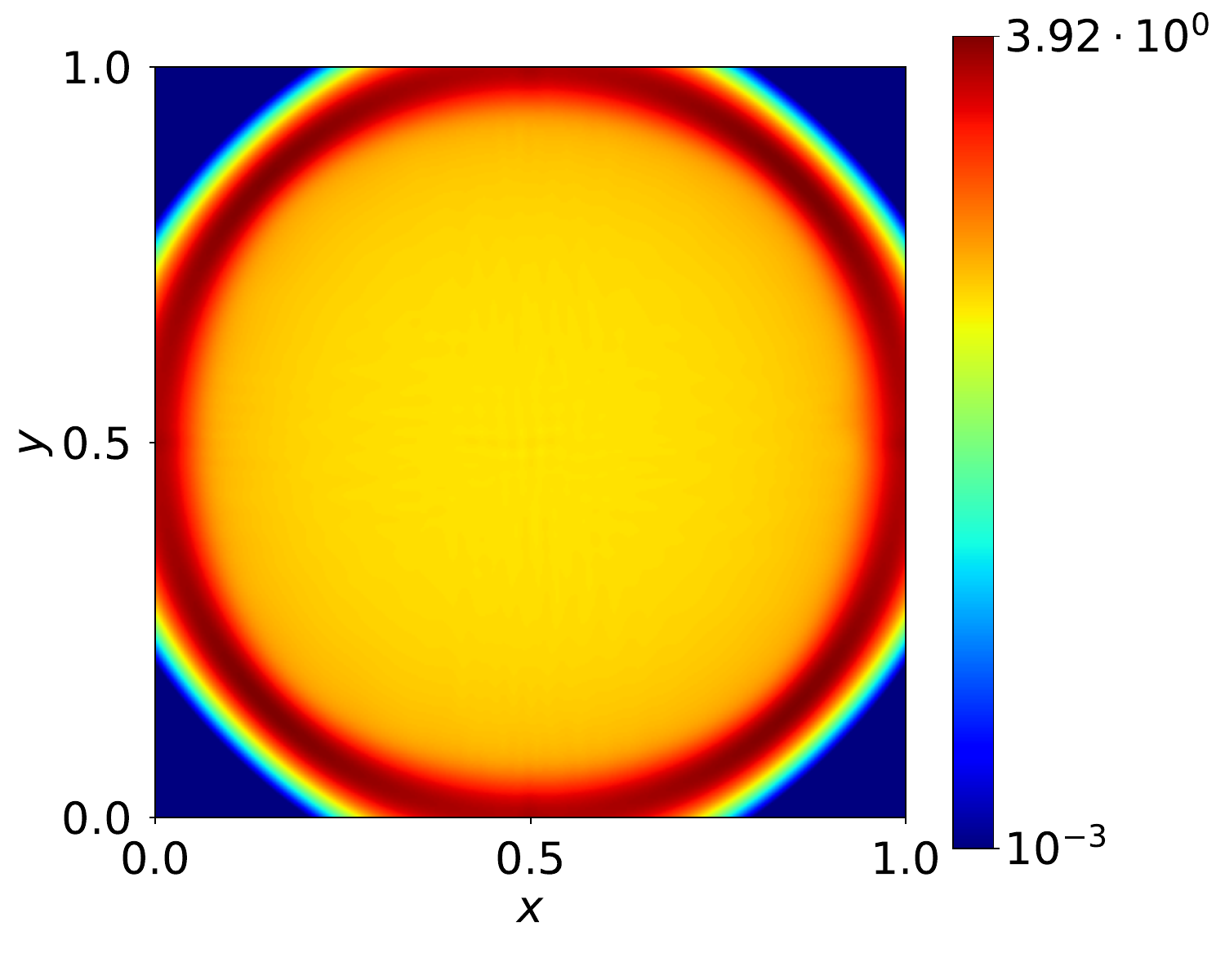}
        \end{minipage}
        \begin{minipage}{0.424\textwidth}
            \centering
            \includegraphics[width=\textwidth]{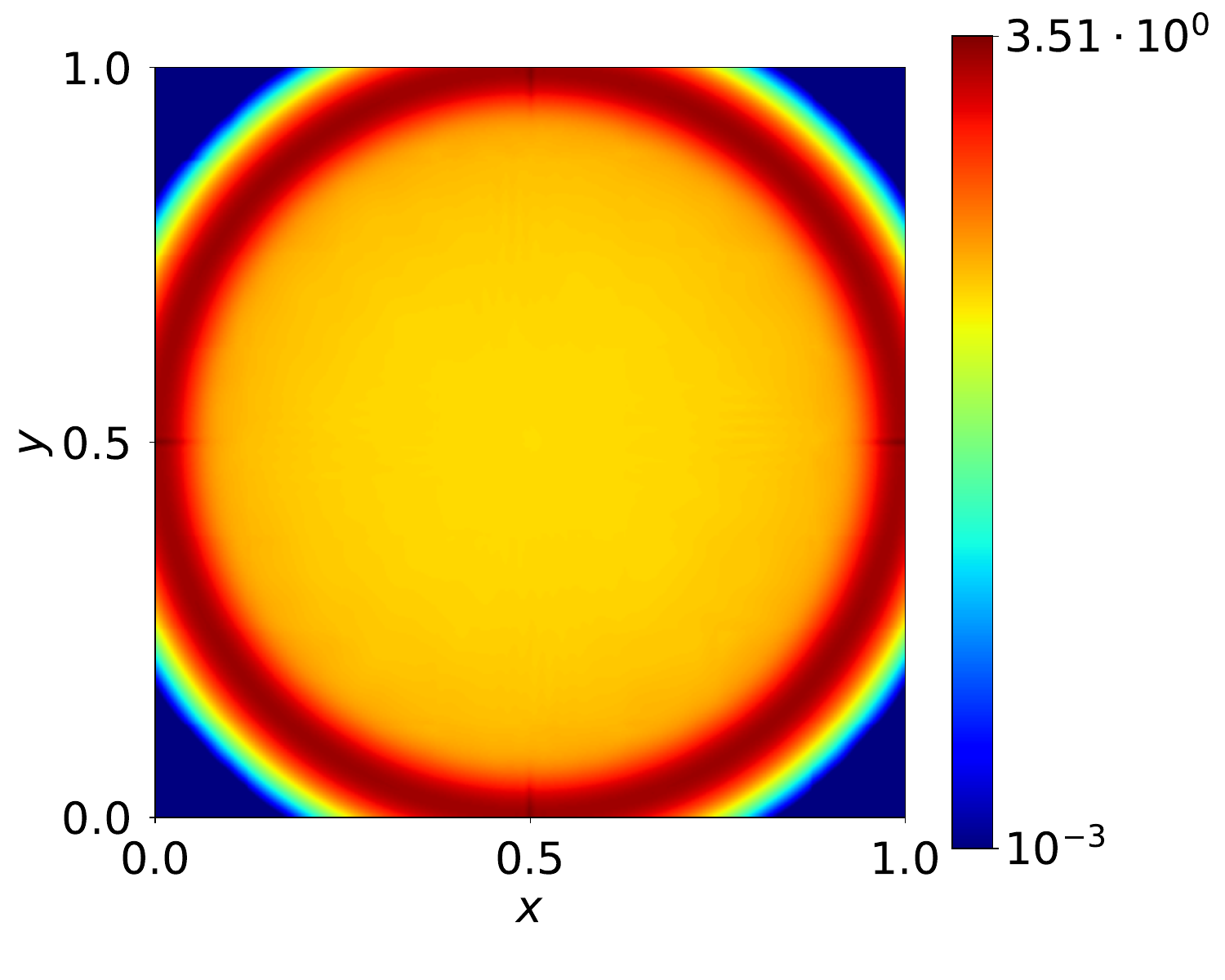}
        \end{minipage}
    \end{subfigure}
    \caption{We show $\rho (t,x,y)$ obtained by solving the line source test problem with the classic dynamical low-rank algorithm (left) and with the domain decomposition low-rank algorithm \ref{alg::ddlr_algorithm} (right). For the latter, the grid is divided into $8 \times 8$ subdomains.}
    \label{fig::linesource_rho_plots}
\end{figure}

When looking at the ranks in figures \ref{fig::linesource_1domain_ranks} and \ref{fig::linesource_domaindecomp_ranks}, we see that the classic low-rank method requires a high rank on the full domain, while for the domain decomposition low-rank method, only a few subdomains require a high rank. The adaptivity is also more natural in the domain decomposition approach, as the rank decreases again in the subdomains where the radiation has already passed through. This is also reflected in the degrees of freedom required by each method, as shown in figure \ref{fig::pointsource_dof}, where the domain decomposition approach requires roughly  a factor of 4 fewer degrees of freedom when comparing the intermediate ranks. 

\begin{figure}[H]
    \centering
    \includegraphics[width=0.44\textwidth]{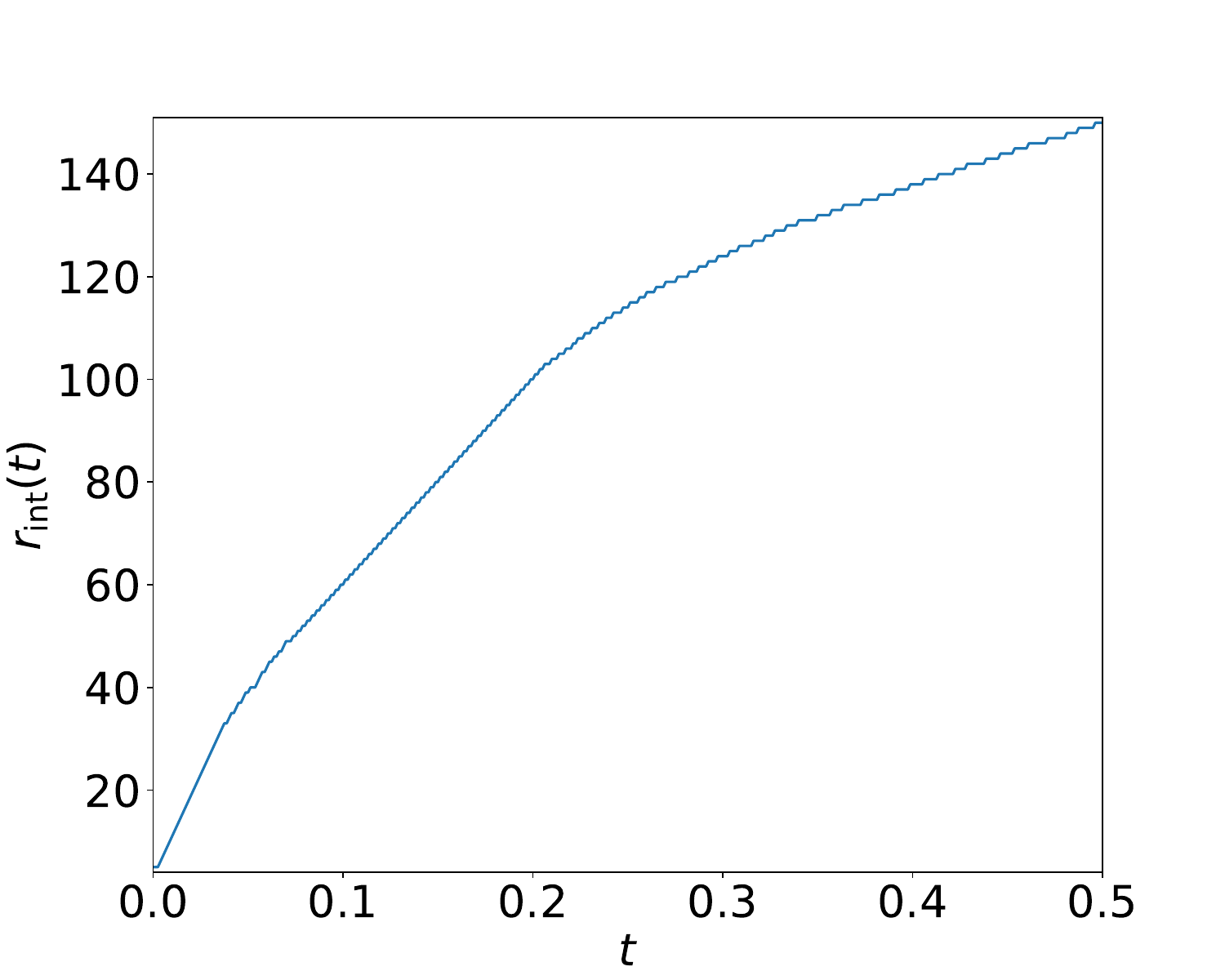}
    \caption{We show the time evolution of the rank for the classic low-rank simulation of the line source test problem.  As an adaptive rank scheme we used the scheme described in \cite{Hochbruck_rankadaptivity} with $\text{tol}=2\cdot10^{-3}$.}
    \label{fig::linesource_1domain_ranks}
\end{figure}

\begin{figure}[H]
    \centering
    \begin{subfigure}{0.44\textwidth}
        \centering
        \includegraphics[width=\textwidth]{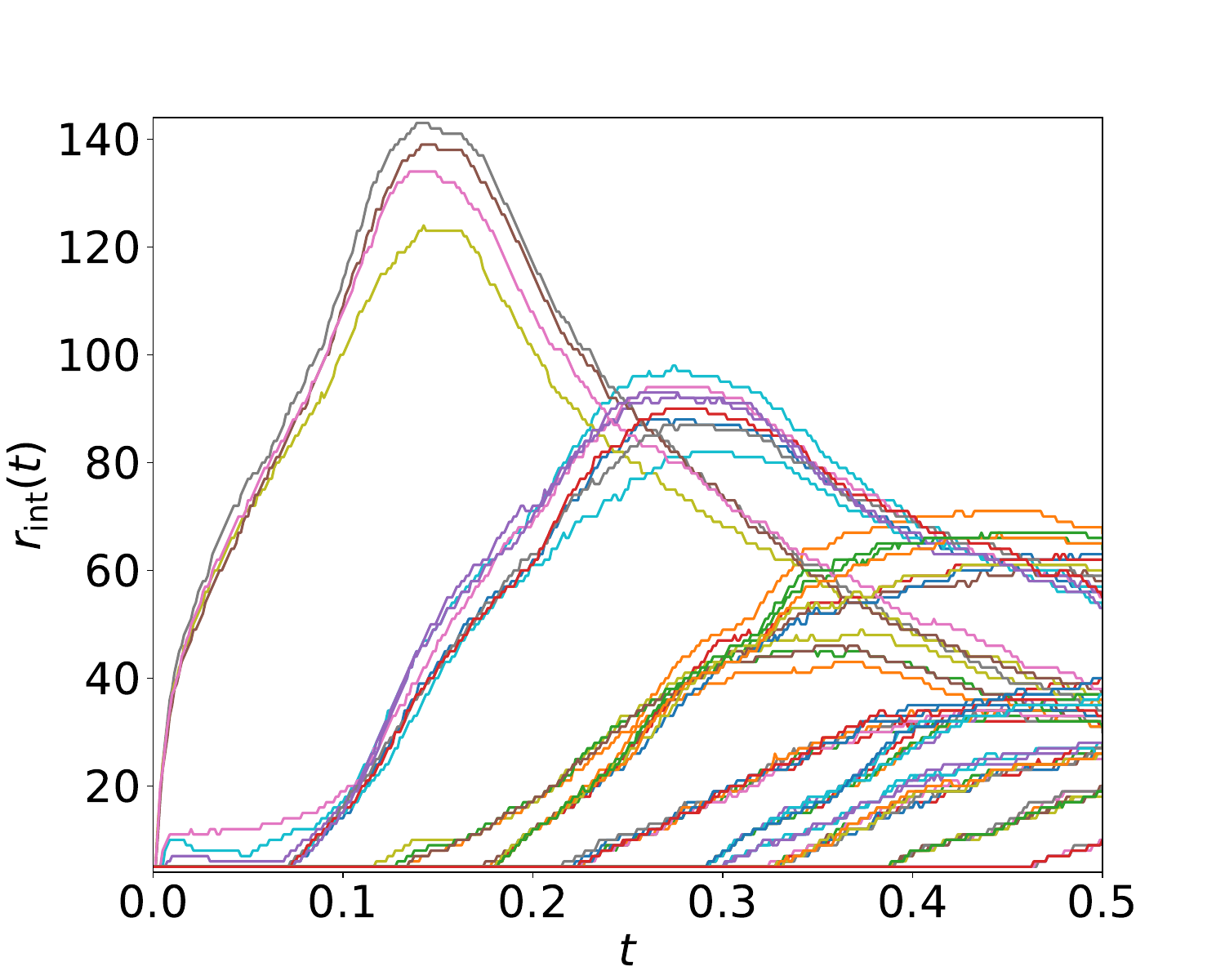}
        \caption{Intermediate ranks}
    \end{subfigure}
    \begin{subfigure}{0.46\textwidth}
        \centering
        \includegraphics[width=\textwidth]{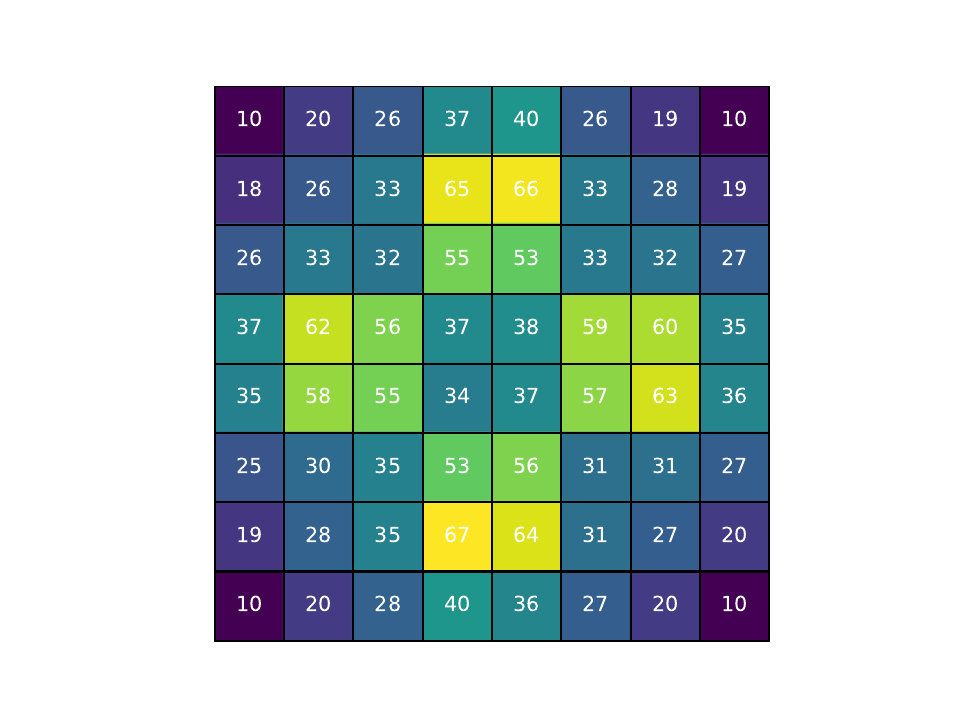}
        \caption{Intermediate ranks at time $t=0.5$}
    \end{subfigure}
    \begin{subfigure}{0.44\textwidth}
        \centering
        \includegraphics[width=\textwidth]{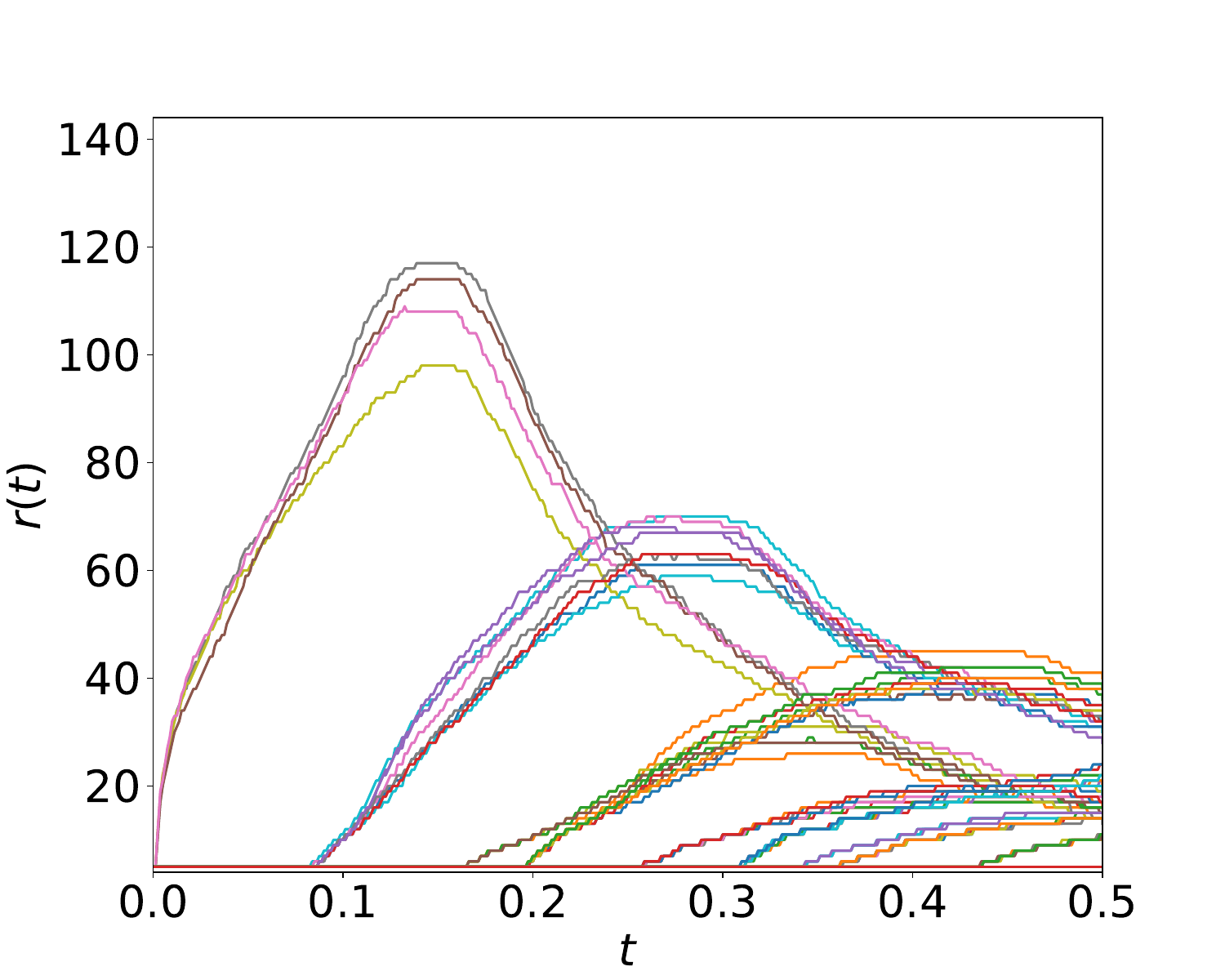}
        \caption{Ranks}
    \end{subfigure}
    \begin{subfigure}{0.46\textwidth}
        \centering
        \includegraphics[width=\textwidth]{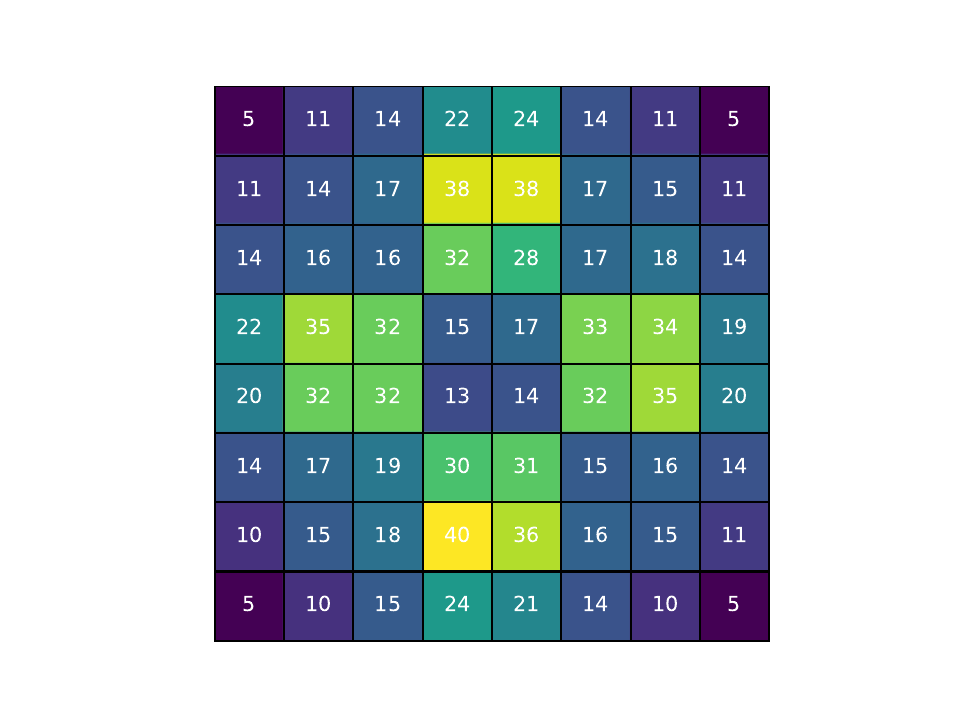}
        \caption{Ranks at time $t=0.5$}
    \end{subfigure}
    \caption{We show the time evolution of the intermediate ranks and ranks after the time step for the domain decomposition dynamical low-rank simulation of the line source test problem. The adaptive rank scheme as proposed in algorithms \ref{alg:augmentation} and \ref{alg:truncate} with $\text{tol}=2\cdot 10^{-4}$ was used.}    
    \label{fig::linesource_domaindecomp_ranks}
\end{figure}

}

{\color{black}
\section{Conclusion}
In this work we have proposed a domain decomposition approach for dynamical low-rank approximation of the radiative transfer equation. The main idea is to use a separate low-rank representation in each subdomain, and to connect the simulations on those subdomains by inflow boundary conditions. To allow the information transfer between subdomains and improve the scheme's efficiency, we have proposed an adaptive rank scheme. We have shown that our method can be applied to challenging test problems, where the total degrees of freedom can be kept significantly lower compared to a classic low-rank approach, while still maintaining good accuracy and allowing parallelization on distributed memory systems. In future work, one could improve load-balancing by using subdomains of smaller size in regions where the rank is expected to be higher, or by grouping subdomains together. Additionally, extending our scheme to an AP scheme might be interesting future work.
}

\section{Acknowledgements}

This work was performed under the auspices of the U.S. Department of Energy by Lawrence Livermore National Laboratory under Contract DE-AC52-07NA27344.

\raggedright \printbibliography

\appendix
\section{Application of upwind discretization to algorithm \ref{alg::ddlr_algorithm}}\label{sec::appendix_upwind}

Our goal in this section is to apply the first order upwind scheme to the K step in both algorithms \ref{alg:projec-advx} and \ref{alg:projec-advy}. 
We will only show the derivation for the $x$ advection step, since the $y$ advection step can be done analogously.

To begin, we will rewrite equation \eqref{eq::advx_evolK} in a discretized setting using matrices in linear algebra notation. We get
\begin{equation}\label{eq::Kstep_appendixA}
  \partial_t K = - c_{\text{adv}} D_X K C,
\end{equation}
where $C = V^T \text{diag}(\sin(\phi)) V \cdot \Delta \phi$ and $D_X$ is the discretization of $\partial_x$.
In order to obtain the direction of flow we diagonalize the symmetric and small matrix $C$. This gives
\begin{equation*}
  C = P \Lambda P^{-1},
\end{equation*}
where $\Lambda$ is a diagonal matrix with eigenvalues on its diagonal and $P$ contains the corresponding eigenvectors. Plugging this into equation \eqref{eq::Kstep_appendixA} we get
\begin{equation*}
  \partial_t K = - c_{\text{adv}} D_X K P \Lambda P^{-1}.
\end{equation*}
Next, we express $D_XKP$ by applying the upwind matrices $D_X^-$ and $D_X^+$ to each column of $KP$, depending on the sign of the corresponding eigenvalue in $\Lambda$. More precisely, if $\Lambda[i,i]>0$, then
\begin{equation*}
  D_XKP[:,i] = D_X^- KP[:,i],
\end{equation*}
and if $\Lambda[i,i]<0$, then
\begin{equation*}
  D_XKP[:,i] = D_X^+ KP[:,i].
\end{equation*}

\end{document}